\documentclass{amsart}

\usepackage{amssymb}
\usepackage{amsfonts}
\usepackage{amscd}
\usepackage{verbatim}
\usepackage{epsfig}

\input xy
\xyoption{all}

\DeclareFontFamily{OT1}{rsfs}{}
\DeclareFontShape{OT1}{rsfs}{n}{it}{<-> rsfs10}{}
\DeclareMathAlphabet{\mathscr}{OT1}{rsfs}{n}{it}

\newtheorem{thm}{Theorem}[subsection]
\newtheorem{lem}[thm]{Lemma}
\newtheorem{prop}[thm]{Proposition}
\newtheorem{cor}[thm]{Corollary}
\newtheorem{rem}[thm]{Remark}
\newtheorem{rems}[thm]{Remarks}

\newtheorem{sth}{Theorem}[thm]

\theoremstyle{definition}

  \newtheorem{defi}[thm]{Definition}
  
  \newtheorem{srem}[sth]{Remark}
  
  \newtheorem{exam}[thm]{Example}

  \newtheorem{ack}{Acknowledgements} 

\numberwithin{equation}{thm}

\newcommand{\Cref}[1]{Corollary~\textup{\ref{#1}}}
\newcommand{\Dref}[1]{Definition~\textup{\ref{#1}}}
\newcommand{\Eref}[1]{Example~\textup{\ref{#1}}}
\newcommand{\Lref}[1]{Lemma~\textup{\ref{#1}}}
\newcommand{\Pref}[1]{Proposition~\textup{\ref{#1}}}
\newcommand{\Rref}[1]{Remark~\textup{\ref{#1}}}
\newcommand{\Rsref}[1]{Remarks~\textup{\ref{#1}}}
\newcommand{\Sref}[1]{Section~\textup{\ref{#1}}}
\newcommand{\Ssref}[1]{Subsection~\textup{\ref{#1}}}
\newcommand{\Tref}[1]{Theorem~\textup{\ref{#1}}}

\def\bilap#1{\hbox to 0pt{\hss#1\hss}}
  \def\Rarrow#1{\bilap{\hbox to#1{\rightarrowfill}}}
  \def\Larrow#1{\bilap{\hbox to#1{\leftarrowfill}}}
  \def\Equals#1{\bilap
                   {\raise 4pt\hbox
                     {\vrule width#1 height.5pt}%
                    \kern-#1\raise 1pt\hbox
                     {\vrule width#1 height.5pt}%
                   }}

\newcommand{\EQAL}[1]%
{\,\begin{picture}(#1,0)%
\put(0,3){\line(1,0){#1}}%
\put(0,1){\line(1,0){#1}}%
\end{picture}\,}%

\newcommand{\vlto}[1]%
{\,\begin{picture}(#1,3)%
\put(0,2){\vector(1,0){#1}}%
\end{picture}\,}%

\newcommand{\vllarrow}[1]%
{\,\begin{picture}(#1,3)%
\put(#1,2){\vector(-1,0){#1}}%
\end{picture}\,}%

\newcommand{\dirlm}[1]%
   {
      {\lim\hskip-1.58em\lower.65ex
        \hbox{$
                 {}_{\stackrel{\lower1ex\hbox
                                         {$\scriptstyle -\!\!\!\longrightarrow$}
                                       }{\vbox to0pt{\vss\vskip.6ex
                                             \hbox{$\scriptstyle{}^{#1}$}\vss}}
                    }
             $}
      }
\:}

\newcommand{\subdirlm}[1]%
   {
      {\lim\hskip-1.5em\lower.6ex
        \hbox{$
                    {}_{\stackrel{\lower1ex\hbox
                                            {$\scriptstyle\longrightarrow$}
                                 }{ ^{#1} }
                       }
              $}
      }
\:}

\newcommand{\inlm}[1]%
    {
       {\lim\hskip-1.58em\lower.65ex
         \hbox{$
                  {}_{\stackrel{\lower1ex\hbox
                                         {$\scriptstyle \longleftarrow\!\!\<-$}
                               }{\vbox to0pt{\vss\vskip.6ex
                                             \hbox{$\scriptstyle{}^{#1}$}\vss}}
                     }
              $}
       }
\:}

\def\hz#1{{\hbox to 0pt{#1}}}
\def\Iso{\vbox to 0pt{\vss\hbox{$\widetilde{\phantom{nn}}$}\vskip-7pt}}

\def\>{\mspace {1mu}}
\def\<{\mspace{-1mu}}
\def\({{\textup(}}
\def\){{\textup)}}

\newcommand{\fm}{{\mathfrak{m}}}

\newcommand{\fp}{{\mathfrak{p}}}
\newcommand{\fq}{{\mathfrak{q}}}
\newcommand{\X}{{\mathscr X}}

\newcommand{\Y}{{\mathscr Y}}

\newcommand{\W}{{\mathscr W}}
\newcommand{\I}{{\mathscr I}}
\newcommand{\J}{{\mathscr J}}
\newcommand{\eA}{{\mathscr A}}
\newcommand{\eB}{{\mathscr B}}
\newcommand{\eC}{{\mathscr C}}
\newcommand{\eD}{{\mathscr D}}
\newcommand{\eE}{{\mathscr E}}
\newcommand{\eH}{{\mathscr H}}

\newcommand{\eS}{{\mathscr S}}
\newcommand{\eM}{{\mathscr M}}
\newcommand{\eN}{{\mathscr N}}
\newcommand{\co}{{\mathscr O}}
\newcommand{\eF}{{\mathscr F}}
\newcommand{\eG}{{\mathscr G}}
\newcommand{\ssf}{{\mathsf{s}}}
\newcommand{\Ssf}{{\mathsf{S}}}
\newcommand{\Spec}{{\mathrm {Spec}}}
\newcommand{\Spf}{{\mathrm {Spf}}}
\newcommand{\cm}{{\mathbf{cm}}}

\newcommand{\A}{{\mathcal A}}

\newcommand{\cD}{{\mathcal D}}
\newcommand{\De}{{\Delta}}

\newcommand{\cDt}{\cD_{\<\mathrm t}}

\newcommand{\F}{{\mathcal F}}
\newcommand{\G}{{\mathcal G}}

\newcommand{\Hr}{{\mathrm H}}

\newcommand{\M}{{\mathcal M}}

\newcommand{\eR}{{\mathscr R}}
\newcommand{\eK}{{\mathscr K}}
\newcommand{\cR}{{\mathcal R}}

\newcommand{\D}{{\mathbf D}}
\newcommand{\K}{{\mathbf K}}
\newcommand{\C}{{\mathbf C}}

\newcommand{\bbF}{{\mathbb F}}
\newcommand{\bbFc}{{\bbF_c}}

\newcommand{\rbbF}{{{\bbF}^{\>r}}}

\newcommand{\rbbFc}{{\bbF^{\>r}_{\<\<\rm c}}}
\newcommand{\Coz}{{\mathrm {Coz}}}
\newcommand{\Cozst}{{\mathbf{coz}}}
\newcommand{\Cozs}[1]{{\mathrm {Coz}}_{#1}}
\newcommand{\Cozt}[1]{{{\mathbf{coz}}}_{#1}}
\newcommand{\Ed}[1]{E_{#1}}

\newcommand{\vc}{{\vec{\mathrm{c}}}}

\newcommand{\wDqc}{ \widetilde
          {\vbox to6.5pt{\vss\hbox{$\mathbf D$}}}
    _{\mkern-1.5mu\mathrm {qc}} }

\newcommand{\wDqcp}{\wDqc^{\lower.5ex\hbox{$\scriptstyle+$}}}

\newcommand{\Dvc}{\D_{\<\vc}}
\newcommand{\Dqct}{\D_{\mkern-1.5mu\mathrm{qct}}}

\newcommand{\Dc}{\D_{\mkern-1.5mu\mathrm c}}
\newcommand{\Dcs}{{\Dc\<\!\!{}^*}}

\newcommand{\qc}{{\mathrm{qc}}}

\newcommand{\R}{{\mathbf R}}

\newcommand{\bL}{{\mathbf L}}
\newcommand{\Hom}{{\mathrm {Hom}}}

\newcommand{\Homb}{{\mathrm {Hom}}^{\bullet}}
\newcommand{\Ac}{\A_{\mathrm c}}
\newcommand{\Aqc}{\A_{\qc}}
\newcommand{\Avc}{\A_{\vec {\mathrm c}}}
\newcommand{\Aqct}{\A_{\mathrm {qct}}\<}

\newcommand{\At}{\A_{\mathrm t}\<}
\newcommand{\Dt}{\D_{\mathrm t}\<}

\newcommand{\fs}{f^!}

\newcommand{\ga}[1]{\gamma^!_{\hbox{$\scriptstyle{#1}$}}}

\newcommand{\sh}[1]{{#1}^{\sharp}}
\newcommand{\psh}[1]{{#1}^{(\sharp)}}

\newcommand{\pshr}[1]{{#1}^{(!)}}

\newcommand{\Tr}[1]{{\mathrm {Tr}}_{{#1}}}

\newcommand{\sTr}[1]{{\mathrm {Tr}}^{(\sharp)}_{{#1}}}
\newcommand{\ttr}[1]{{\mathrm {\tau}}_{#1}}
\newcommand{\Ext}{\eE{xt}}
\newcommand{\sHom}{\eH{om}}
\newcommand{\sHomb}{\eH{om}^{\bullet}}
\newcommand{\iGp}[1]{{\varGamma_{\<\!#1}'}}
\newcommand{\iG}[1]{{\varGamma_{\<\!#1}^{\phantom\prime}}}

\newcommand{\set}{\!:=}

\newcommand{\iso}%
{{\mkern8mu\longrightarrow \mkern-25.5mu{}^\sim\mkern17mu}}
\newcommand{\osi}%
{{\mkern8mu\longleftarrow \mkern-24.5mu{}^\sim\mkern16mu}}
\newcommand{\Otimes}{\underset
   {\vbox to 0pt {\vskip-1ex\hbox{$\scriptscriptstyle=$}\vss}}
     {\otimes}\vadjust{\kern.4pt}} 
\newcommand{\BG}{\boldsymbol{\varGamma}}

\newcommand{\BL}{{\boldsymbol\Lambda}}

\newcommand{\smcirc}%
   {{\raise.15ex\hbox to.7em{$\hss \scriptstyle\circ\hss$}}}

\title[Duality theory and Cousin complexes]%
{Applications of duality theory to Cousin complexes}

\author[S.\,Nayak]{Suresh Nayak}%
{\smash{}}
\address{Chennai Mathematical Institute\\
                 Sipcot IT Park\\
                 Siruseri, TN-603103, INDIA}
\email{snayak@cmi.ac.in}

\author[P.\,Sastry]{Pramathanath Sastry}%
{\smash{}}
\address{Department of Mathematics\\
                 East Carolina University,
                 Greenville, NC 27858, USA }
\email{sastryp@ecu.edu}
\date{\today}

\begin{document}
%{\bf\hfill \today}

\begin{abstract}{We use the anti-equivalence between Cohen-Macaulay complexes and
coherent sheaves on formal schemes to shed light on some older results and prove new
results. We bring out the relations between a coherent
sheaf $\eM$ satisfying an $S_2$ condition and the lowest cohomology
$\eN$ of its ``dual" complex. We show that if a scheme has a Gorenstein complex satisfying certain
coherence conditions, then in a finite \'etale neighborhood of each point, it has a dualizing
complex. If the scheme already has a dualizing complex, then we show that the
Gorenstein complex must be a tensor product of a dualizing complex and a vector bundle of
finite rank. We relate the various results in \cite{dcc} on Cousin complexes to dual results
on coherent sheaves on formal schemes.}
\end{abstract}

\maketitle

%\setcounter{tocdepth}{2}
%\tableofcontents
%\newpage
%\input sec1.tex

%\newpage

\section{Introduction}

The theme which lurks behind the various results in this paper
is the (anti) equivalence between Cohen-Macaulay complexes and
coherent sheaves proven in
\cite[p.\,108,\,Prop.\,9.3.1 and Cor.\,9.3.2]{lns} and restated here
in \Pref{prop:Ac-CM} and \Pref{prop:Ac-Coz}. The Cohen-Macaulay complexes
we just referred to are with respect to a fixed codimension function and
satisfy certain coherence conditions, which for ordinary schemes amount
to requiring that all cohomology sheaves are coherent. In all the sections, except
\Sref{s:gorenstein}, the formal schemes
involved are also required to satisfy conditions---e.g.~they should carry
``c-dualizing complexes" (see \Dref{def:dualizing} below).

This anti-equivalence is the unifying thread that runs through the
three main topics of this paper.  It was first observed by Yekutieli and
Zhang in \cite[Thm.\,8.9]{yz} for ordinary schemes 
of finite type over a regular scheme, and later in greater generality
by Lipman and the authors in \cite{lns}.  
We first give a short 
description of each topic we deal with before embarking on a more
detailed discussion putting our results in context. Here is the
brief version:

1) We explore symmetries between a coherent sheaf (on an ordinary scheme)
satisfying an ``$S_2$ condition" with respect to a codimension function
(cf.~\Dref{def:s2}) and an associated ``dual" coherent sheaf (which also
is shown to satisfy the same $S_2$ condition). The example to keep in mind
is the symmetry between the structure sheaf of an $S_2$ scheme and a
canonical module on the scheme (cf.~\cite[p.\,19,\,Thm.\,1.4]{dt} and
\cite[Thm.\,4.4]{kw}).

2) We give a relationship between Gorenstein complexes and
dualizing complexes (both with respect to a fixed codimension function).

3) We find an alternate approach to some of the results in \cite{dcc} when
our Cousin complexes involved satisfy certain coherence conditions (which,
as before, translate on an ordinary scheme to usual coherence conditions).
\pagebreak 
And in this approach we do not need to assume that the maps involved
(between formal schemes) are composites of compactifiable maps. It was
A.~Yekutieli who made the suggestion (to the second
author) that the results in \cite{dcc} should be re-examined
in light of the above mentioned duality between 
Cohen-Macaulay complexes and coherent sheaves. It should, however, be
pointed out that the proofs in the present paper only give an illusion of being
simpler for  we need the deeper results of \cite{dcc}, which deals with a larger
category. However, our proofs are illuminating, since they interpret operations on
Cousin complexes in terms of natural operations on the dual category of coherent
sheaves.

%\newpage

Let us examine each of these topics in somewhat greater detail. All
schemes involved (formal or ordinary) are assumed to be noetherian
and carrying a c-dualizing complex (forcing them to be of finite
Krull dimension). We use the following notations
\begin{itemize}
\item For a scheme $\X$, and $\co_{\X}$-modules $\eA$ and $\eB$, we often
write $\Hom(\eA,\,\eB)$ and $\sHom(\eA,\,\eB)$ for $\Hom_{\co_\X}(\eA,\,\eB)$
and $\sHom_{\co_\X}(\eA,\,\eB)$ respectively  (here $\sHom$ is the sheafified
version of $\Hom$).
\item Similarly, when no confusion is likely to arise, we write $\eA\otimes\eB$ for
$\eA\otimes_{\co_\X}\eB$.
\item For an $\co_\X$ complex $\eF$, and an integer $p$, $\eF^p$ will denote the
$p$-th graded piece of $\eF$.  
\item For two $\co_\X$ complexes $\eF$ and $\eG$, $\Homb(\eF,\,\eG)$ denotes
$\Homb_{\co_\X}(\eF,\,\eG)$,  i.e., $\Homb$ gives a complex of $\Gamma(\X,\,\co_\X)$-modules. 
In contrast, $\Hom(\eF,\,\eG)$ denotes the ``external $\Hom$", i.e. 
$\Hom(\eF,\,\eG)=\Hom_{\C(\X)}(\eF,\,\eG)$ the $\Gamma(\X,\,\co_\X)$-module of $\co_\X$ co-chain maps from $\eF$ to $\eG$. It is well known that the latter is the
module of $0$--cocycles of the former.
\item With $\eF$ and $\eG$ as above, $\sHomb_{\co_\X}(\eF,\,\eG)$ and 
$\sHom_{\C(\X)}(\eF,\,\eG)$ denote
the sheafified versions of 
$\Homb_{\co_\X}(\eF,\,\eG)$ and $\Hom_{\C(\X)}(\eF,\,\eG)$ respectively. As before, when the
context is clear, we write $\sHomb$ for $\sHomb_{\co_\X}$ and $\sHom$ 
for $\sHom_{\C(\X)}$. The relationship between $\sHomb$ and $\sHom$ is analogous to the relation between
$\Homb$ and $\Hom$. We identify the $\co_\X$-module $\sHom(\eF,\,\eG)$ 
with the sheaf of $0$-cocycles in $\sHomb(\eF,\,\eG)$.
\end{itemize} 

\subsection{$\De$-$S_2$ complexes} Let $X$ be an ordinary scheme and
let $\eR$ be a dualizing complex on $X$ which we assume (without loss
of generality) is residual. Let $\De\colon |X|\to {\mathbb{Z}}$ be
the associated codimension function (so that $\eR=\Ed{\De}\eR$, where
$\Ed{\De}$ is the Cousin functor associated with 
$\De$ (see \S\S\,2.3 below)). Recall that 
 $X$ is $S_2$ if and only if the natural map $\co_X\to \Ed{\eH}(\co_X)$ gives an
 isomorphism on applying $H^0$, where, $\Ed{\eH}(\co_X)$ is the Cousin complex
 of $\co_X$ with respect to the ``height filtration" $\eH=(H_i)$ given by
  $H_i=\{x\in X\,\vert\, \dim{\co_{X,x}}\ge i\}$. 
 One defines the notion of a $\De$-$S_2$ module
along the above lines (cf.~\Dref{def:s2}). Let $\eM$ be such a module,
which by definition is coherent. Let $\eN=\sHom(\Ed{\De}\eM,\eR)$.
We show that $\eN$ is also coherent and $\eM$ and $\eN$ share the
following symmetries, where ``=" denotes functorial isomorphisms
(cf.~\Tref{thm:main1}).
\begin{enumerate}
\item[(i)] $\eM = \sHom(\Ed{\De}\eN,\eR)$. 
(Note: $\eN\set\sHom(\Ed{\De}\eM,\eR)$.)
\item[(ii)] $\Ed{\De}\eM = \sHomb(\eN,\eR)$, $\Ed{\De}\eN=\sHomb(\eM,\eR)$.
\item[(iii)] $\eM=H^0(\sHomb(\eN,\eR))$, $\eN=H^0(\sHomb(\eM,\eR))$.
\end{enumerate}
If $S_2(\De)$ is the category of all $\De$-$S_2$ modules 
(viewed as a full subcategory of the category of coherent $\co_X$-modules) 
and $\text{coz}_\Delta^2$ represents the essential image of $S_2(\De)$ 
under $\Ed{\De}$, then
the situation is summarized by the following weakly ``commutative" diagram.
\stepcounter{thm}
\begin{equation*}\label{diag:wk}\tag{\thethm}
\xymatrix{
S_2(\Delta) \ar[rr]^{\textup{dualize}} 
\ar[d]^{\downarrow E_{\Delta}}
 && \textup{coz}^2_{\Delta} \ar[ll] \ar[d]^{\downarrow H^0} \\
\textup{coz}^2_{\Delta} \ar[u]^{H^0 \uparrow} 
\ar[rr]_{\textup{dualize}} 
 && S_2(\Delta) \ar[ll] \ar[u]^{E_{\Delta} \uparrow}
}
\end{equation*}
In terms of the discussion above the diagram, if $\eM\in S_2(\De)$ is an object in the northwest vertex, 
then its ``dual",  $\eN\in S_2(\De)$ occurs in the southeast vertex.
If $\De(\fp) ={\textup{ht}}(\fp)$ ($\fp\in \Spec{(A)}$), and $\eM=\co_X$ is $\De$-$S_2$, then $\eN$
is a canonical
module and the above relations have been established by Dibaei, Tousi
\cite{dt} and Kawasaki \cite{kw} as we pointed out earlier.

\subsection{Gorenstein complexes} 
The study of Gorenstein modules over a local ring $A$ 
was initiated by Sharp in \cite{sh1} where their first properties
were established. A non-zero finitely generated $A$-module $G$ is Gorenstein
if---when regarded as a complex---it is a Gorenstein complex in the sense of \cite[p.\,248]{RD} (see (a), (b), (c) below for an 
extension to formal schemes). In commutative algebraic terms, a non-zero
finitely generated 
$A$-module is Gorenstein if its Cousin complex 
(with respect to the height filtration) is
an injective resolution of $G$.  If $A$ admits a Gorenstein module, 
Sharp shows, $A$
is Cohen-Macaulay, the associated height function is a codimension
function on $X=\Spec{(A)}$, $\Hom(G,G)$ is free of rank $r^2$, $r>0$.
The positive integer $r$ is called the Gorenstein rank of $G$.
The module $G$ (regarded as a complex) is a dualizing complex if 
and only if $r=1$. If $A$ has a Gorenstein module then it has 
one of rank $r=1$ if and only if $A$ is the homomorphic 
image of a Gorenstein ring, if and only if
$A$ has a dualizing complex. In \cite{ffgr}, Fossum, Foxby, Griffith
and Reiten show that if $G$ is Gorenstein of minimal rank, then every
Gorenstein module on $A$ is of the form $G^s$ for some $s\ge 1$. 
This last result was anticipated in \cite{sh2} by Sharp in the instance
when $A$ is a complete Cohen-Macaulay ring, so that, by Cohen's structure 
theorem, $A$ is the homomorphic image of a Gorenstein ring, and whence
has a Gorenstein module of rank $r=1$, necessarily of minimal rank. 
(Cf.~also \cite{sh4} for related results.) 
In addition to the above mentioned results in \cite{ffgr}, 
Fossum {\emph{et.~al.}} also show that if $A$ has 
a Gorenstein module, then some standard {\'e}tale
neighborhood of $A$ has a Gorenstein module of rank $r=1$ 
(i.e.~a Gorenstein module which is also a dualizing complex). 

Consider a pair $(\X,\De)$ where $\X$ is a formal scheme, universally catenary,
of finite Krull dimension and $\De$ a codimension function on $\X$. Let $\Dc(\X)$ 
be the full subcategory of $\D(\X)$ (the derived category of the
category of complexes of $\co_\X$-modules) consisting of objects whose cohomology
sheaves are coherent, and let $\Dcs(\X)$ denote its essential image under the functor
$\R\iGp{\X}|_{\Dc}$, where (with $\I$ a defining ideal of $\X$)
\[
\iGp{\X} \set \dirlm{n}\sHom_{\co_\X}(\co_\X/\I^n,\,{\boldsymbol{-}}).
\]
(If $\X$ is ordinary, then $\Dcs(\X)=\Dc(\X)$.)
A complex $\eG$ is said to be c-Gorenstein on $(\X,\De)$ (or $\De$-Gorenstein) if
\begin{enumerate}
\item[(a)] $\eG\in\Dc(\X)$.
\item[(b)] $\R\iGp{\X}\eG \cong \Ed{\De}\R\iGp{\X}\eG$ in $\D(\X)$, 
i.e.~ $\R\iGp{\X}\eG$ is Cohen-Macaulay
on $(\X,\De)$.
\item[(c)] $\Ed{\De}\R\iGp{\X}\eG$ consists of $\Aqct(\X)$-injectives, where $\Aqct(\X)$
is as in \S\S\,\ref{ss:cat}.
\end{enumerate}
For the rest of this discussion, for simplicity, we assume that our
complex $\R\iGp{\X}\eG$ is non-exact on every connected component
of $\X$, equivalently $E\R\iGp{\X}\G\neq 0$ on any 
connected component of $\X$.
In this paper, using this result, we show that if
$\X$ has a c-dualizing complex, then 
\[ \eG \cong \eD\otimes {\mathcal{V}} \tag{*}\]
where $\eD$ is a c-dualizing complex whose associated codimension
function is $\De$ and ${\mathcal{V}}$ is a coherent locally free 
$\co_\X$-module (cf.~\Tref{thm:gorenstein}). 
Note that it follows that $\R\sHomb(\eG,\eG)$ is
isomorphic in $\D(\X)$ to ${\mathcal{V}}^*\otimes{\mathcal{V}}$, i.e.
to a coherent locally free $\co_\X$-module of rank $r^2$, where $r$ is
the rank of ${\mathcal{V}}$. Since $\X$ is not assumed to be connected,
we have to interpret $r$ as a locally constant, positive integer valued
function. 

Suppose we drop the assumption that $\X$ has a c-dualizing complex. Can $r$ 
(the ``rank" of $\eG$) still be defined?  In \Pref{prop:gorenstein} (and its proof)
we show that 
$\R\sHomb(\eG,\eG)$ is isomorphic
(in $\D(\X)$) to a coherent locally free sheaf ${\mathcal{W}}$ of rank $r^2$ 
where $r$ is a positive integer valued function.  
In fact, for a point $x\in \X$, $r(x)$ is the number of copies of the injective
hull $E(x)$ of the residue field $k(x)$ (thought of as a $\co_{\X,x}$-module)
in the injective module $G(x)={\Hr}^{\De(x)}_x(\R\iGp{\X}\eG)$. The result
implies that this number (of copies of $E(x)$ in $G(x)$)
is constant on connected components of $\X$, something which
is not {\emph{a priori}} obvious. Further, when $r=1$, $\eG$ is c-dualizing.
We  study the (possibly) non-commutative
$\co_\X$-algebra ${\eA}=\sHom(E\R\iGp{\X}\eG,E\R\iGp{\X}\eG)$ 
(isomorphic as a coherent sheaf to ${\mathcal{W}}$), for it sheds
light on the existence of {\'e}tale open sets of $\X$ on which
$\eG$ ``untwists" and reveals itself in the form (*). In fact one 
can show that $\eA$ is a sheaf of Azumaya algebras (see \Pref{prop:azumaya}), whose splitting
implies the existence of a dualizing complex (see \Tref{thm:azumaya1}). 
This generalizes \cite[p.\,209,\,Cor.\,(4.8)]{ffgr}. One consequence is this: 
if $(\X,\,\De)$ has a c-Gorenstein
complex, then $\Ed{\De}(\eF) \in \Dcs$ if $\eF\in\Dcs$ (cf. \Tref{thm:azumaya2}\,(b)). In particular if
$\X$ is an ordinary scheme (so that $\Dcs=\Dc$) possessing a c-Gorenstein complex and
if $\eF$ has coherent cohomology, then so does its Cousin complex with respect $\De$.

Now suppose $\X=\Spec{(A)}$, where $A$ is a local ring of dimension $d$.
Suppose further that $\X$ has  a dualizing complex
and  that $\co_\X$ is $S_2$. It is not hard to show that this
forces $\co_\X$ to be $\De$-$S_2$ where $\De$ is the codimension function $\fp\mapsto d-\dim{A/\fp}$. 
Moreover, in this case it is well-known (see, for example, \cite[p.\,23,\,Rmk.,2.1]{dt}) that
$\De$ is the ``height function" $\fp\mapsto{\text{ht}}_A(\fp)$.
 Let $\eD$ be a dualizing complex whose associated codimension 
function is $\De$ (under our hypothesis, such a $\eD$ exists). 
Let $\eK\set H^0(\eD)$. Now  (*) combined with \Tref{thm:main1} gives us
that if $\eG$ is Gorenstein with respect to
the $\De$, then $\eN\set H^0(\eG)$ is also $S_2$ and 
$\eN=\eK\otimes {\mathcal{V}}$. We believe this gives a more natural 
interpretation of \cite[p.\,125,\,Thm.\,3.3]{dib}. 

In \Ssref{ss:sharp} we discuss (very briefly) the
relationship between various results in this paper (especially \Tref{thm:gorenstein} and \Tref{thm:azumaya2})
and Gorenstein {\emph{modules}}.

\subsection{Duality theory} The paper \cite{dcc} is concerned with
studying ``the gap" between the Cousin complex $\sh{f}\eF$
constructed in \cite{lns} and the twisted inverse image $f^!\eF$
(for a map $f\colon \X\to \Y$ which is a composite of compactifiable maps
and a torsion Cousin complex $\eF$ on
$\Y$). This is done via a functorial map $\gamma_f^!(\eF)\colon
\sh{f}\eF \to f^!\eF$ which is defined for a pseudo-proper map by the universal
properties of $f^!$ (since Sastry shows that for such a map there is a map of
complexes ${\text{Tr}}_f\colon f_*\sh{f}(\eF)\to \eF$), and then for maps $f$ which are
composites of compactifiable maps. The main result is that
$E(\gamma_f^!)$ is an isomorphism, and the hardest technical step is in showing that
$\gamma_f^!$ is an isomorphism when $f$ is smooth. From this a number of results
follow: among them the just mentioned fact that $E(\gamma_f^!)$ is an isomorphism,
$\gamma_f^!$ is a functorial isomorphism when $f$ is flat, and that
$\gamma_f^!(\eF)$ is an isomorphism when $\eF$ is residual (or more generally, in the
language of the present paper, when $\eF$ is t-Gorenstein). This last result is crucial
in establishing an explicit isomorphism between the functor $f^!$ (restricted to complexes
satisfying certain coherence conditions) and a functor $f^{(!)}$ defined along the lines
of Grothendieck's original treatment of his duality on ordinary schemes, as laid out
in \cite{RD} (i.e., by first ``dualizing" on the base using a residual complex $\eR$, then 
applying $\bL f^*$, and then ``dualizing" on the source using $\sh{f}(\eR)$, with the identification
of $\sh{f}(\eR)$ with $f^!\eR$ via $\gamma_f(\eR)$ being needed crucially).
Now, Cousin complexes are equivalent to
Cohen-Macaulay complexes. Therefore there is a duality (i.e.~ an
anti-equivalence) between Cousin complexes in $\Dcs$ and coherent
sheaves. It is natural to ask for dual notions corresponding
to $\sh{f}$ and $f^!$ (restricted to Cousins in $\Dcs$). We show
that the corresponding functors on coherent sheaves are $f^*$ and
$\bL f^*$. \Tref{thm:CMness}\,(iii) and (iv) together with 
\Tref{thm:sh-(sh)} should be regarded as the precise
formulation of this statement.
Thus, if $\eF$ is Cousin on $\Y$ and in $\Dcs(\Y)$, and $\eM$
the associated coherent sheaf (under our duality), then the gap
between $\sh{f}\eF$ and $f^!\eF$ is equivalent---in the dual
situation---to the gap between $f^*\eM$ and $\bL f^*\eM$. The
comparison map $\gamma_f^!\colon \sh{f} \to f^!$ corresponds
to the natural transformation $\bL f^* \to f^*$ on coherent 
$\co_\Y$-modules. If $f$ is flat, this means that the gap can be
closed for all Cousins $\eF$ in $\Dcs(\Y)$ (i.e.~for all coherent
$\eM$ on $\Y$) and vice-versa. This gives a natural interpretation
of the result in \cite[p.\,182,\,7.2.2]{dcc} (cf.~\Tref{thm:dcc-flat} together
with \Tref{thm:sh-(sh)}). 
In general, the condition that $\sh{f}\eF\cong f^!\eF$ imposes 
conditions on the pair $(f,\,\eF)$, whence on $(f,\,\eM)$. 
We interpret this in terms of Tor-independence (cf.~\Dref{def:torind}
and \Lref{lem:torind}).
There are drawbacks to the approach taken in this paper. 
We have to restrict ourselves to complexes satisfying certain 
coherence conditions (they should be in $\Dcs$) 
and to schemes carrying c-dualizing complexes, whereas Sastry works with any
torsion Cousin complex, and without the use of c-dualizing complexes.
Even with these restrictions, we point out to the readers that any seeming simplicity in the
proofs here is nullified by the realization that we do not use $f^\sharp$,
$f^!$ and $\gamma_f^!$ for our ``simpler" proofs, but use instead the analogous functors (and natural maps)
$f^{(\sharp)}$, $f^{(!)}$, and  $\gamma_f^{(!)}$. To show (as we do in\Tref{thm:sh-(sh)} and
Diagram \eqref{diag:zeta-phi}) that the analogous
functors and maps are actually isomorphic, one requires some of the deeper results
in \cite{dcc}. So the point of the results on Grothendieck duality in this paper is to establish the 
dictum ``$\sh{f}$ is to
$f^!$ as $f^*$ is to $\bL f^*$" i.e.  ``$\gamma_f^!$ is dual to $\bL f^* \to f^*$".
There are gains in doing this, for, after replacing $f^!$ by its variant
$f^{(!)}$ \cite[\S\,9]{dcc}, $\sh{f}$ by its variant $\psh{f}$, and $\gamma_f^!$ by $\gamma_f^{(!)}$,
we are able to extend the results in 
\cite{dcc} to arbitrary pseudo-finite
type maps between the allowed schemes, whereas in \cite{dcc}, 
Sastry had to restrict himself to maps which were composites of 
compactifiable maps. 

\section{Preliminaries}

In this paper, all schemes---ordinary or formal---are noetherian.

\subsection{Basic Terminology; Cousin complexes}\label{ss:codim} A codimension function on a 
formal scheme $\X$ is a map 
\[ \De\colon |\X|\to {\mathbb{Z}} \]
such that $\De(x')=\De(x)+1$ for every immediate specialization $x'$ of a point $x\in \X$. Here,
$|\X|$ denotes the set of points underlying the scheme $\X$. Let $\bbF$ denote the category
whose objects are noetherian universally catenary formal schemes admitting a codimension
function and whose morphisms $\X'\to \X$ are formal scheme maps which are {\emph{essentially
of pseudo-finite type}}, i.e., if $\I\subset \co_\X$ and $\J\subset \co_{\X'}$
are ideals of definition of $\X$ and $\X'$ respectively, with 
$\J\supset \I\co_{\X'}$, then the map ${\mathbf{Spec}}(\co_{\X'}/\J)\to
{\mathbf{Spec}}(\co_\X/\I)$ is a localization of a finite type map of ordinary
schemes.

%after modding out defining ideals in $\X'$ and $\X$, the map is a
%localization of a finite type map of ordinary schemes.

As in \cite{lns}, we often work with a slightly more refined category $\bbFc$ consisting of objects
$(\X,\De)$, where $\X\in \bbF$ and $\De$ is a codimension function on $\X$, whose morphisms
$(\X',\,\De')\to (\X,\,\De)$ are maps $f\colon \X'\to \X$ in $\bbF$ such that for $x'\in\X'$ and $x=f(x')$,
$\De(x)-\De'(x')$ is equal to the transcendence degree of the residue field $k(x')$ of $x'$ over the
residue field $k(x)$ of $x$.

Let $\X\in \bbF$ and let $x\in\X$. For any abelian group $G$, $i_xG$ is the extension by zero
to $\X$ of the constant sheaf ${\overline{G}}$ modeled on $G$, on the closure $\overline{\{x\}}$ of
$\{x\}$.

Let $(\X,\De)\in\bbFc$. A complex $\eF$ of $\co_\X$-modules is called a {{\emph{Cousin complex
on}}}  $(\X,\De)$ or a $\De$-Cousin complex (or simply a Cousin complex with respect to $\De$)
if, for each $n\in{\mathbb{Z}}$, $\eF^n$ is the direct sum of a family of $\co_\X$-modules
$(i_xF_x)_{x\in\X,\De(x)=n}$, where $F_x$ is an $\co_{\X,x}$-module.
We refer the reader to \cite[pp.\,36--44,\,\S\S\,3.2 and 3.3]{lns} 
for  more elaborate definitions of Cousin (and Cohen-Macaulay) complexes with respect to $\De$.
We do make fleeting references to Cousin complexes with respect to descending filtration 
$\eH=(H_i)_{i\ge i_\circ}$ of closed sets $H_i$ in $\X$. In this case $\eF$ is Cousin with respect to
$\eH$ if $\eF^n$ is
the direct sum of a family $(i_xF_x)_{x\in\partial{H_i}}$, where $F_x$ is an $\co_{\X,x}$-module and
$\partial{H_i}= H_i-H_{i+1}$. However our emphasis will always be on the more special Cousin 
complexes which are associated to a codimension function $\De$.

\subsection{Categories of complexes}\label{ss:cat}
For a formal scheme $\X$, let $\A(\X)$ be the category of
$\co_\X$-modules, and $\Aqc(\X)$ (resp.~$\Ac(\X)$, resp.~$\Avc(\X)$)
the full subcategory of $\A(\X)$ whose objects are the quasi-coherent
(resp.~ coherent, resp.~ $\dirlm{}$'s of coherent) $\co_\X$-modules.
As in \cite[p.\,6,\,1.2.1]{dfs}, we define the torsion functor $\iGp{\X}\colon
\A(\X)\to \A(\X)$ by the formula
\[
\iGp{\X} \set \dirlm{n}\sHom_{\co_\X}(\co_\X/\I^n,\,{\boldsymbol{-}})
\]
where $\I\subset \co_{\X}$ is an ideal of definition of $\X$. The
definition of $\iGp{\X}$ is independent of the choice of $\I$. Note
that $\iGp{\X}$ is a subfunctor of the identity functor.

An object $\eM\in\A(\X)$ is called a {\emph{torsion}} $\co_\X$-module
if $\eM=\iGp{\X}(\eM)$. We denote by $\At(\X)$ (resp.~$\Aqct(\X)$) the
full subcategory of $\A(\X)$ consisting of torsion (resp.~quasi-coherent
torsion)
$\co_{\X}$-modules.

Let $\C(\X)$ be the category of $\A(\X)$-complexes, $\K(\X)$ the associated
homotopy category, and $\D(\X)$ the corresponding derived category, obtained
from $\K(\X)$ by inverting quasi-isomorphisms.

For any full subcategory $\A_{\dots}(\X)$ of $\A(\X)$, denote by
$\C_{\dots}(\X)$ (resp.~ $\K_{\dots}(\X)$, resp.~ $\D_{\dots}(\X)$)
the full subcategory of $\C(\X)$ (resp.~ $\K(\X)$, resp.~ $\D(\X)$)
whose objects are those complexes whose cohomology sheaves all lie
in $\A_{\dots}(\X)$, and by $\D_{\dots}^+(\X)$ (resp.~ $\D_{\dots}^-(\X)$,
resp.~ $\D_{\dots}^b(\X)$) the full subcategory of $\D_{\dots}(\X)$
whose objects are complexes $\eF\in \D_{\dots}(\X)$ such that the
cohomology $H^m(\eF)$ vanishes for all $m\ll 0$ (resp.~ $m\gg 0$,
resp.~ $|m|\gg 0$). We often write $\Dc\>$, $\Dqct\>$, \dots for $\Dc(\X)$,
$\Dqct(\X)$, \dots when there is no danger of confusion.

The essential image of $\R\iGp{\X}|_{\Dc}$ is of considerable interest
to us, and as in \cite[\S\S\,2.5,\,p.\,24,\,second paragraph]{dfs} 
we denote it by $\Dcs(\X)$. In greater
detail, $\Dcs(\X)$ is the full subcategory of $\D(\X)$ such that
$\eE\in \Dcs(\X) \Leftrightarrow \eE\cong \R\iGp{\X}\eF$ with $\eF\in\Dc(\X)$.
It is immediate that $\Dcs(\X)$ is a triangulated subcategory of $\D$
or~$\Dqct$.
%In this paper, the importance of $\Dcs(\X)$ stems from the duality
%(in the presence of mild hypotheses on $\X$) between Cohen-Macaulay
%complexes (with respect to a fixed codimension function) in $\Dcs(\X)$
%and $\Ac(\X)$. This duality was first observed by Yekutieli and
%Zhang in \cite{yz} for ordinary schemes (where $\Dc(\X)=\Dcs(\X)$)
%of finite type over a regular scheme, and later in greater generality
%by Lipman and the authors in \cite{lns}.

\subsection{Dualizing complexes} As shown in \cite[p.\,26,\,Lemma\,2.5.3]{dfs},
the notion of a dualizing complex
on an ordinary scheme breaks up into two related notions on a formal
scheme. We recall here the definitions and first properties from 
\cite[p.\,24,\,Definition\,2.5.1]{dfs}.

\begin{defi}\label{def:dualizing} A complex $\eR$ is a {\emph{c-dualizing
complex on}} $\X$ if
\begin{enumerate}
\item $\eR\in \Dc^+$.
\item The natural map $\co_\X \to \R\sHomb(\eR,\,\eR)$ is an isomorphism.
\item There is an integer $b$ such that for every 
coherent {\emph{torsion}} sheaf $\eM$
and every $i>b$, it holds that $\Ext^i(\eM,\,\eR)\set H^i\R\sHomb(\eM,\,\eR)=0$.
\end{enumerate}
A complex $\eR$ is a {\emph{t-dualizing complex}} on $\X$ if
\begin{enumerate}
\item $\eR\in \Dt^+$.
\item The natural map $\co_\X \to \R\sHomb(\eR,\,\eR)$ is an isomorphism.
\item There is an integer $b$ such that for every coherent {\emph{torsion}}
sheaf $\eM$ and for every $i>b$, $\Ext^i(\eM,\,\eR)\set H^i\R\sHomb(\eM,\,\eR)
=0$.
\item For some ideal of definition $\J$ of $\X$, $\R\sHomb(\co_\X/\J,\,\eR)\in
\Dc(\X)$ (equivalently,
$\R\sHomb(\eM,\,\eR)\in\Dc(\X)$ for every coherent torsion sheaf $\eM$.)
\end{enumerate}
\end{defi}

Yekutieli was the first to consider t-dualizing complexes in \cite{y}, where they were
simply called ``dualizing".
We note from \cite[2.5.3 and 2.5.8]{dfs} that $\X$ has a c-dualizing
complex if and only if $\X$ has a t-dualizing complex which lies
in $\Dcs$. In greater detail, if $\eR$ is a c-dualizing complex, then
$\R\iGp{\X}\eR \in \Dcs$ is a t-dualizing complex. Conversely, if $\eR$
is a t-dualizing complex that lies in $\Dcs$, then
$\R\sHomb(\R\iGp{\X}\co_\X,\eR)$ is c-dualizing. 
For an ordinary scheme, 
$\Dt=\D$ and $\Dc=\Dcs$ and the notions of a c-dualizing complex
and a t-dualizing complex
coincide with the usual notion of a dualizing complex. We point out that according
to \cite[p.\,106,\,Prop.\,9.2.2]{lns}, a t-dualizing complex lies in $\Dqct^b(\X)$.

\begin{exam}\label{ex:dualizing}  Let $X$ be an ordinary scheme 
and $\kappa \colon \X \to X$ its completion along a closed subscheme $Z$.
Then for any dualizing complex $\cR$ on $X$, $\kappa^*\cR$ 
is c-dualizing on $\X$ and 
$\R\iGp{\X}\kappa^*\cR \cong \kappa^*\R\iG{Z}\cR$ is t-dualizing
and lies in $\Dcs$ \cite[p.\,25,\,2.5.2(2)]{dfs}. 
In particular, if $k$ is a field and $\X$ is the formal spectrum 
of $A \set k[[X_1,\ldots,X_n]]$ equipped with 
the $\fm$-adic topology where $\fm = (X_1,\ldots,X_n)$, 
(which implies that~$\X$ consists of a single point)
then a c-dualizing complex on $\X$ is given by $A$ 
while a t-dualizing complex
is given by the injective hull of $k = A/\fm$.
\end{exam}

For a fixed t-dualizing complex $\eR$ on $\X$ define the
dualizing functor $\cDt={\cDt}(\eR)\colon \D\to \D$ by
\[\cDt\set\R\sHomb({\boldsymbol{-}},\,\eR).\] 
If $\eR \in \Dcs$---equivalently, if
$\X$ has a c-dualizing complex---then according to 
\cite[p.\,28,\,Prop.\,2.5.8]{dfs}
\begin{enumerate}
\item $\eE\in \Dcs \Leftrightarrow \cDt\eE\in\Dc$.
\item $\eF\in \Dc \Leftrightarrow  \cDt\eF\in\Dcs$.
\item If $\eF$ is in either $\Dc(\X)$ or $\Dcs(\X)$, the natural
map is an isomorphism:
\stepcounter{thm}
\begin{equation*}\label{eq:Dt}\tag{\thethm}
\eF \iso \cDt\cDt\eF.
\end{equation*}
\end{enumerate}

The above facts can be summarized as follows:

\begin{prop}\cite[p.\,28,\,Prop.\,2.5.8]{dfs}\label{prop:Dc-Dcs} Let 
$\X$ be a formal scheme
with a t-dualizing complex $\eR\in\Dcs(\X)$. Then
the functor $\cDt$ induces, in either direction, an anti-equivalence
of categories between $\Dc(\X)$ and $\Dcs(\X)$.
\end{prop}

Regarding $\Ac(\X)$ as a full subcategory of $\Dc(\X)$, 
\cite[p.\,108,\,Cor.\,9.3.2]{lns}
characterizes the essential image of $\Ac(\X)$ in $\Dcs(\X)$
under the above anti-equivalence, and this characterization 
underpins most of the results in this paper. 
We describe this in subsection \Ssref{ss:anti}

\subsection{Cohen-Macaulay and Cousin complexes} Let $(\X,\De)\in\bbFc$. We say that a complex
$\eF\in\D^+(\X)$ is {\emph{Cohen-Macaulay}} with respect to $\De$ (or Cohen-Macaulay
on $\De$) if, for any $x\in\X$, $\Hr^i_x(\eF)=0$ for $i\neq \De(x)$. Here $\Hr^i_x\eF$ is defined
to be (the ``hyper local cohomology" module) $\Hr^i(\R\Gamma_x\eF)$ (sometimes denoted
${\mathbb{H}}^i_x(\eF)$). We (again) refer the reader to \cite[pp.\,36--44,\,\S\S\,3.2 and 3.3]{lns} 
for  more elaborate definitions
of Cousin complexes and Cohen-Macaulay complexes. 

Cohen-Macaulay complexes and Cousin complexes are intimately related. In 
\cite[Theorem\,3.9]{suominen} Suominen shows (and this, in a more general form, is the main 
result of that paper) that the full subcategory of $\D(\X)$ consisting of Cohen-Macaulay complexes on
$(\X,\,\De)$ is equivalent to the full subcategory of $\C(\X)$ consisting of Cousin complexes
on $(\X,\,\De)$ via the restriction of the localization functor $Q\colon \K(\X)\to \D(\X)$ to the
category of Cousin complexes on $(\X,\De)$.
 The category of $\De$ Cousin complexes can be
regarded as a full subcategory of $\K(\X)$ since any two maps between $\De$ Cousin complexes
which are homotopic to each other are actually equal.\footnote{In
fact the results of Suominen show that Cousin complexes on $(\X,\,\De)$ can be regarded as
a full subcategory of $\D(\X)$.}  A pseudo-inverse for
$Q$ restricted to $\De$ Cousin complexes is provided by the restriction of the Cousin
functor $\Ed{\De}\colon \D^+(\X) \to \C(\X)$ to $\De$ Cohen-Macaulay complexes.  (See 
\cite[Prop.\,3.3.1,\,p.\,42]{lns}.)

For this paper, we are not interested in all Cohen-Macaulay or all Cousin complexes on $(\X,\De)$. In\cite{lns} and \cite{dcc} the interest was often in Cohen-Macaulay (resp. Cousin) complexes in
$\Dqct(\X)$. In this paper our interests are more special. We will concentrate on Cohen-Macaulay
complexes in $\Dcs (\subset \Dqct)$. To that end we denote by $\cm(\X,\De)$ the full subcategory
of $\Dcs(\X)$ consisting of Cohen-Macaulay complexes on $(\X,\De)$ and we denote by 
$\Cozt{\De}(\X)$ the full subcategory of the category of Cousin complexes 
on $(\X,\De)$ which corresponds to $\cm(\X,\De)$ under Suominen's equivalence 
above. Note that any complex in $\Cozt{\De}(\X)$ consists of $\Aqct$-modules,
see the paragraph above Lemma 3.2.2 in \cite[page 40]{lns}.

In \cite{dcc} we used the symbols ${\mathrm{CM}}^*$ and $\Coz^*_{\De}$ for 
$\cm$ and $\Cozt{\De}$ respectively.

%We remind the
%reader that there is an equivalence of categories between the
%full subcategory ${\mathrm{Cou}}(\X;\De)\subset \C(X)$ consisting
%of $\De$-Cousin complexes and the full subcategory 
%$\D^+(\X;\De)_{\cm}\subset \D^+(\X)$ of $\De$-CM complexes. Its brief
%description is as follows. First note that ${\mathrm{Cou}}(\X;\De)$ is
%also a full subcategory of $\K$. Indeed (see 
%\cite[top of p.\,42,\,\S\S\,3.3]{lns}), a map of $\De$-Cousin complexes homotopic
%to zero is already the zero map. If $Q=Q_\X\colon \K\to \D$ denotes
%the usual localization functor, then $\D^+(\X;\De)_{\cm}$ is the 
%essential image of ${\mathrm{Cou}}(\X;\De)$ in $\D^+$ under $Q$,
%and the resulting functor
%\[
%Q|_{{\mathrm{Cou}}}\colon {\mathrm{Cou}}(\X;\De)\to \D^+(\X;\De)_{\cm}
%\]
%is an equivalence of categories. An  inverse equivalence is given
%by the {\emph{restriction}} of the Cousin functor
%\[
%\Ed{\De}\colon \D^+(\X) \to {\mathrm{Cou}}(\X;\De)
%\]
%to $\D^+(\X;\De)_{\cm}$ (cf. \cite[Thm.\,3.9]{suominen} and 
%\cite[p.\,42,\,Prop.\,3.3.1]{lns}).

%We set $\cm(\X;\De)\set \D^+(\X;\De)_{\cm}\cap\Dqct$ and 
%$\cm^*(\X;\De)\set \cm(\X;\De)\cap\Dcs$. On the Cousin side
%we first set
%$\Cozs{\De}(\X)\set {\mathrm{Cou}}(\X;\De)\cap\Cqct(\X)$
%and note that in view of \cite[p.\,40, (12)]{lns},
%$\Cozs{\De}(\X)$ corresponds to $\cm(\X;\De)$
%through $Q$ and $\Ed{\De}$. Next we define
%$\Cozt{\De}(\X)$ to be the full subcategory of $\Cozs{\De}(\X)$
%which corresponds to $\cm(\X;\De)$ under the equivalence
%above.

\subsection{An anti-equivalence}\label{ss:anti} We are now in a position to identify the subcategory of 
$\Dcs(\X)$ which corresponds to $\Ac(\X)\subset \Dc(\X)$ under the anti-equivalence
of \Pref{prop:Dc-Dcs}.
First, given a t-dualizing complex $\eR$ on $\X$, one has an associated
codimension function $\De_\eR$ \cite[p.\,106,\,9.2.2(ii)(b)]{lns}. Moreover,
$\eR$ {\emph{is Cohen-Macaulay with respect to $\De_{\eR}$}} 
\cite[Prop.\,9.2.2(iii)(a)]{lns}.
According to \cite[p.\,108,\,Prop.\,9.3.1 and Cor.\,9.3.2]{lns} we have:

\begin{prop}\cite{lns}\label{prop:Ac-CM} Let $\X$ be a formal scheme
with a t-dualizing complex $\eR\in\Dcs(\X)$. Let
$\De=\De_\eR$ be the associated codimension function. Then
the functor $\cDt$ induces, in either direction, an anti-equivalence
between $\Ac(\X)$ and $\cm(\X;\De)$. Thus there
exists a commutative diagram as follows,
with $\equiv$ denoting equivalence of categories,
the vertical arrows being the standard inclusions,
and $C^\circ$ denoting the category opposite to the category $C$:
\[
\xymatrix{
\Dc(\X) \ar[rr]^{\equiv}_{\cDt} & & \Dcs(\X)^\circ \\
\Ac(\X) \ar[u] \ar[rr]^{\equiv}_{\cDt} & & \cm(\X;\De)^\circ \ar[u]
}
\]
\end{prop}
\Pref{prop:Ac-CM} was first proved by Yekutieli and Zhang 
\cite[Thm.\,8.9]{yz} for ordinary schemes of finite type 
over noetherian finite dimensional regular rings.

\subsection{Residual complexes}
Since $\cm(\X;\De)$ is equivalent to $\Cozt{\De}(\X)$, one can
restate the anti-equivalence between $\Ac(\X)$ and $\cm(\X;\De)$
in terms of $\Cozt{\De}(\X)$. The resulting anti-equivalence
between $\Ac(\X)$ and $\Cozt{\De}(\X)$ can be stated entirely
in terms of complexes, i.e., within $\C(\X)$ rather than $\D(\X)$.
As a first step toward this, we discuss the notion of a
residual complex on a formal scheme.

On an ordinary scheme, we refer to \cite[p.\,304]{RD} for a
definition of a residual complex. 
Following \cite[p.\,104,\,9.1.1]{lns}, by a residual
complex on a formal scheme~$\X$ we mean a complex
$\eR$ of $\At$-modules such that there exists a defining 
ideal $\I\subset\co_\X$
with the property that for any $n>0$, $\sHomb(\co_\X/\I^n,\,\eR)$
is residual on the ordinary scheme $(\X,\,\co_\X/\I^n)$.
(cf.~[Ye], \cite[\S9, footnotes]{lns}) 
The residual complex~$\eR$ induces a codimension function~$\De=\De_\eR$
on~$\X$, and $\eR$ is $\De$-Cousin consisting of $\Aqct$-modules
\cite[p.\,106,\,Prop.\,9.2.2]{lns}.
Moreover, if $\X$ admits a residual complex, then $\X$ is
universally catenary (since the corresponding statement is
true for ordinary schemes).

According to \cite[Prop.\,9.2.2(iii)]{lns}, if $\eD$ is t-dualizing and
$\De=\De_\eD$, then $\eR\set \Ed{\De}\eD$ is a residual complex. 
%(see \cite[\S\,9.1]{lns} for a definition). 
Moreover, since $\eD$ is Cohen-Macaulay on $(\X,\,\De)$, there is a canonical
isomorphism between $\eD$ and $Q\eR$ 
(see \cite[p.\,42,\,3.3.1 and 3.3.2]{lns}).
Moreover, it is immediate from the definition of $\De_\eD$ that
$\De_\eD=\De_\eR$.
Since the presence
of a t-dualizing complex forces $\X$ to be of finite Krull dimension
\cite[p.\,106,\,Prop.\,9.2.2(ii)]{lns}, $\eR$ must be a bounded complex.
Conversely, if $\eR$ is a bounded residual complex, then $Q\eR$
is t-dualizing \cite[Prop.\,9.2.2(ii) and (iii)]{lns}.

We need a little more terminology which will facilitate discussions on
Cousin complexes. As in \cite{dcc}, let $\rbbF$ denote the full subcategory
of $\bbF$ whose objects $\X$  admit a
bounded residual complex $\eR$ (necessarily a t-dualizing complex)
such that $Q\eR\in\Dcs(\X)$. Note that the presence of such an
$\eR$ on $\X\in\rbbF$ is equivalent to the existence 
of a c-dualizing complex on 
$\X$. Next consider the full subcategory $\rbbFc$ of $\bbFc$ consisting
of pairs $(\X,\,\De)\in \bbFc$ with $\X\in\rbbF$. We remind the reader that
a morphisms $(\X,\De')\to (\Y,\,\De)$ in
$\rbbFc$ is therefore a map $f\colon \X\to\Y$ in $\rbbF$ and such that
for $x\in \X$ and $y=f(x)$, $\De(y)-\De'(x)$ is equal to the transcendence
degree of the residue field $k(x)$ of $x$ over the residue field $k(y)$
of $y$. In other words, if $\sh{f}\De$ is defined by the formula
\stepcounter{thm}
\begin{equation*}\label{eq:cod-f}\tag{\thethm}
\sh{f}\De(x)\set \De(y) - {\mathrm{tr.deg.}}_{k(y)}k(x) \qquad \qquad
(x\in \X, y\set f(x))
\end{equation*}
then $\De'=\sh{f}\De$. One checks that $\sh{f}\De$ is always a codimension
function on $\X$.
If $(\X,\,\De)\in\rbbFc$ then a Cohen-Macaulay (resp.~Cousin) complex
on $(\X,\,\De)$ is a complex in $\cm(\X;\De)$ (resp.~$\Cozt{\De}(\X)$).

\subsection{Cousin complexes and coherent sheaves}
%For the rest of this paper {\emph{assume that all schemes occurring
%admit c-dualizing complexes,}} (so that every t-dualizing complex lies
%in $\Dcs$), i.e. all schemes are in $\rbbF$ (see sentence following
%\Dref{def:dualizing}).
%
As seen so far,
for any~$\X$ admitting a c-dualizing complex (equivalently, admitting 
a t-dualizing complex in~$\Dc^*$) with associated codimension
function $\De$, the category $\cm(\X, \De)$ is closely related to
two abelian categories, namely, $\Ac(\X)$ via dualizing and
$\Cozt{\De}(\X)$ via the cousin functor $E_{\De}$. We now relate
these two abelian categories directly.

Fix $(\X,\,\eR)$ with $\eR$ a residual complex on the formal scheme $\X$
and set $\De=\De_\eR$ %$\Coz(\X)=\Cozs{\De}(\X)$ 
and $\Cozst(\X)=\Cozt{\De}(\X)$. By \cite[p.\,104,\,Lemma\,9.1.3]{lns} 
and \cite[p.\,123]{RD},
we see that $\eR$ is a complex of $\Aqct(\X)$-injectives. 
Let $\eF$ be a complex of $\Avc$-modules. (Recall that $\Avc$-modules are 
direct limits of coherent modules and there are inclusions 
$\Ac \subset \Avc$, $\Aqct \subset \Avc$.)
For any such $\F$ we shall now make the identification
\[
\cDt\eF = \sHomb_{\co_\X}(\eF,\,\eR).
\]
This can be justified using \cite[page 27, 2.5.6]{dfs} 
(where the proof holds only when~$\eF$ is a complex of $\Avc$-modules and
not for $\eF \in \Dvc$ as claimed). 
%For simplicity we drop
%occurrences of the prefix $Q$ from the objects involved. 
Using this version of~$\cDt$, we make the following three observations:

\noindent
1) For $\eM\in\Ac(\X)$, the complex
\[ \eM'\set \sHomb_{\co_\X}(\eM,\,\eR)\]
lies in $\Cozst(\X)$. 
%and \[\eM' = \cDt\eM.\] The first assertion 
That $\eM'$ is a Cousin complex follows easily from the definitions
but to see that its image in the derived category lies in $\Dc^*$
we need \Pref{prop:Ac-CM}, and the equivalence
between $\Cozst(\X)$ and $\cm(\X;\De)$.

\noindent
2) For $\eF\in\Cozst(\X)$, the $\co_{\X}$-module
\[\eF^*\set \sHom_{\C(\X)}(\eF,\,\eR)\]
lies in $\Ac(\X)$  and there is a natural quasi-isomorphism
\stepcounter{thm}
\begin{equation*}\label{iso:*Dt}\tag{\thethm}
\eF^*\to \cDt\eF.
\end{equation*}
Indeed, note that for any object $\eG\in\Cozst(\X)$, we
have
\[\sHom(\eG,\,\eR) = H^0(\sHomb(\eG,\,\eR)) = H^0(\cDt \eG).\]
(To see the first equality, note that the only $\C(\X)$-map
$\eG\to \eR$ homotopic to zero is the zero map, for $\eG$ and $\eR$ are
$\De$-Cousin.)
The assertions for $\eF^*$ follow from
\Pref{prop:Ac-CM} and the fact that $\Cozst(\X)$ is equivalent to
$\cm(\X;\De)$. %(via $Q$ and $\Ed{\De}$).

\noindent
3) The operations $\eM\mapsto \eM'$ and $\eF\mapsto \eF^*$ are inverse
operations. In greater detail:
\begin{enumerate}
\item[(i)] For $\eM\in\Ac(\X)$, the natural map in $\Ac(\X)$ given by
``evaluation" is an isomorphism
\stepcounter{thm}
\begin{equation*}\label{eq:M}\tag{\thethm}
\eM \iso (\eM')^*.
\end{equation*}
Indeed, in $\Dc$, the above map is equivalent to \eqref{eq:Dt}.
\item[(ii)] For $\eF\in\Cozst(\X)$, the natural map in $\Cozst(\X)$
given by ``evaluation" is an isomorphism
\stepcounter{thm}
\begin{equation*}\label{eq:F}\tag{\thethm}
\eF \iso (\eF^*)'.
\end{equation*}
As in (i), this follows from \eqref{eq:Dt} for objects in $\Dcs$.
\end{enumerate}

Note that the correspondences $\eM\mapsto \eM'$ and $\eF\mapsto \eF^*$ 
are functorial, defining contravariant functors ${\boldsymbol{{-}'}}$
and ${\boldsymbol{{-}^*}}$. Here then is the restatement of
\Pref{prop:Ac-CM}:

\begin{prop}\label{prop:Ac-Coz} The functors
\[{\boldsymbol{{-}^*}}\colon \Cozst(\X) \to \Ac(\X)^\circ\]
and
\[{\boldsymbol{{-}'}}\colon \Ac(\X)^\circ \to \Cozst(\X)\]
are pseudoinverses via \eqref{eq:M} and \eqref{eq:F}, and therefore
set up an anti-equivalence of categories between $\Cozst(\X)$ and $\Ac(\X)$.
\end{prop}

\begin{rem}\label{rem:star-prime}{\emph{The functors $\boldsymbol{{-}'}$ 
and $\boldsymbol{{-}^*}$ depend upon the choice of $\eR$ (as we will make
explicit later in this remark).  It will be clear from the context
what the underlying residual complex is. There will be occasions when
we deal with maps $f\colon (\X,\De')\to (\Y,\De)$ in $\rbbFc$, with
a residual complex $\eR$ on $\Y$ and a residual complex
$\sh{f}\eR$ on $\X$, but even here it will be clear from the context,
which residual complex is being used and when. As an example, for 
the symbol $(f^*\eF^*)'$, it is to be assumed that the ``upper star"
occurring as a superscript of $\eF$ is with respect to $\eR$ and the
``prime" outside the parenthesis is with respect to $\sh{f}\eR$.
As for the dependence on $\eR$, if $F_\eR\colon \Ac(\X)^\circ\to \Cozst(\X)$
and $G_\eR\colon \Cozst(\X) \to \Ac(\X)^\circ$ denote
$\sHomb({\boldsymbol{-}},\eR)$ and $\sHom({\boldsymbol{-}},\eR)$
respectively, and $\eR'$ is another residual complex whose associated
codimension function is also $\De$, then 
$F_{\eR'}\cong F_\eR \otimes{\mathcal{L}}$ and
$G_{\eR'}\cong G_\eR\otimes{\mathcal{L}}$ where ${\mathcal{L}}=\sHom(\eR,\eR')$.
Note that ${\mathcal{L}}$ is an invertible $\co_\X$-module with inverse
$\sHom(\eR',\eR)$ and we have the relation 
$\eR'\cong \eR\otimes{\mathcal{L}}$.}}
\end{rem}

\section{The $S_2$ condition}

\subsection{The map $\ssf(\eG)$}  We first state a part of 
\cite[p.\,109,\,Prop.\,9.3.5]{lns} in a form that is useful to us. The content of the
Proposition is that the Cousin complex $\Ed{\De}(\eM)$ of a complex $\eM$ depends
{\em only on the zeroth cohomology of its dual} $\cDt(\eM)$.

\begin{prop}\cite{lns}\label{prop:theta} Let $\X\in\rbbF$, $\eR$ 
a residual complex which is a $t$-dualizing complex in $\Dc^*\X$, 
and let $\De=\De_\eR$. 
If $\theta\colon\eF\to \eG$
is a map in $\Dc(\X)$ such that $H^m(\theta)\colon H^m(\eF)\to H^m(\eG)$
is an isomorphism, then the induced map
\[\Ed{\De-m}(\cDt\theta)\colon \Ed{\De-m}(\cDt\eG) \iso \Ed{\De-m}(\cDt\eF).\]
is an isomorphism of Cousin complexes in $\Cozt{\De-m}(\X)$. 
In particular, truncation on $\eF$ induces natural isomorphisms
\[
\Ed{\De}(\cDt\eF) \cong \Ed{\De}(\cDt(H^0\eF))
= \Ed{\De}((H^0\eF)') \cong (H^0\eF)'.
\]
%In particular, the truncation maps $\tau_{\le m}, \tau_{\ge m}$ 
%on $\eF$ induce natural isomorphisms 
%$\Ed{\De-m}(\cDt((H^m\eF)[-m])) \iso \Ed{\De-m}(\cDt\eF)$.
\end{prop}

The proof of $\Ed{\De-m}(\cDt\theta)$ being an isomorphism
is contained in the opening paragraph of the proof
of \cite[p.\,109,\,Prop.\,9.3.5]{lns}. The fact that $\Ed{\De-m}(\cDt\theta)$ 
is in $\Cozt{\De-m}(\X)$ follows from the last part of the statement of
{\emph{loc.cit.}} The last isomorphism of the above proposition 
holds because $(H^0\eF)'$ is already a Cousin complex.

We would like to define the notion of an $S_2$ module with respect to
a codimension function. For this we need to recall certain parts
of \cite[pp.\,108--111,\,\S\S\,9.3]{lns}, especially as it relates to 
Corollary\,9.3.6 of {\emph{loc.cit}}. Let $(\X,\De)$ be 
an ordinary scheme in $\rbbFc$
and fix a residual complex $\eR$ on $\X$ that is dualizing 
and with $\De_\eR=\De$.
Let $\eM \in \Ac(\X)$ and for simplicity, 
assume that 
\[ {\mathrm{min}}\{n\,|\, H^n\eM'\neq 0\} = 0.\] 
%Let $0\neq\eG\in{\Dcs}^-(\X)$ and for simplicity, assume that 
%\[ {\mathrm{min}}\{n\,|\, H^n\cDt\eG\neq 0\} = 0.\]
%Set $H\set H^0(\cDt\eG)$ and let $\theta\colon H \to \cDt\eG$
%be the obvious canonical map in $\Dc^+$.
Let $\theta\colon H^0(\eM') \to \eM'$
be the obvious canonical map in $\Dc^+$ induced by truncation.
Then~$\theta$ induces a $\D(\X)$ map 
(cf.~\cite[p.\,109,\,9.3.6]{lns})
\stepcounter{thm}
\begin{equation*}\label{eq:SG}\tag{\thethm}
\ssf(\eM)\colon \eM \to \Ed{\De}(\eM)
\end{equation*}
defined by the natural maps
\[
\eM \iso \eM'{}^* \xrightarrow{\;\cDt(\theta)\;} (H^0(\eM'))' 
\xrightarrow[\;\ref{prop:theta}\;]{\cong} \Ed{\De}(\eM).
\]
%defined by the commutativity of
%\[\xymatrix{ \Ed{\De}(\eM) \ar[r]^-{\Iso} & \Ed{\De}(\eM'') 
%\ar[rr]^{\Iso}_{E\theta'} & & \Ed{\De}{H'} \\
%\eM \ar[u]^{\ssf(\eM)} \ar[r]^-{\Iso} & \eM''
%\ar[rr]_{\theta'} & & H' \ar[u]_{\wr} }\]
%We point out that $H'$ is already a Cousin 
%complex and from this we deduce the vertical
%isomorphism on the right (cf. \cite[p.\,42,\,Prop.\,3.3.1]{lns}).
%The second horizontal arrow on the top row is an isomorphism
%by \Pref{prop:theta}. 
In fact $\ssf(-)$ is defined for any $0\neq\eG\in{\Dcs}^-(\X)$.
We refer to \cite[p.\,109,\,Cor.\,9.3.6]{lns} for more details.

\subsection{Coherent $S_2$ sheaves on ordinary schemes} For the rest
of this section, {\emph{all schemes are ordinary}} and, as before, lie in
$\rbbF$, which translates---in this situation---to the existence of
a dualizing complex on that scheme. We first recall the ``classical" definitions
on a commutative ring.

Fix a noetherian commutative ring $A$. Recall that a finite $A$-module $M$ is said to
satisfy Serre's condition $(S_n)$ if ${\text{depth}}_{A_\fp}(M_\fp) \ge \min\{n,\,\dim{M_\fp}\}$ 
for every prime ideal $\fp$ of $A$. According to \cite[p.516,\,Example\,4.4]{sh5}, $M$ satisfies
condition $(S_n)$ if and only if its Cousin complex $\Ed{\eH}(M)$ with respect to the ``height filtration of 
$M$" $\eH=(H_i)$ with $H_i=\{\fp\in{\text{Supp}}_A(M)\,\vert\, {\text{ht}}_M(\fp)\ge 1\}$
 is exact at its 
$k$th term for $1\le k \le n-2$ and if the natural map
$M \to H^0(\Ed{\eH}(M))$ is an isomorphism\footnote{Strictly speaking, 
Sharp does not follow Hartshorne's recipe of constructing a Cousin
complex of $M$ for the height filtration $\eH$. But Sharp's 
Cousin complexes also live on the skeletons induced $\eH$ and hence
by the uniqueness part of \cite[Page 232, IV Prop 2.3]{RD} 
and by \cite[pp.\,500--501,\,\S\,1.1]{sh5}, these constructions coincide.} 
(see also \cite[p.\,621,\,Theorem\,(2.2)]{sh3}). In what
follows, we prefer to use the symbol $S_n$ without the traditional parenthesis for $(S_n)$. 

This motivates the following definition.

\begin{defi}\label{def:s2} Let $(X,\De)\in \rbbFc$ and suppose $\eR$ 
is a residual complex on $(X,\De)$. We say $\eM\in\A(X)$ is an $S_2$-module 
on $(X,\De)$ (or $S_2$ on $(X,\De)$; or simply $\De$-$S_2$) if 
\begin{enumerate}
\item[(a)] $\eM\in \Ac(X)$;
\item[(b)] ${\mathrm{min}}\{n\,|\,H^n(\eM')\neq 0\}=0$;
\item[(c)] With $\ssf(\eM)\colon \eM \to \Ed{\Delta}(\eM)$ the map
in \eqref{eq:SG}, we have
\[H^0(\ssf(\eM))\colon \eM \iso H^0\Ed{\Delta}\eM.\]
\end{enumerate}
Let the full subcategory of $\Ac(\X)$ consisting of $\De$-$S_2$ modules
be denoted $S_2(\De)$. (Note: By definition, a $\De$-$S_2$-module cannot be zero.)
\end{defi}

The connection between the notions of $\De$-$S_2$-modules and of $S_2$-modules is given
in \Pref{prop:deltaS2}. But first we wish to make a couple of remarks.
%\pagebreak[3]
%\newpage
\begin{rems}\label{rem:neg} 
 {\em{1) In spite of appearances, the $\De$-$S_2$ condition does not depend on
$\eR$, but only on $\De$. This is seen in two steps. First,
if ${\mathscr S}$ is a second residual complex, with associated
codimension function $\De$, then ${\mathscr S}=\eR\otimes{{\mathcal L}}$,
with ${{\mathcal L}}$ an invertible sheaf on $X$. This means condition (b)
above does not depend on $\eR$. Second, the map $\ssf(\eM)$ is independent
of $\eR$, for it has the property that any map $\eM\to \eF$ with 
$\eF\in\cm(X,\De)$ factors uniquely through 
$\ssf(\eM)$ (cf.~\cite[p.\,109,\,9.3.6(i)]{lns}).

2) If $\eM$ satisfies (a) and (b) of \Dref{def:s2}, then the $\D(X)$-map $\ssf(\M)$ 
can be uniquely realized as a $\C(X)$-map. Indeed, since 
$\Hr^i_x\eM=0$ for $i<0$, $\Ed{\De}\eM$
has no non-zero components in negative degrees. We should point out that in this case,
the map $\ssf(\eM)$ has a concrete realization, namely it is given by the natural localization
map $\eM\to \oplus_x i_x(\eM_x)$, as $x$ ranges over points with $\De$ value $0$. This can be
seen by base changing to the scheme $\Spec{(\co_{X,x})}$ where $\De(x)=0$, and using the fact
that $\ssf(\eM)$ behaves well with respect to localizations.

3) From the definition of $\ssf(\M)$ we see that any $\eM$ satisfying
(a) and (b) of \Dref{def:s2} is $\De$-$S_2$ iff the natural map
$\eM \to H^0((H^0(\eM'))')$, induced by applying $H^0 \circ \cDt$ to 
the natural truncation map $H^0(\eM') \to \eM'$, is an isomorphism. 

%3) Let $A$ be a local ring  such that $X=\Spec{(A)}$ possesses a dualizing complex (assumed
%residual, to keep with commutative algebraic conventions).
%Let $M$ be a finite non-zero  $A$-module
%and $\eM$ the corresponding coherent $\co_X$-module. 
%Let $\De$ be a codimension function on $X$ and suppose $D^\bullet$ is the dualizing complex
% which has $\De$ as its associated codimension function (i.e.,  if $m=\De(\fm_{{}_A})$, then
%$D^m$ is the injective hull of the $A$-module $k(\fm)=A/\fm_{{}_A}$). According to \Pref{prop:deltaS2}
%in the Appendix, $\eM$ is $\De$-$S_2$ if and only if the following conditions are satisfied:
%\begin{itemize}
%\item $M$ is $S_2$; and
%\item ${\text{ht}}_M(\fp)=\De(\fp)$ for $\fp\in{\text{Supp}}_A(M)$.
%\end{itemize}
%Then, using 
%\Dref{def:s2} and the fact that $M$ is $S_2$ if and only if it is the zeroth cohomology of a Cousin
%complex (see statement above \Dref{def:s2}), it  is easy to see that
%$\eM$ is $\De$-$S_2$ if and only if the following conditions are satisfied:
%\begin{itemize}
%\item $M$ is $S_2$; and
%\item ${\text{ht}}_M(\fp)=\De(\fp)$ for $\fp\in{\text{Supp}}_A(M)$; and
%\item $\Hom_A(M,\,Z^0(D^\bullet))\neq 0$.
%\end{itemize}
%The stronger conditions imposed on $S_2(\De)$-modules yields greater symmetries between
%$M$ and $\Hom_A(M,Z^0(D^\bullet))$ (see \Rsref{rem:Db-T-Kw}). 
}}
\end{rems}
We now give the connection between $\De$-$S_2$ and $S_2$ modules. 
If $A$ is a noetherian ring, we transplant notations and concepts for quasi-coherent sheaves on 
$X=\Spec{(A)}$ to
modules on $A$ in an obvious way, and the notations are self-explanatory. 
Thus if $\De$ is a codimension function on $X$, and $M$ is
an $A$-module, then $\Ed{\De}(M)$ means the complex of global sections of $\Ed{\De}(\eM)$ where
$\eM$ is the quasi-coherent $\co_X$-module corresponding to $M$. 
There are certain commutative algebra conventions we will follow whenever we move to that mode.
Complexes of $A$-modules will have  ``bullets" as superscripts. A dualizing complex
$D^\bullet$ of $A$-modules is a complex whose corresponding complex of quasi-coherent
$\co_X$-modules is residual. 
If $D^\bullet$ is a dualizing complex, then the associated codimension function $\De$ is
defined by the property that $\Hr^i_{\fp}(D^\bullet_\fp)=0$
for $i\ne\De(\fp)$ and $\Hr_{\fp}^{\De(\fp)}(D^\bullet_\fp)=J_A(A/\fp)$, where $J_A(A/\fp)$ 
is the injective hull of the $A$-module $A/\fp$, and $D^\bullet_\fp$
means the localization $D^\bullet$ at $\fp$. Note that, by the conventions of commutative algebra
$D^\bullet$ is Cousin with respect to its associated codimension function.
Finally if $\De$ is a codimension function on $\Spec{(A)}$, and
$C^\bullet$ a Cousin complex of $A$-modules with respect to $\De$ 
(i.e., $C^\bullet=\Ed{\De}(C^\bullet)$), then
\[ C^\bullet(\fp)\set H^{\De(\fp)}(\Gamma_\fp C^\bullet_\fp) \qquad \qquad (\fp\in\Spec{(A)}). \]

\begin{lem}\label{lem:deltaS2}
Let $A$ be a noetherian ring such that $\Spec{(A)}$ has a codimension function $\De$
and let $M$ be a non-zero finitely generated $A$-module. Then
 ${\rm{ht}}_M=\De\vert_{{\rm{Supp}}(M)}$ if and only if 
 ${\rm{Min}}(M) = \{\fp\in{\rm{Supp}}(M)\,\vert\, \De(\fp)=0\}$.
  \end{lem}
 
 \proof This is obvious. 
 \qed
 
 \begin{prop}\label{prop:deltaS2} Let $A$ be a noetherian ring with a dualizing complex $D^\bullet$ 
 whose associated codimension function is $\De$. Then $M\in S_2(\De)$
 if and only if $M$ is $S_2$ and ${\rm{ht}}_M=\De\vert_{{\rm{Supp}}(M)}$.
\end{prop}

\proof In what follows $\eH=(H_i)$ is the filtration on $\Spec{(A)}$ given by 
$H_i=\{\fp\,\vert\, {\rm{ht}}_M(\fp)\ge i\}$.
Suppose $M\in S_2(\De)$. The map $\ssf= \ssf(M)\colon M \to \Ed{\De}(M)$ of \eqref{eq:SG} is
then given by $M\to \oplus_{\De(\fp)=0}M_\fp$ (see 2) of \Rsref{rem:neg}), and $H^0(\ssf)$ is an isomorphism. It follows that
${\rm{Min}}(M) = \{\fp\in {\rm{Supp}}(M)\,\vert\,\De(\fp)=0\}$. 
Therefore by \Lref{lem:deltaS2}, 
${\rm{ht}}_M=\De\vert_{{\rm{Supp}}(M)}$, whence the Cousin complex 
$\Ed{\eH}(M)$ of $M$ with
respect to the filtration $\eH$ of $M$, is the Cousin complex $\Ed{\De}(M)$. Since $H^0(\ssf)$ is an
isomorphism, it follows that $M$ is $S_2$ by \cite[p.\,516,\,Example\,4.4]{sh5}. 

Conversely, suppose $M$ is $S_2$ and ${\rm{ht}}_M=\De\vert_{{\rm{Supp}}(M)}$. Then $\Ed{\eH}(M)=
\Ed{\De}(M)$, and the natural map $M\to \Ed{\De}(M)$ (i.e., the one induced by 
$M\to \oplus_{\De(\fp)=0}M_\fp$) yields an isomorphism on applying $H^0$. It therefore only remains
to show that $H^0(\Homb_A(M,\,D^\bullet))\neq 0$. Now note that $\Homb_A(M,\,D^\bullet)$ has no
negative terms. Indeed, if $\De(\fp)<0$, then $\fp\notin {\rm{Supp}}(M)$, and therefore we have
$\Hom_A(M,\,D^\bullet(\fp))=\Hom_{A_\fp}(M_\fp,\,D^\bullet(\fp))=0$. 
Also note that for any $\fp$ such that $\De(\fp)=0$ 
(e.g., by \Lref{lem:deltaS2}, if $\fp\in{\rm{Min}}(M)$), 
$D^\bullet_\fp$ has no positive terms and $D^0_\fp = D^\bullet(\fp)$.
It follows that for any such $\fp$ 
%if $\De(\fp)=0$ (e.g., by \Lref{lem:deltaS2}, if $\fp\in{\rm{Min}}(M)$),
we have,
\[ \Homb_{A}(M,\,D^\bullet)_\fp = 
\Homb_{A_\fp}(M_\fp,\,D^\bullet(\fp)).\]
The last complex is a one term complex, concentrated at degree $0$, and equal in that degree
to the $A_{\fp}$-module 
$\Hom_{A_\fp}(M_\fp,\,D^\bullet(\fp))$. Therefore, $H^0(\Homb_A(M,\,D^\bullet))_\fp=\Hom_{A_\fp}(M_\fp,\,D^\bullet(\fp))$,  i.e., to the Matlis dual of the finitely generated $A_\fp$-module $M_\fp$, which is
non-zero whenever $\fp\in{\rm{Min}}(M)$. 
\qed

\begin{rem}\label{rem:deltaS2}{\em{Let $A$ be a local ring  such that $X=\Spec{(A)}$ possesses a dualizing complex (assumed residual, to keep with commutative algebraic conventions) and let 
$d=\dim{A}$. Suppose $A$ has a canonical
module $K$, i.e., a module such that ${\widehat{K}}=\Hom_A(\Hr_\fm^d(A), J_A(k(\fm))$, where
$J_A(k(\fm)$ is the injective hull of the $A$-module $k(\fm)$.
Let $\De_{{}_F}$ be the {\em fundamental
codimension function} on $X$, i.e., $\De_{{}_F}(\fp)=d-\dim{(A/\fp)}$. Let $D^\bullet$ be a dualizing
complex with codimension function $\De_{{}_F}$. It is well known that $K\cong H^0(D^\bullet)$.
For  $M\in\{A, K\}$ one checks that $M$ is $S_2$ if and only if 
$M$ is $(\De_{{}_F})$-$S_2$. Indeed if $A$ is $S_2$, $A$ is equidimensional, whence
${\text{ht}}_A=\De_{{}_F}$ (see \cite[1.9]{ayma}). If $K$ is $S_2$, then
$K$ is a submodule of $D^0$, whence ${\rm{Min}}(K)\subset \{\fp\,\vert\,\De_{{}_F}(\fp)=0\}$. On the
other hand, $K_\fp$ is a canonical module for $A_\fp$ for every $\fp\in\Spec{(A)}$, whence
${\rm{Supp}}(K) = \Spec{(A)}$. It follows that ${\rm{Min}}(K)= \{\fp\,\vert\,\De_{{}_F}(\fp)=0\}$, and from
this it is easy to check that ${\text{ht}}_K=\De_{{}_F}$.
}}
\end{rem}

\begin{lem}\label{lem:neg} Let $\eM$ be $S_2$ on $(X,\De)$ and $\eR$
a residual complex with $\De$ as its associated codimension function.
Then the Cousin complex $\eM'\in\Cozt{\De}(X)$ has no non-vanishing
terms in negative degrees.
\end{lem}

\proof Note that the first non-vanishing homology of a Cousin complex $\eC$ (with respect
to $\De$) appears at the first non-zero term of $\eC$. More precisely,
\begin{equation*}\tag{{$*$}}
\min\{n\,\vert\,H^n(\eC)\neq 0\} = \min\{p\,\vert\,\eC^p\neq 0\}.
\end{equation*}
To see this, first let $\eC(x)$ ($x\in X$) be as in \cite[\S\,3.2,\,p.\,37,\,para.\,6]{lns}, i.e., $\eC(x)$ is an
${\widehat{\co}}_{X,x}$-module isomorphic to $\Hr^{\De(x)}_x(\eC)$ and $\eC^p=
\oplus_{\De(x)=p}i_x\eC(x)$. Equality ($*$) translates to:
\[\min\{n\,\vert\,H^n(\eC)\neq 0\}=\min\{\De(x)\,\vert\,\eC(x)\neq 0\}.\]
Now suppose $x$ is a point with least $\De$-value satisfying $\eC(x)\neq 0$. Then $\eC_x$ is
a one-term complex (the term being $\eC(x)$), whence
\[H^{\De(x)}(\eC)_x \cong H^{\De(x)}(\eC_x) = \eC(x)\neq 0\]
proving ($*$). If $\eM$ is $S_2$ with respect to $\De$, then applying ($*$) to $\eC=\eM'$, we get the
result. \qed

%For $x\in X$, let $\eM'(x)$ be as in 
%\cite[\S\S\,3.2,\,p.\,37,\,para.\,6]{lns}, i.e.,
%$\eM'(x)$ is an ${\widehat{\co}}_{X,x}$-module canonically isomorphic
%to $\Hr^{\De(x)}_x(\eM')$ and $(\eM')^p=\oplus_{\De(x)=p}i_x\eM'(x)$.
%Suppose $x\in X$ with $\De(x)<0$. It is enough to show that
%$\eM'(x)=0$. Consider the natural flat map in $\rbbF$: 
%$f\colon X'\set \Spec{\co_{X,x}}\to X$, and let $\eC=f^*\eM'$. Then
%$\eC$ is Cousin with respect to $\De'\set \sh{f}\De$, given by
%$\De'(y)=\De(f(y))$ (cf.~\cite[p.\,14,\,2.1.2]{lns}). Moreover, since all
%points in $X'$ have a negative $\De'$ value, therefore $H^i\eC=0$
%when $i\ge 0$. On the other hand, since $H^i\eC = f^*H^i(\eM')$,
%$H^i\eC=0$ when $i<0$. Being a Cousin complex all of whose cohomology sheaves
%are zero, $\eC\cong 0$ in $\D(\X)$ and this means $\eC=0$.
%But $\eC(y)=\eM'(f(y))$ for all $y\in X'$, whence $\eM'(x)=0$.\qed
%
\begin{defi}Let $\eM$ be $S_2$ on $(X,\De)$. We say $\eM$ is
{\emph{Cohen-Macaulay up to degree $m$ on}} $(X,\De)$ (or
$\De$-CM up to degree $m$) if $\ssf(\eM)_x\colon \eM_x\to \Ed{\De}(\eM)_x$
is a quasi-isomorphism for every $x\in X$ with $\De(x)\le m$. The full
subcategory of $S_2(\De)$  of modules which are $\De$-CM up 
to degree $m$ will be denoted $\cm(\De)_{\le\!m}$.
\end{defi}

We remark that $\cm(\De)_{\le 2} = S_2(\De)$. This follows from the 
fact that for $p>0$, the support of $H^p(\Ed{\De}(\eM))$ has 
codimension at least $p+2$. At the other extreme, we use 
$\cm(\De)_{\le\!\infty}$ to denote the subcategory of all 
$\De$-Cohen-Macaulay modules.

We are in a position to state and prove the first of our main theorems, 
namely, \Tref{thm:main1}. We wish to make a few orienting remarks in order to
understand the Theorem's relationship to the results of Dibaei and Tousi 
\cite[p.\,19,Thm.\,1.4]{dt} and of Kawasaki \cite[Thm.\,4.4]{kw} (see \Rsref{rem:Db-T-Kw}).
Fix a residual complex $\eR$ on $(X,\,\De)$.
Let $\eM\in S_2(\De)$ and $\eN\set(\Ed{\De}\eM)^*$. The Theorem
is concerned with  certain symmetric relations between $\eM$
and $\eN$. The first assertion is that $\eN \in S_2(\De)$. According
to the Theorem, stripped of its category theoretic language,
the relations between $\eM$ and $\eN$ are as follows 
(where we write equalities for functorial isomorphisms to reduce clutter):
\begin{enumerate}
\item[(i)] $\Ed{\De}(\eN) =\eM'$; $\Ed{\De}(\eM)=\eN'$.
\item[(ii)] $\eN=H^0(\eM')$; $\eM=H^0(\eN')$.
\item[(iii)] $\eM=\Ed{\De}(\eN)^*$ (note $\eN\set \Ed{\De}(\eM)^*$).
\item[(iv)]  The following are equivalent:(a) $\eM \in\cm(\De)_{\le m}$, 
(b) $\eN\in\cm(\De)_{\le m}$, (c) $H^0(\eM')\in\cm(\De)_{\le m}$, 
and (d) $H^0(\eN')\in\cm(\De)_{\le m}$.
\end{enumerate}
If $X=\Spec{(A)}$,  and $S_2$ in the usual
sense, then, in \Rsref{rem:Db-T-Kw},  we  say more on the connections with the 
just cited results of Dibaei, Tousi \cite{dt}
and Kawasaki \cite{kw}. We would also like to draw the reader's attention to
\cite[p.\,110,\,9.3.7]{lns}.

Finally some orienting remarks concerning the notations used in the theorem. Let
$\Ed{\De}^*$ denote the composite ${\boldsymbol{-}}^*\circ \Ed{\De}\colon \Dcs(\X) \to \Ac(\X)$,
i.e., $\Ed{\De}^*(\eF) = (\Ed{\De}(\eF))^*$ for $\eF\in\Dcs(\X)$. One may regard $\Ac(\X)$, whence
$S_2(\De)$, as a subcategory of $\Dcs(\X)$. By ``restricting" $\Ed{\De}^*$ to $S_2(\De)$ we get
a functor
\[T\colon S_2(\De)\to \Ac(\X)\]
given (at the level of objects) by $\eM\mapsto (\Ed{\De}\eM)^*=\Ed{\De}^*(\eM)$. Similarly
we have the functor $H^0\circ{\boldsymbol{(-')}}\colon\Ac(\X)\to \Ac(\X)$. Restricting to $S_2(\De)$ 
we get
\[U\colon S_2(\De)\to \Ac(\X)\]
given at the level of objects by $\eM\mapsto H^0(\eM')$. 

The theorem asserts (among other things) that $T$ (resp. $U$) takes values in $S_2(\De)$. In
such a case, by standard ``extension" and ``restriction" dualities, if $i\colon S_2(\De)\to \Ac(\X)$
is the inclusion functor, we get endofunctors  ${\mathbb{T}}$ and ${\mathbb{U}}$ on $S_2(\De)$
defined uniquely by $i\circ{\mathbb{T}}=T$ and $i\circ{\mathbb{U}}=U$.

\begin{thm}\label{thm:main1} Let $\eR$ be a residual complex on $(X,\,\De)$
and let ${\boldsymbol{{-}^*}}$ and ${\boldsymbol{{-}'}}$ be computed with
respect to $\eR$. Let $i\colon S_2(\De) \to \Ac(\X)$ be the natural
embedding.
\begin{enumerate}
\item[(a)] The contravariant functors $T\colon S_2(\De)\to \Ac(\X)$ and $U\colon S_2(\De)\to \Ac(\X)$ 
take values in $S_2(\De)$.
\item[(b)] Let ${\mathbb{T}}\colon S_2(\De) \to S_2(\De)$ and 
${\mathbb{U}}\colon S_2(\De) \to S_2(\De)$ be the contravariant functors
defined by $i\circ {\mathbb{T}} = T$ and $i\circ{\mathbb{U}}=U$. Then
\begin{equation}\label{iso:T-S-1}
{\mathbb{T}} \iso {\mathbb{U}}
\end{equation}
or, equivalently,
\begin{equation}\label{iso:T-S-2}
T \iso U.
\end{equation}
\item[(c)] The contravariant functor ${\mathbb{T}}$ (and therefore
${\mathbb{U}}$) is an anti-equivalence of categories and is its
own pseudo-inverse, i.e.~,
\begin{equation}\label{iso:Tsq}
{\mathbb{T}}^2 \cong {\boldsymbol{1}} \cong {\mathbb{U}}^2
\end{equation}
\item[(d)] There is a functorial isomorphism
\begin{equation}\label{iso:ET}
\Ed{\De}\circ T \iso {\boldsymbol{{-}'}}|_{S_2(\De)}
\end{equation}
such that the following diagram commutes:
\[
\xymatrix{
{\boldsymbol{({-}')^*}}|_{S_2(\De)} \ar[rr]^-{\Iso}_{\phantom{XXX}\eqref{eq:M}}
& & i \\
(\Ed{\De}T)^* \ar[u]^{{\eqref{iso:ET}}^*}_\wr \ar@{=}[rr] & & i{\mathbb{T}}^2
\ar[u]^{\wr}_{\eqref{iso:Tsq}}
}
\]
(Note that therefore \eqref{iso:ET} and \eqref{iso:Tsq} determine each other.)
\item[(e)] ${\mathbb{T}}\eM\in\cm(\De)_{\le\!m}\Longleftrightarrow
\eM\in\cm(\De)_{\le\!m}\Longleftrightarrow {\mathbb{U}}\eM\in\cm(\De)_{\le\!m}$.
\end{enumerate}
\end{thm}

\begin{rems}\label{rem:Db-T-Kw} {\em{For the more commutative algebraic minded readers,
here are the re-interpretations in terms of rings and complexes of modules. As before,
the notations used are the obvious transplants of the notations we have used for schemes and
quasi-coherent sheaves, and are self-explanatory.
In what follows
$A$ is a local ring of dimension $d$,
possessing a dualizing complex $D^\bullet$, and $M\neq 0$ is an $A$-module..
Let $\De\colon \Spec{(A)}\to {\mathbb{Z}}$
be the associated codimension function.   We remind the
reader that it is standard in commutative algebra to assume that $D^\bullet$ 
is Cousin with respect to $\De$ (or, equivalently, $D^\bullet$ is residual).
This implies 
\[D^i=\bigoplus_{\substack{\fp\in\Spec{(A)},\\ \De(\fp)=i}}J_A(A/\fp),\] 
where,  $J_A(A/\fp)$ is the injective hull of the $A$-module $A/\fp$.
The presence of a  codimension function $\De$ on $\Spec{(A)}$ implies the presence of a 
{\em fundamental codimension function} $\De_{{}F}$ given by
$(\fp)\mapsto d-\dim{A/\fp}$. If the codimension function of the $D^\bullet$ is $\De_{{}F}$, then
$D^\bullet$ is called a {\em fundamental dualizing complex}
\cite[p.120,\,1.1]{dib}. We denote the fundamental dualizing complex $D^\bullet_{{}_F}$.

%\begin{enumerate}
 
1) Suppose $M$ is $S_2$ with respect to $\De$ (not necessarily equal to $\De_{{}_F}$).
 \Tref{thm:main1} asserts the existence of another finitely generated module, also $S_2$
 with respect to $\De$, in some fundamental sense, {\em dual} to $M$. This module can be
 described in two ways (and the fact that two descriptions are the ``same" is one of the statements
 of \Tref{thm:main1}). 
 The first way is as follows:  Consider the complex $T^\bullet=\Homb_A(\Ed{\De}(M),\,D^\bullet)$ and set $T(M)=H^0(T^\bullet)$. The second way is to consider 
$U^\bullet=\Homb_A(M,\,D^\bullet)$ and set
 $U(M)=H^0(U^\bullet)$. Then $T(M)$ and $U(M)$ are (functorially) isomorphic. In what follows,
 call the common module $N$.
 
2) The module $M$ can be recovered from $N$ as either $T(N)$ or as $U(N)$. This establishes
 a ``symmetry" between $M$ and $N$.
 
3) By \cite[p.109,\,Prop.\,9.3.5]{lns} the Cousin complexes $\Ed{\De}(M)$ and $\Ed{\De}(N)$ 
have finitely generated cohomology modules. (See also \cite[Thm.\,3.2]{dt}.) In fact
$\Ed{\De}(M)$ can be identified with $\Homb_A(N,\,D^\bullet)$ and $\Ed{\De}(N)$ can
be identified with the complex $\Homb_A(M,\,D^\bullet)$.
 
4)The above can be summarized by the diagram \eqref{diag:wk} which we reproduce
(reminding the reader that ${\rm{coz}}^2_\De$ is the essential image of $S_2(\De)$ under
$\Ed{\De}$):
\begin{equation*}
\xymatrix{
S_2(\Delta) \ar[rr]^{\textup{dualize}} 
\ar[d]^{\downarrow E_{\Delta}}
 && \textup{coz}^2_{\Delta} \ar[ll] \ar[d]^{\downarrow H^0} \\
\textup{coz}^2_{\Delta} \ar[u]^{H^0 \uparrow} 
\ar[rr]_{\textup{dualize}} 
 && S_2(\Delta) \ar[ll] \ar[u]^{E_{\Delta} \uparrow}
}
\end{equation*}
Note that from the diagram, completing a clockwise circuit starting from either the 
northwest corner, or the 
the southeast corner amounts to saying $U(U(M))\cong M$ for $M\in S_2(\De)$;
and completing a counter-clockwise circuit starting from the same two vertices, 
amounts to saying $T(T(M))\cong M$ for $M\in S_2(\De)$. If we imagine $M$ as an object
in the northwest corner, then its ``dual" $N$ occurs in the southeast corner (by following the
transformations of $M$ along any of the routes possible). At the southwest corner, we have
the transform of  $M$ along the immediate the south pointing arrow  and the transform of $N$ 
along the immediate west pointing direction. The ``equality" of these transforms amounts to
the ``equality" $\Ed{\De}(M)=\Homb_A(N,\,D^\bullet)$. We leave it to the reader to use the diagram
to work out other possible relations between $M$, $N$, $\Ed{\De}(M)$, and $\Ed{\De}(N)$.

5) The relationship between \Tref{thm:main1} and \cite[p.\,19,\,Thm.\,1.4]{dt} is complicated
 and needs to be explored in greater detail than we do in this
paper. Here is what we can say. Let $k=\min\{j\,|\, {\text{Supp}}_A(M)\cap{\text{Ass}}_A(D_{{}_F}^j)\neq 0\}$, and let $\De\!=\!\De_{{}_M}\!=\!\De_{{}_F}\!-\!k$. Define $U^\bullet$ and $U(M)$ as we did in 1), 
but without assuming $M$ is $\De$-$S_2$.
%By our choice of $\De$ it is not hard to see that $\Hom_A(M,D^j)=0$ for $j< 0$. 
The re-interpretation of \cite[Thm.\,1.4]{dt} in our language is then as follows.
In {\it loc.\,cit.} it is shown that there is a natural (and unique) map of complexes 
$\Ed{\De}(U(M))\to \Homb(M,D^\bullet)$ which lifts the identity on $U(M)$. Consider the
condition
\begin{equation*}\tag{\dag}
{\rm{Min}}(M)={\rm{Ass}}(U(M)).
\end{equation*}
Condition ($\dag$) is readily seen to be equivalent to our preferred condition 
${\rm{ht}}_M=\De\vert_{{\rm{Supp}}(M)}$ (see {\em{loc.\,cit.}}\,(ii)). Thus, according to
\Pref{prop:deltaS2}, $M$ is $\De$-$S_2$ if and only if $M$ is $S_2$ and ($\dag$) holds.
Under the assumption that $M$ is $S_2$ and ($\dag$) holds, 
%It is easy to check that $\De(\fp)=
%{\text{ht}}_M(\fp)$ for all $\fp\in{\text{Supp}}(M)$ if and only if ${\text{Min}}_A(M)={\text{Ass}}_A(S(M))$.
%In particular the Cousin complex of $M$ with respect to the height filtration $\eH$ of $M$ agrees
%with the Cousin complex of $M$ with respect to the codimension function $\De$ in this case,
%whence $M$ is $S_2$ if and only if the natural map $M\to \Ed{\De}(M)$ yields an isomorphism on
%taking the zeroth cohomology. Further, the equality ${\text{Min}}_A(M)={\text{Ass}}_A(S(M))$ yields
%$S(M)\neq 0$, i.e., $H^0(\Homb_A(M,D^\bullet))\neq 0$Thus in this situation (when ${\text{Min}}_A(M)
%={\text{Ass}}_A(S(M))$) we have $M$ is $S_2$ if and only if it is $\De$-$S_2$, since 
%$H^j(\Homb_A(M,D^\bullet)=0$ for $j<0$.  
the identification
 \[\Ed{\De}(U(M))=\Homb_A(M,D^\bullet)\] 
 of 3) above is established by different methods in
 \cite[p.\,19,Thm.\,1.4]{dt}(iv).
 However, the symmetric relationship between $M$ and $U(M)$ (= $K_M^I$ in the
 notation of {\it{loc.\,cit.}}) is not fully explored, nor the fact that
 $U(M)$ is also $T(M)$. It should be pointed out that there are other results (concerning $M$ and
 $U(M)$) in {\it{loc.\,cit.}} that are of considerable interest (see also \cite[Thm.\,4.4]{kw}).
 
6) As noted in \Rref{rem:deltaS2}, $A$ (resp.~$K$) is $S_2$ if and only if it is $\De_{{}_F}$-$S_2$.
Since $M=A$ (resp.~$M=K$) yields $N=K$ (resp.~$N=A$), it follows that $A$ is $S_2$ if and
only if $K$ is $S_2$. In this case $\Ed{\De_{{}_F}}=D^\bullet_{{}_F}$ and $\Ed{\De_{{}_F}}(A)=
\Homb_A(K,\,D^\bullet_{{}_F})$ (cf. \cite[p.110,\,Example\,9.3.7]{lns} and 
 \cite[p.126,\,Thm.\,4.6]{dib}). This observation has consequences later, when we generalize
 \cite[p.\,125,\,Thm.\,3.3]{dib} (see \Rref{rem:Db-gor}).
 
 7) As a matter of record we point out the well known fact that
 if $A$ is $S_2$ then $A$ has no embedded primes and is equidimensional. In
 particular, the fundamental codimension function in this case coincides with 
the height function on $\Spec{A}$, $\fp\mapsto {\text{ht}}(\fp)$ (see, e.g., \cite[p.23,\,Rmk.\,2.1]{dt}). 
%\end{enumerate}
}}
\end{rems}

%\begin{thm}\label{thm:main1} Let $\eM$ be $S_2$ on $(X,\De)$ and $\eR$
%a residual complex with $\De_\eR=\De$. Let $\eN=(\Ed{\De}\eM)^*$.
%Then
%\begin{enumerate}
%\item[(i)] There is a functorial isomorphism 
%$\Ed{\Delta}(\eN)\iso \eM'$.
%\item[(ii)] $\eN$ is $S_2$ on $(X,\De)$. In particular, from (i),
%there is a functorial isomorphism of coherent sheaves $\eN\iso H^0(\eM')$.
%\item[(iii)] There is a functorial isomorphism $\Ed{\De}(\eM)\iso \eN'$.
%In particular, from (ii), there is a functorial isomorphism of coherent
%sheaves $\eM\iso H^0(\eN')$.
%\item[(iv)] For $m\ge 0$, $\eM$ is $\De$-CM up to degree $m$ if and only if
%$\eN$ is $\De$-CM up to degree $m$.
%\end{enumerate}
%\end{thm}
%
%\begin{rems}{\emph{
%\begin{enumerate}
%\item From \cite[p.\,109,\,Prop.\,9.3.5]{lns} one sees that the complex
%$E\eM=\Ed{\De}\eM$ has coherent cohomology 
%(since $\Dc=\Dcs=\Dt$ on ordinary 
%schemes). This means that $\eN$ is defined and coherent.
%\item In view of 2) of Rsref\,\ref{rem:neg}, we will always regard
%$\ssf(\eM)$ as a map in $\C(X)$ whenever $\eM$ is $\De$-$S_2$.
%\end{enumerate}
%}}
%\end{rems}

%Before we give a proof of \ref{thm:main1}, we need a lemma.
%\begin{lem}\label{lem:s2dual}
%For any $\eM \in S_2(\De)$, we have $H^0(\eM') \in S_2(\De)$.
%\end{lem}

\noindent{\em Proof of \ref{thm:main1}.} 
Let $\eM \in S_2(\De)$.
By Proposition \ref{prop:theta}, there is a natural isomorphism
$E_{\Delta}\eM \cong (H^0(\eM'))'$ and hence by dualizing we obtain 
an isomorphism
\[
T\eM = (E_{\Delta}\eM)^* \cong H^0(\eM') = U(\eM)
\]
which we take as our choice for \eqref{iso:T-S-2}.
(Recall that since we are on an ordinary scheme, 
we have $\Dc^* = \Dc$.)

%Note that $\eM'$ is a Cousin complex. 
Next we claim that the natural map $U(\eM) = H^0(\eM') \to \eM'$
satisfies the universal property of $\ssf(U(\eM))$ as 
mentioned in 1) of \Rsref{rem:neg}. Suppose 
$\eF$ is Cohen-Macaulay and $\alpha \colon H^0\eM' \to \eF$ a map 
of complexes. Upon dualizing we get the following 
natural maps.
\[
\begin{CD}
\eM @. \cDt\eF \\
@V{\ssf(\eM)}VV  @VV{\quad\cDt(\alpha)\quad}V\\
E\eM @<{\;\;\quad\cong\quad\;\;}<< (U(\eM))'
\end{CD}
\]
The bottom row is an isomorphism of Cousin complexes while the top
row consists of objects in $\Ac$. Therefore, if $\eM$ is $\De$-$S_2$, so that 
$\eM = H^0E\eM$, then there exists a unique map 
$\cDt\eF \to \eM$ that makes the above diagram commute.
Thus, dualizing back we obtain our claim. 

It follows from the universal property of $\ssf(U(\eM))$ that 
there is a unique isomorphism $E(U(\eM)) \iso \eM'$ over $U(\eM)$.
This shows that $U(\eM)$ is also $\De$-$S_2$ and moreover 
via $U \iso T$ we also get
a choice for the isomorphism \eqref{iso:ET} of
part\,(d). (Note that $E$ of a module has no negative terms
and hence $\eM'$ satisfies \ref{def:s2}(b).)
Thus we obtain parts (a)--(d) of the theorem.

It remains to prove (e). In view of part (c), 
it suffices to prove any one of the implications in (e).
We shall prove that $\eM \in \cm(\De)_{\le m}
\Longrightarrow \mathbb T \eM \in \cm(\De)_{\le m}$.

We first remark that if 
$0 \ne \eN \in \cm(\De)_{\le \infty}$, i.e., $\eN$ is Cohen-Macaulay, then  
so is $\mathbb U \eN \cong \mathbb T \eN$. 
Indeed, by \ref{prop:Ac-CM}, the Cousin complex $\eN'$ has only one 
non-vanishing homology and hence $\mathbb U \eN$ is $\D$-isomorphic 
to $\eN'$ which is Cohen-Macaulay.

Let $\eM \in \Ac$. Pick a point $x\in X$ and set 
$X'_x=\Spec{\co_{X,x}}$. Let $f_x \colon X'_x \to X$
be the canonical flat map. The constructions 
of $E$ and $\ssf$ behave well with base-change to $X_x$.
Since Cohen-Macaulayness of a module is equivalent
to $\ssf$ being a quasi-isomorphism, 
one checks that 
$\eM \in \cm(X,\De)_{\le m}$
iff for all $x \in X$ such that $\De(x)\le m$, 
it holds that $f_x^*\eM$ is Cohen-Macaulay on $X_x$.
Also note that $f_x^*$ ``commutes'' with the functors
$\ssf$ and $\mathbb T$ and for the corresponding notions
on $X_x$ we use the symbols $\ssf_x$, $\mathbb T_x$
respectively.

Thus we have
\begin{align*}
\eM \in \cm(X, \De)_{\le m} &\Longleftrightarrow
f_x^*\eM \textup{ is CM on }X_x \textup{ for all } x \in X, \De(x)\le m \\
&\Longleftrightarrow 
\mathbb T_xf_x^*\eM \textup{ is CM on }X_x 
\textup{ for all }x \in X, \De(x)\le m \\
&\Longleftrightarrow 
f_x^*\mathbb T\eM \textup{ is CM on }X_x 
\textup{ for all }x \in X, \De(x)\le m \\
&\Longleftrightarrow 
\mathbb T\eM \in \cm(X, \De)_{\le m}. 
\end{align*}

\qed

\section{Connections with Grothendieck duality}

In this section we use \Pref{prop:Ac-CM} (or, equivalently, 
\Pref{prop:Ac-Coz}) to extend and make transparent some of the results
in \cite{dcc}, e.g.~the result that a map $f$ is flat if and only
if $f^!$ transforms Cohen-Macaulay complexes to appropriate Cohen-Macaulay
complexes (cf.~\cite[Theorems\,7.2.2 and 9.3.12]{dcc}).

\subsection{Pull back of Cousin complexes} We fix, for the rest of this 
discussion, a map $f\colon (\X,\De')\to (\Y,\,\De)$ in $\rbbFc$ and a residual
complex $\eR$ on $(\Y,\,\De)$. Let  $\Cozs{\De}$ (resp. $\Cozs{\De'}$) represents the category
of all Cousin complexes in $\Dqct(\Y)$ (resp. $\Dqct(\X)$). This category is larger than
$\Cozt{\De}$, the category Cousin complexes in $\Dcs(\Y)$. Let 
\stepcounter{thm}
\begin{equation}\label{eq:sharp}\tag{\thethm}
\sh{f}\colon \Cozs{\De}(\Y) \to \Cozs{\De'}(\X)
\end{equation}
be the functor constructed in \cite[p.\,10,\,Main Theorem]{lns}. Very briefly, if
$\sh{f}\De$ is the codimension function in \eqref{eq:cod-f},  with $f$ is smooth and $d$
 its relative dimension{\footnote{In affine terms, this means, if $(A,I)\to (B,J)$ is a pseudo-finite
 map of
 adic rings, then the map is formally smooth---in the sense of lifting idempotents---and
 for any prime ideal $\fp\subset A$, and $\fq\in\Spec{(B)}$ lying over $\fp$, the integer $d$ is given by
 $d= \dim{B_\fq/\fp B}+{\text{tr.deg.}}_{k(\fp)k(\fq)}$ \cite[p.22,\,Def.\,2.4.2 and p.28,\, Def.\,2.6.2]{lns}.}}, 
 then $\sh{f}(\eF)$ is isomorphic to the Cousin complex
$\Ed{\sh{f}\De}(\R\iGp{\X}(\bL f^*\eF\otimes_{\co_\X}\Omega^d_f[d]))$. If $f$ is a closed immersion, then
$\sh{f}(\eF)$ is isomorphic to the $\sh{f}\De$-Cousin complex $\sHomb(\co_\X,\,\eF)$. A general
$f\in \rbbFc$, can locally be written as a closed immersion followed by a smooth map, giving an idea
how one might construct $\sh{f}\eF$ in this case. One has to show the results are independent of
factorizations (into closed immersions followed by smooth maps), and have 2-functorial properties (i.e.
$\sh{(gf)} \cong \sh{f}\sh{g}$ plus ``associativity"), and this is what is done in \cite{lns}. In \cite{dcc} it is
shown that $\sh{f}(\eF)$ is a concrete approximation of the twisted inverse image $f^!\eF$ of
Alonso, Jerem\'{\i}as and Lipman \cite{dfs}. In the event $f$ is pseudo-proper, $\sh{f}(\eF)$ represents
the functor $\eG\mapsto\Hom_{\C(\Y)}(f_*\eG,\,\eF)$ on $\Cozs{\De'}$ (see \cite[p.185,\,Thm.\,8.1.10]
{dcc}), and if $f$ is an open immersion then $\sh{f}=f^*$.

For $\eF$ in
$\Cozt{\De}(\Y)$, define
\stepcounter{thm}
\begin{equation}\label{def:(sharp)}\tag{\thethm}
\psh{f}_\eR(\eF) \set \sHomb_\X(f^*\eF^*,\,\sh{f}\eR) = (f^*\eF^*)'
\end{equation}
where ``upper star" is with respect to $\eR$ and ${\boldsymbol{-'}}$ 
is with respect
to $\sh{f}\eR$. Since $\eF^*$ is a coherent $\co_\Y$-module by
\Pref{prop:Ac-Coz}, $\psh{f}_\eR\eF$ is in $\Cozt{\De'}(\X)$. Thus we
have a functor
\[
\psh{f}_\eR\colon \Cozt{\De}(\Y) \to \Cozt{\De'}(\X).
\]
The functor $\psh{f}_\eR$ makes transparent many of the relationships
established between the twisted inverse image functor $f^!$ and $\sh{f}$
in \cite{dcc}. We will show in \Tref{thm:sh-(sh)} that $\psh{f}_\eR$
is essentially $\sh{f}|_{\Cozt{}}$. But first, we would like to show
that $\psh{f}_\eR$ is independent of $\eR$.

\begin{prop} Let $f$ be as above. There is a family
of isomorphisms
\begin{equation}\label{iso:psiRR'}
\psi_{\eR,\eR'} = \psi_{f,\eR,\eR'}\colon \psh{f}_{\eR'} \iso \psh{f}_\eR,
\end{equation}
one for each pair of residual complexes $\eR, \eR'$ on $(\Y,\De)$
such that $\psi_{\eR,\eR'}\circ\psi_{\eR',\eR''}=\psi_{\eR,\eR''}$ 
(cocycle condition) for any  three residual complexes $\eR$,
$\eR'$, $\eR''$ on $(\Y,\De)$.
\end{prop}

\proof The proof rests on the fact that there are isomorphisms between
$\eR'$ and $\eS\set \eR\otimes {\mathcal{L}}$, where ${\mathcal{L}}$ is the
coherent invertible $\co_\Y$-module $\sHom(\eR,\,\eR')$, and that
isomorphisms between $\eR'$ and $\eS$ are indexed by units in the ring
$\Gamma(\Y,\,\co_\Y)$ (since $\Hom(\eR',\eR')=\Gamma(\Y,\,\co_\Y)$).

We first make the identification
\[ \psh{f}_{\eR} = \psh{f}_\eS \]
via the canonical identifications 
$\sHom(\eF,\eS)=\sHom(\eF,\eR)\otimes{\mathcal{L}}$, 
$\sh{f}\eS=\sh{f}\eR\otimes f^*{\mathcal{L}}$,  and
$\sHomb({\eM}\otimes f^*{\mathcal{L}}, \sh{f}\eR\otimes f^*{\mathcal{L}})
=\sHomb(\eM,\,\sh{f}\eR)$ for a coherent sheaf $\eF$ on $\Y$ and
a coherent sheaf $\eM$ on $\X$.

Next, pick an isomorphism $\alpha\colon \eR'\iso \eS$. Then $\alpha$
induces an isomorphism
\[\psi_\alpha\colon \psh{f}_{\eR'} \iso \psh{f}_\eS (=\psh{f}_\eR). \]
In greater detail, $\psi_\alpha=q_\alpha^{-1}p_\alpha = s_\alpha r_\alpha^{-1}$
where $p_\alpha, q_\alpha, r_\alpha, s_\alpha$ are the maps induced by
$\alpha$ in the commutative diagram below (where $\sHomb = \sHomb_\X$):
\[
\xymatrix{
\sHomb(f^*\sHom_\Y(\eF,\eR'),\sh{f}\eR')
\ar[rr]^{\Iso}_{p_\alpha}
& & 
\sHomb(f^*\sHom_\Y(\eF,\eR'),\sh{f}\eS) \\
\sHomb(f^*\sHom_\Y(\eF,\eS),\sh{f}\eR')
\ar[u]_{\wr}^{r_\alpha}
\ar[rr]_{s_\alpha}
  & &  \
\sHomb(f^*\sHom_\Y(\eF,\eS),\sh{f}\eS) \ar[u]^{\wr}_{q_\alpha}
}
\]
Suppose $\beta\colon \eR'\iso \eS$ is another isomorphism. We claim that
$\psi_\alpha=\psi^{\phantom{X}}_\beta$. Note that there exists a (unique) unit 
$a\in \Gamma(\Y,\,\co_\Y)$ such that $\alpha=a\beta$, so that
$p_\alpha=ap^{\phantom{X}}_\beta$ and $q_\alpha=aq^{}_\beta$. It follows that
$q_\alpha^{-1}p_\alpha=q_\beta^{-1}p^{}_\beta$. This proves the
claim. Setting $\psi_{\eR,\eR'}$ equal to $\psi_\alpha$, it is not
difficult to establish the cocycle rules.
\qed

From the proposition we deduce a well defined functor
\begin{equation}\label{eq:(sharp)}
\psh{f}\colon \Cozt{\De}(\Y) \to \Cozt{\De'}(\X)
\end{equation}
independent of $\eR$, together with isomorphisms
\begin{equation}\label{iso:sigma-R}
\sigma^{\phantom{X}}_\eR\colon \psh{f}_\eR \iso \psh{f}
\end{equation}
such that $\sigma_\eR^{-1}\circ 
\sigma^{\phantom{X}}_{\eR'} = \psi_{\eR,\eR'}$.

\subsection{Grothendieck duality} For $f$ and $\eR$ as above, in
\cite[\S\,9]{dcc}, functors $\pshr{f}_\eR$ and $\pshr{f}$ are
constructed,\footnote{more precisely $\pshr{|f|}_\eR$ and $\pshr{|f|}$
are constructed, where $|f|\colon\X\to \Y$ is the map underlying
$f\colon (\X,\De')\to (\Y,\De)$.} more or less along the lines
that $\psh{f}_\eR$ and from it $\psh{f}$ are constructed. In slightly
greater detail, if $\eF$ is an object in $\Dcs(\Y)\cap \D^+(\Y)$, then
\[\pshr{f}_\eR\eF \set \cDt'\circ\bL f^*\circ \cDt (\eF)\]
where $\cDt$ (resp.~$\cDt'$) is the dualizing functor in \eqref{eq:Dt}
associated to $\eR$ (resp.~$\sh{f}\eR$)\footnote{$\sh{f}\eR$ is also
residual \cite[p.\,105,\,Prop.\,9.1.4]{lns}.}. It is not hard to see
that $\pshr{f}_\eR$ is an object in $\Dcs(\X)\cap\D^+(\X)$ (see
\cite[\S\S\,9.2,\,p.\,187]{dcc}, especially the discussion after
(9.2.1)). The passage from $\pshr{f}_\eR\eF$ to $\pshr{f}\eF$ is
identical to the passage from $\psh{f}_\eR$ to $\psh{f}$, and one
has functorial isomorphisms
\[ \theta_\eR = \theta_{\!f\!,\>\eR}\colon \pshr{f}_\eR \iso \pshr{f} \]
such that, for a second residual complex $\eR'$ on $(\Y,\,\De)$,
\[\phi_{\eR,\eR'} (= \phi_{f,\eR,\eR'}) \set \theta_\eR^{-1}\theta_{\eR'}
\colon \pshr{f}_{\eR'} \iso \pshr{f}_\eR\]
satisfies cocycle rules.

A couple of minor irritants need to be quickly addressed. In
\cite[\S\,9]{dcc}, the source and target of $\pshr{f}_\eR$
and $\pshr{f}$ are complicated subcategories of $\D(\Y)$ and
$\D(\X)$ respectively. For our purposes, it suffices to observe
that the source contains $\Dcs(\Y)\cap \D^+(\Y)$. Thus in this
paper, we regard $\pshr{f}_\eR$ and $\pshr{f}$ as functors with
source $\Dcs(\Y)\cap\D^+(\Y)$ and target $\Dcs(\X)\cap\D^+(\X)$:
\[\pshr{f}_\eR\cong \pshr{f}\colon \Dcs(\Y)\cap\D^+(\Y) \to
\Dcs(\X)\cap\D^+(\X).\]
A second point needs to be made. As in \cite{dcc}, we reserve the
notation $f^!$ (as opposed to $\pshr{f}$) for the twisted inverse
image functor obtained in \cite{dfs} (cf.~[{\it Ibid,}\,p.\,2, Thm.\,2 and
beginning of \S\S\,1.3]) for pseudo-proper maps,
and extended to composites of compactifiable maps in 
\cite[p.\,261,\,7.1.3]{suresh}.
We point out that $f^!$ and $\pshr{f}$ are canonically isomorphic
when both are defined \cite[p.\,190,\,Thm.\,9.3.10]{dcc}.

\subsection{Cousin complexes and duality} Let $f$, $\eR$, $\cDt$,
$\cDt'$ be as in the previous section. As for the symbols
${\boldsymbol{-^*}}$ and ${\boldsymbol{-}'}$, the context will
determine the interpretation (see \Rref{rem:star-prime}). To put
a fine point to it, if $\eG$ is in $\Cozt{\De}(\Y)$, then 
$\eG^*=\sHom(\eG,\eR)$, whereas if $\eG\in \Cozt{\De'}(\X)$,
then $\eG^*=\sHom(\eG,\sh{f}\eR)$. Similarly, ${\eM}'$ is
$\sHomb(\eM,\eR)$ or $\sHomb(\eM,\sh{f}\eR)$ depending on
whether $\eM$ is a coherent $\co_\Y$-module or a coherent
$\co_\X$-module.

We denote by ${\overline{Q}}_\Y$ the localization functor
\[ {\overline{Q}}_\Y\colon \Cozt{\De}(\Y) \to \Dcs(\Y)\cap \D^+(\Y).\]
We would like to understand the effect of duality on Cousin complexes.
In other words, we wish to study the functor
\[ \pshr{f}{\overline{Q}}_\Y\colon \Cozt{\De}(\Y) \to \Dcs(\X)\cap\D^+(\X).\]
In order to describe the above functor more explicitly in terms
of $\eR$, we set
\[\pshr{f}_{[\eR]}\set \cDt'\circ\bL{f^*}\circ Q_\Y(-)^*. \]
By \eqref{iso:*Dt} we have a canonical isomorphism 
$\cDt{\overline{Q}}_\Y\iso Q_\Y(-)^*$ of functors on 
$\Cozt{\De}(\Y)$. This induces a series of isomorphisms

\stepcounter{thm}
\begin{equation*}\label{iso:coz(!)}\tag{\thethm}
\pshr{f}_{[\eR]} \iso \pshr{f}_\eR\circ{\overline{Q}}_\Y \iso 
\pshr{f}\circ{\overline{Q}}_\Y.
\end{equation*}
It is convenient---as we will see---to study $\pshr{f}{\overline{Q}}_\Y$
via $\pshr{f}_{[\eR]}$. The behavior of $\pshr{f}_{[\eR]}$ with respect
to ``change of residual complexes" obviously follows the
behavior of $\pshr{f}_\eR{\overline{Q}}_\Y$ with respect to such
a change. In other words, if $\eR'$ is another residual complex on
$(\Y,\De)$, we have an isomorphism of functors
\[\phi_{[\eR,\eR']}\colon \psh{f}_{[\eR']} \iso \psh{f}_{[\eR]} \]
which is compatible with $\phi_{\eR,\eR'}$ and the first arrow
in \eqref{iso:coz(!)}.

The behavior of $\psh{f}{\overline{Q}}_\Y$ is studied through a
comparison map $\pshr{\gamma}_f\colon {\overline{Q}}_\X\psh{f} \to
\pshr{f}{\overline{Q}}_\Y$ which is a more down to earth version
of the comparison map in \cite[p.\,163,\,(4.1.4.1)]{dcc} when we
restrict our attention to $\Cozt{\De}(\Y)$ (instead of $\Cozs{\De}(\Y))$.
Here is how it is defined. Recall that if $\eM\in \Ac(\X)$, then there
is an obvious functorial map $\bL{f^*}Q_\Y\eM\to Q_\X{f^*}\eM$. This
induces a natural transformation
\stepcounter{thm}
\begin{equation*}\label{eq:gamma*}\tag{\thethm}
\gamma_f^*\colon \bL{f^*}Q_\Y \to Q_\X{f^*}(-)^*
\end{equation*}
between functors on $\Cozt{\De}(\Y)$. Set
\[\pshr{\gamma}_{f,\eR}\set \cDt'\gamma_f^*\colon {\overline{Q}}_\X\psh{f}_\eR
\to \pshr{f}_{[\eR]}.\]
As can be easily checked from the definitions, this map behaves well
with respect to change of residual complexes on $(\Y,\De)$, i.e.
\[\phi_{[\eR,\eR']}\pshr{\gamma}_{f,\eR'}=
\pshr{\gamma}_{f,\eR}{\overline{Q}}_X\psi_{\eR,\eR'}.\]
We therefore have a well-defined comparison map
\stepcounter{thm}
\begin{equation*}\label{eq:(gamma)}\tag{\thethm}
\pshr{\gamma}_f\colon {\overline{Q}}_\X\circ\psh{f} \to \pshr{f}\circ
{\overline{Q}}_\Y.
\end{equation*}

\subsection{Tor-independence} The following definition does not
need $\X$, $\Y$ or $f$ to be in $\rbbF$.

\begin{defi}\label{def:torind}A pair $(f,\,\eM)$, with $f\colon\X\to \Y$
a map of formal schemes and $\eM$ an object of $\Ac(\Y)$, is
said to be a {\it{tor-independent pair}} if the following holds
for every $x\in\X$ (with $y=f(x)$, $A=\co_{\Y,y}$, $B=\co_{\X,x}$
and $M=\eM_y$):
\[ {\mathrm{Tor}}^A_i(B,\,M)=0 \qquad\qquad(i>o).\]
In other words, $(f,\,\eM)$ is tor-independent if and only if the
natural map $\bL{f^*}Q_\Y\to Q_\X{f^*}$ in $\Dc(\X)$ is an isomorphism
on $\eM$:
\[\bL{f^*}Q_\Y\eM \iso Q_\X{f^*}\eM.\]
\end{defi}

\begin{srem}\label{rem:torind} Note that $\eM$ is a flat $\co_\Y$-module
if and only if $(f,\,\eM)$ is tor-independent for every $f$. In fact,
$\eM$ is a flat $\co_\Y$-module if and only if $(f,\eM)$ is tor-independent
for every closed immersion $f\colon\X\to \Y$.
\end{srem}

\begin{lem}\label{lem:torind} Let $f\colon(\X,\,\De')\to (\Y,\,\De)$ be a map in
$\rbbFc$,  $\eF$ an object in $\Cozt{\De}(\Y)$, and $\eR$ a
residual complex on $(\Y,\De)$. For $x\in\X$, let $y=f(x)$, 
$M=(\eF^*)_y$ and $A$, $B$ the local rings at $x$ and $y$. 
Then for every integer $i$
\[\Hr^i_x(\pshr{f}\eF) \cong 
\Hom_B({\mathrm{Tor}}^A_{i-\De'(x)}(B,M),\sh{f}\eR(x)).\]
In particular, $\pshr{f}\eF\in\cm(\X;\De')$ if and only
if $(f,\,\eF^*)$ is a tor-independent pair (since the
right side is the Matlis dual of the finitely generated
$B$-module ${\mathrm{Tor}}^A_{i-\De'(x)}(B,M)$.
\end{lem}

\proof Since $\sh{f}\eR$ is residual, whence injective, we
have by \cite[p.\,33,\,(5.2.1)]{ajl} 
\begin{align*}
\R\Gamma_x\pshr{f}\eF & \cong \Homb_B((\bL{f^*}\eF^*)_x,\,\Gamma_x\eR)\\
& \cong \Homb_B(B{\overset{\bL}{\otimes}}M,\,\eR(x)[-\De'(x)])\\
& \cong \Homb_B(B{\overset{\bL}{\otimes}}M[\De'(x)],\,\eR(x)).
\end{align*}
Since $\eR(x)$ is an injective $B$-module, applying $\Hr^i$ to both
sides, we get the result.\qed

\begin{thm}\label{thm:CMness} Let $f\colon (\X,\De')\to (\Y,\De)$
and $\eF$, $\eR$ be as in the lemma above. The following are equivalent
\begin{enumerate}
\item[(i)] $\pshr{f}\eF$ is Cohen-Macaulay with respect to $\De'$;
\item[(ii)] $(f,\,\eF^*)$ is a tor-independent pair;
\item[(iii)] The map 
\[\gamma^*_f(\eF)\colon\bL{f^*}Q_\Y\eF^* \to Q_\X{f^*}\eF^*\]
of \eqref{eq:gamma*} is an isomorphism;
\vspace{1.5pt}
\item[(iv)] The map 
\[\pshr{\gamma}_f(\eF)\colon {\overline{Q}}_\X\psh{f}\eF \to 
\pshr{f}\overline{Q}_\Y\eF\]
of \eqref{eq:(gamma)} is an isomorphism.
\end{enumerate}
\end{thm}
\proof Evidently (i), (ii) and (iii) are equivalent. Since 
$\pshr{\gamma}_{f,\eR}(\eF)$ is the ``dual" of $\gamma^*_f(\eF)$
with respect to the residual complex $\sh{f}\eR$, clearly (iv) is
equivalent to (iii).\qed

\Tref{thm:CMness} gives us a way of reproving (and allows for a better
understanding of) \cite[p.\,191,\,Thm.\,9.3.12]{dcc} 
(cf. \cite[p.\,182,\,Thm.\,7.2.2]{dcc}). Moreover, coupled with 
\cite[p.\,191,\,Thm.\,9.3.13]{dcc} it allows for subtle twist on
that theorem on Gorenstein complexes. We should point out that
there is a typographical error in {\emph{loc.cit.}}---the hypothesis
on $\F$ should be $\F\in\Cozt{\De}(\Y)$ and not $\F\in\Cozs{\De}(\Y)$.

\begin{thm}\cite[9.3.12 and 7.2.2]{dcc}\label{thm:dcc-flat} 
Let $f$ and $\eR$ be as above. Then the following are equivalent
\begin{enumerate}
\item[(i)] $f$ is flat;
\item[(ii)] $\pshr{f}\eF$ is Cohen-Macaulay with respect to $\De'$ for every
$\eF\in\cm(\Y;\De)$;
\item[(iii)] The map of functors
\[\pshr{\gamma}_f\colon {\overline{Q}}_\X\psh{f} \to 
\pshr{f}{\overline{Q}}_\Y\]
is an isomorphism.
\end{enumerate}
\end{thm}

\proof This follows immediately from \Tref{thm:CMness} and the fact
that $f$ is flat if and only if $\gamma^*_f\colon\bL{f^*}Q_\Y\eF^*
\to Q_\X{f^*}\eF^*$ is an isomorphism for every $\eF\in\Cozt{\De}(\Y)$.
We point out that the essential image of $\Cozt{\De}(\Y)$ under
${\boldsymbol{-'}}$ is $\Ac(\Y)$ according to \Pref{prop:Ac-Coz}.\qed

\smallskip

We now move to examining another statement in \cite{dcc}.
In \cite[p.\,178,\,Thm.\,6.3.1]{dcc} it is shown that the
Cousin of the map $\gamma_f^!$ is an isomorphism. It is
much simpler to prove the analogous statement for 
$\pshr{\gamma}_f$. It is worth pointing out that this analogous statement gains strength
only when content is poured into it by showing that $\pshr{\gamma}_f$ is ``isomorphic" to
$\gamma_f^!$, i.e., by showing Diagram \eqref{diag:zeta-phi} commutes. Note this involves showing that $\psh{f}$ is isomorphic to $\sh{f}$ (when restricted
to Cousin complexes in $\Dcs$) and that
$f^!$ is isomorphic to $\pshr{f}$ (when restricted to $\Dc$). 
The latter is proven in \cite{dcc} using the fact that $\gamma_f^!$ is
an isomorphism for residual complexes (in itself needing a careful examination of $\gamma_f^!$)
and the former needs \Tref{thm:sh-(sh)}. In other words, the ``simpler" proof gains content only when
many more complicated proofs are brought in to set the context. 
Here then is the statement analogous to \cite[p.\,178,\,Thm.\,6.3.1]{dcc}. Let $f\colon (\X,\De')\to
(\Y,\De)$ be a map in $\rbbFc$, and set
\[f^{(E)}\set E_{\De'}\pshr{f}{\overline{Q}}_\Y\colon\Cozt{\De}(\Y)
\to \Cozt{\De'}(\X)\]
and $\gamma_f^{(E)}$ to be the composite
\[\psh{f} \iso E_{\De'}(\psh{f})\xrightarrow{E(\pshr{\gamma}_f)} f^{(E)}.\]
We then have
\begin{prop}\label{prop:gammaE}
The functorial map
\[\gamma_f^{(E)}\colon \psh{f} \to f^{(E)}\]
is an isomorphism.
\end{prop}
\proof Fix a residual complex $\eR$ on $(\Y,\,\De)$. It is enough to
show that the functorial map $E(\pshr{\gamma}_\eR)\colon E(\psh{f}_\eR) \to E(\pshr{f}_\eR)$
is an isomorphism or what amounts to the same thing, that
\[\Hr_x^{\De'(x)}(\pshr{\gamma}_\eR)\colon \Hr_x^{\De'(x)}(\psh{f}_\eR)
\to \Hr_x^{\De'(x)}(\pshr{f}_\eR)\]
is an isomorphism for every $x\in\X$. Fixing such an $x$, we see---as in
the proof of \Lref{lem:torind}---that after taking Matlis duals this
amounts to showing that for $\eF\in\Cozt{\De}(\Y)$, the natural map
\[H^0((\gamma_f^*))_x\colon H^0(\bL{f^*}\eF^*)_x \to
H^0(f^*\eF^*)_x\]
is an isomorphism, which it clearly is. \qed

\section{The pseudofunctor ${\boldsymbol{\psh{-}}}$ vs.~the pseudofunctor 
${\boldsymbol{\sh{-}}}$}

In this section we show that $\psh{f}\eF$ is naturally isomorphic
to $\sh{f}\eF$ when $\eF\in\Cozt{\De}(\Y)$ and $f\colon (\X,\De')\to
(\Y,\De)$ is a map in $\rbbFc$.  But first we wish to understand
the behavior of $\psh{(fg)}$ for a composite of two maps $f$ and
$g$ with respect to $\psh{f}$ and $\psh{g}$.

\subsection{Variance properties} We assume familiarity with the notion
of a {\emph{contravariant pseudofunctor}} defined for example in
\cite[p.\,45]{lns}. Indeed the main focus of \cite{lns} is to construct
$\sh{f}$ for suitable maps $f$ in such a way that the assignments
$(\Y,\,\De)\mapsto \Cozs{\De}(\Y)$ and
$f\mapsto \sh{f}$ define a pseudofunctor $\boldsymbol{\sh{-}}$. It
turns out that the assignments $(\Y,\,\De)\mapsto\Cozt{\De}(\Y)$,
$(\Y,\De)\in\rbbFc$, and $f\mapsto \psh{f}$, $f$ a map in $\rbbFc$,
are pseudofunctorial. To see this, let
\[(\W,\,\De'')\xrightarrow{g} (\X,\,\De')\xrightarrow{f} (\Y,\De) \]
be a pair of maps in $\rbbFc$. Let $\eR$ be a residual complex on
$(\Y,\,\De)$ and $\eS\set \sh{f}\eR$. The pseudofunctor ${\boldsymbol{\sh{-}}}$
gives an isomorphism
\[C_{g,f}^{\sharp}(\eR) \colon \sh{g}\sh{f}\eR \iso \sh{(fg)}\eR.\]
This together with the isomorphisms $f^*\eF^*\iso {(f^*\eF^*)'}^*= 
(\psh{f}\eF)^*$ (cf.~\Rref{rem:star-prime}) gives an isomorphism
\[C^{(\sharp)}_{g,f,\eR}\colon \psh{g}_\eS\psh{f}_\eR \iso \psh{(fg)}_\eR.\]
The process is completely analogous to the one described
\cite[p.\,136,\,(3.3.15)]{conrad} and \cite[p.\,188,\,(9.2.3)]{dcc} for
${\boldsymbol{\pshr{-}}}$. The isomorphism $C^{(\sharp)}_{g,f,\eR}$
behaves well with respect to change of residual complexes, giving
an isomorphism
\[C^{(\sharp)}_{g,f}\colon \psh{g}\psh{f} \iso \psh{(fg)}.\]
Using the pseudofunctoriality of ${\boldsymbol{\sh{-}}}$ it is easy to
see that the above identification is ``associative", and hence defines
a pseudofunctor ${\boldsymbol{\psh{-}}}$ on $\rbbFc$ with
$\psh{(\Y,\De)}=\Cozt{\De}(\Y)$ for $(\Y,\,\De)\in \rbbFc$. Since, as
we briefly noted, the process is identical to the process of constructing
the pseudofunctor ${\boldsymbol{\pshr{-}}}$, with $f^*, g^*$ and $(fg)^*$
replacing $\bL{f^*}, \bL{g^*}$ and $\bL{(fg)}^*$ in the construction
in \cite[p.\,188,\,(9.2.3)]{dcc}, we have the following proposition
(cf. \cite[p.\,163,\,Thm.\,4.1.4(d)]{dcc}):

\begin{prop}\label{prop:(sh)(shr)} With $f,g$ as above, the following
diagram commutes: 
$$
\xymatrix{
{\overline{Q}}_\W\psh{g}\psh{f} \ar[d]_{\pshr{\gamma}_g(\psh{f})} 
\ar[r]^-{C^{(\sharp)}_{g,f}}
& {\overline{Q}}_\W\psh{(fg)} \ar[dd]^{\pshr{\gamma}_{fg}} \\
\pshr{g}{\overline{Q}}_\X\psh{f} \ar[d]_{\pshr{g}(\pshr{\gamma}_f)} & \\
\pshr{g}\pshr{f}{\overline{Q}}_{\Y} \ar[r]^{C^{(!)}_{g,f}} & 
\pshr{(fg)}{\overline{Q}}_{\Y}
}
$$
where the map $C^{(!)}_{g,f}$ is the map in \cite[p.\,188,\,(9.2.3)]{dcc}.
\end{prop}

\subsection{${\boldsymbol{\psh{-}}}$ vs. ${\boldsymbol{\sh{-}}}$} For
$f\colon (\X,\,\De')\to (\Y,\,\De)$ and $\eF\in\Cozt{\De}(\Y)$ define
a map
\[\zeta = \zeta_f(\eF) \colon \sh{f}\eF \to \psh{f}\eF\]
as follows. Pick a residual complex $\eR$ on $(\Y,\,\De)$. Functoriality
of $\sh{f}$ gives a map $\Gamma(\Y,\,\co_\Y)$-modules
$\Hom(\eF,\eR)\to\Hom(\sh{f}\eF,\sh{f}\eR)$ which is well behaved
with respect to Zariski localizations of $\Y$. In other words we have
a map of $\co_\Y$-modules
\[\eF^*=\sHom(\eF,\eR) \to f_*\sHom(\sh{f}\eF,\sh{f}\eR)=f_*((\sh{f}\eF)^*)\]
inducing a map of coherent $\co_\X$-modules
\[\xi = \xi_f(\eF)\colon f^*(\eF^*) \to (\sh{f}\eF)^*.\]
The natural isomorphism $\sh{f}\eF \iso {(\sh{f}\eF)^*}'$ of
\Pref{prop:Ac-Coz} followed by ${\xi}'$ gives us a
map
\[\zeta_\eR\colon \sh{f}\eF \to \psh{f}_\eR\eF\]
which one checks (from the definitions) is independent of $\eR$, i.e.
\[\psi_{[\eR,\eR']}(\eF)\circ\zeta_{\eR'} = \zeta_\eR.\]
We therefore get a well defined map of functors
\stepcounter{thm}
\begin{equation*}\label{map:zeta}\tag{\thethm}
\zeta_f \colon \sh{f}|_{\Cozt{\phantom{.}}(\Y)} \to \psh{f}
\end{equation*}
If $g\colon (W,\,\De'')\to (\X,\,\De')$ is a second map, it is
easy to check from the definitions that the diagram
\stepcounter{thm}
\begin{equation*}\label{diag:zeta}\tag{\thethm}
\xymatrix{
\sh{g}\sh{f}\eF \ar[r]^{\sh{g}\zeta_f} \ar[d]_{\sh{C}_{g,f}}^{\wr} 
& \sh{g}\psh{f}\eF
\ar[r]^{\zeta_g} & \psh{g}\psh{f}\eF \ar[d]_{\wr}^{\psh{C}_{g,f}}\\
\sh{(fg)}\eF \ar[rr]_{\zeta_{fg}} && \psh{(fg)}\eF
}
\end{equation*}
commutes for every $\eF\in\Cozt{\De}(\Y)$.

\subsection{Traces} Let $f\colon (\X,\,\De')\to (\Y,\,\De)$ be a 
pseudo-proper map in $\rbbFc$. According to \cite[p.\,146,\,(2.2.4)]{dcc}
and \cite[p.\,156,\,Thm.\,2.4.2(b)]{dcc}, for every $\eF\in\Cozs{\De}(\Y)$
we have a trace map
\[\Tr{f}(\eF)\colon f_*\sh{f}\eF \to \eF.\]
If $\eF\in \Cozt{\De}(\Y)\subset \Cozs{\De}(\Y)$, then we
define, as a counterpart to $\Tr{f}$, 
\stepcounter{thm}
\begin{equation*}\label{map:Tr-f}\tag{\thethm}
\sTr{f}(\eF)\colon f_*\psh{f}\eF \to \eF
\end{equation*}
as follows. First  pick a residual complex $\eR$ on $(\Y,\,\De)$
and define $\sTr{f,\eR}(\eF)$ as the map which 
makes the following diagram commute (see
also \cite[p.\,189,\,(9.3.5)]{dcc}):
\[
\xymatrix{
f_*\psh{f}_\eR\eF \ar[d]_{\sTr{f,\eR}(\eF)} \ar@{=}[r] &
f_*\sHomb(f^*\eF^*,\sh{f}\eR) \ar[r]^{\Iso} &
\sHomb(\eF^*,\,f_*\sh{f}\eR) \ar[d]^{\Tr{f}(\eR)} \\
\eF \ar[r]^{\Iso} & {\eF^*}' \ar@{=}[r] &
\sHomb(\eF^*,\,\eR)
}
\]
As usual, one checks that this definition is independent of
$\eR$, i.e.~we have a relation 
$\sTr{f,\eR}f_*\psi_{[\eR,\eR']}=\sTr{f,\eR'}$. This
gives \eqref{map:Tr-f}.

We had, just before the above definition, fleetingly drawn the
reader's attention to the trace map in \cite[p.\,189,\,(9.3.5)]{dcc}
\[\ttr{f}^r \colon \R{f_*}\pshr{f}\to {\boldsymbol{1}}.\] 
The point is that the definition of $\sTr{f}$ is almost identical to the
definition of $\ttr{f}^r$, provided we replace $f^*$ by $\bL{f^*}$,
and this gives part\,(iii) of the \Pref{prop:trace} below. Part\,(i)
is immediate from the analogous \cite[p.\,156,\,Thm.\,2.4.2(b)]{dcc} and
part\,(ii) is immediate from the definition of $\sTr{f}$.

\begin{prop}\label{prop:trace} Let $f\colon (\X,\,\De')\to (\Y,\,\De)$
be a pseudo-proper map in $\rbbFc$ and $\eF\in \Cozt{\De}(\Y)$.
\begin{itemize}
\item[(i)] If $g\colon (W,\,\De'')\to (\X,\,\De')$ is a second 
pseudo-proper map then the diagram
\[
\xymatrix{
(fg)_*\psh{g}\psh{f}\eF \ar@{=}[d] \ar[rr]^{\Iso}_{\psh{C}_{g,f}}
& & (fg)_*\psh{(fg)}\eF \ar[dd]^{\sTr{fg}} \\
f_*g_*\psh{g}\psh{f}\eF \ar[d]_{f_*\sTr{g}} & & \\
f_*\psh{f}\eF \ar[rr]_{\sTr{f}} & & \eF
}
\]
commutes (see \cite[p.\,156,\,Thm.\,2.4.2(b)]{dcc}).
\item[(ii)] The diagram 
\[
\xymatrix{
f_*\sh{f}\eF \ar[r]^{f_*\zeta_f} \ar[dr]_{\Tr{f}}
& f_*\psh{f}\eF \ar[d]^{\sTr{f}}\\
& \eF}
\]
commutes. 
\item[(iii)] The diagram (in which we suppress localization
functors like ${\overline{Q}}_\Y$ to avoid clutter)
\[
\xymatrix{
f_*\psh{f}\eF \ar[d]_{\sTr{f}} \ar[r]^{\Iso} & \R{f_*}\psh{f}\eF 
\ar[d]^{\R{f_*}\pshr{\gamma}_f} \\
\eF & \R{f_*}\pshr{f}\eF \ar[l]^{\ttr{f}^r}
}
\]
commutes in $\Dcs(\Y)\cap\D^+(\Y)$.
\end{itemize}
\end{prop}

If $f$ and $\eF$ are as in the Proposition and
\[\Phi_f(\eF) \colon \pshr{f}\eF \iso f^!\eF\]
is the isomorphism in \cite[p.\,190,\,Thm.\,9.3.10]{dcc} then
by \Pref{prop:trace}(ii) and (iii) and the universal properties
of $(f^!,\ttr{f})$ and $(\pshr{f},\ttr{f}^r)$ the following
diagram 
\begin{equation}\label{diag:zeta-phi}
\xymatrix{
\sh{f}\eF \ar[d]_{\ga{f}} \ar[r]^{\zeta_f} & \psh{f}\eF 
\ar[d]^{\pshr{\gamma}_f}\\
\fs\eF & \pshr{f}\eF \ar[l]^{\Phi_f}_{\Iso}
}
\end{equation}
commutes in $\Dcs(\X)\cap\D^+(\X)$, where $\ga{f}$ is the map in 
\cite[p.\,163,\,(4.1.4.1)]{dcc}.

Here is how we compare ${\boldsymbol{\sh{-}}}$ and ${\boldsymbol{\psh{-}}}$.
Recall that a compactifiable map is a map that can be written as
an open immersion followed by a pseudo-proper map.
\begin{thm}\label{thm:sh-(sh)} Let $f\colon (\X,\,\De') \to (\Y,\De)$ be 
a map in $\rbbFc$.
\begin{itemize}
\item[(i)] The map $\zeta_f\colon \sh{f}|_{\Cozt{\De}(\Y)} \to \psh{f}$
is an isomorphism of functors.
\item[(ii)] Diagram \eqref{diag:zeta-phi} continues to commute under the
weaker hypothesis that $f$ is a composite of compactifiable maps.
\end{itemize}
\end{thm}
\proof We first prove (ii). By \eqref{diag:zeta}, \Pref{prop:(sh)(shr)},
\cite[p.\,163,\,Thm.\,4.1.4(d)]{dcc} and [{\emph{Ibid}},\,p.\,190,\,(9.3.10.1)]
the maps $\zeta_f$, $\pshr{\gamma}_f$, $\ga{f}$ and $\Phi_f$ behave
well with respect to composition of maps. Therefore it is enough
to prove that \eqref{diag:zeta-phi} commutes when $f$ is pseudo-proper
and when $f$ is an open immersion. We have already argued that
the diagram commutes when $f$ is pseudo-proper. If $f$ is an open immersion,
all vertices in the diagram can be identified with $f^*\eF$ and all
arrows with the identity map, and hence we are done.

Part\,(i) is equivalent to showing that $\Ed{\De'}(\zeta_f)$ is an isomorphism.
Moreover the question is local on $\X$, and therefore we may assume 
that $f$ is a composite of compactifiable maps. We have proven that in this
case \eqref{diag:zeta-phi} commutes. Applying $\Ed{\De'}$ to this diagram,
and using the fact that $\Ed{\De'}(\ga{f})$ and $\Ed{\De'}(\pshr{\gamma}_f)$
are isomorphisms by \cite[p.\,178\,Thm.\,6.3.1]{dcc} and \Pref{prop:gammaE},
we are done.\qed
\vspace{1.5 pt}

One consequence of \Tref{thm:sh-(sh)} is that every $\C(\X)$\nobreakdash-map
$\sh{f}\eF\to \sh{f}\eR$ is induced by a $\C(\Y)$\nobreakdash-map 
$\eF\to \eR$. More precisely, we have:

\begin{cor} Let $\eF$  be an object in $\Cozt{\De}(\Y)$
and $\eR$ a residual complex on $(\Y,\,\De)$. The natural map 
\[\xi_f \colon f^*\sHom_\Y(\eF,\,\eR) \to \sHom_\X(\sh{f}\eF,\,\sh{f}\eR)\]
is an isomorphism.
\end{cor}

\proof By construction of $\zeta_f$, $\xi_f$ is the dual (with respect
to $\eR$) of $\zeta_f$, which we have shown is an isomorphism.\qed

\section{Gorenstein complexes}\label{s:gorenstein}

In this section we drop the assumption that our formal schemes contain c-dualizing complexes.
In \Ssref{ss:lam-gam} we work with general noetherian schemes.
Starting from \Ssref{ss:gorenstein} we restrict ourselves to schemes and maps in $\bbF$ (which, recall from \Ssref{ss:codim} is the category whose objects are noetherian formal schemes which are universally catenary and admit codimension functions). 

\subsection{The homology localization functor and the derived torsion functor}\label{ss:lam-gam}
The purpose of this subsection is to recall (and draw out) the various ways in which the
``homology localization functor" of \cite{dfs} and \cite{ajl} interacts with the derived torsion functor,
as summarized in \cite[pp.\,69--70,\,Remarks\,6.3.1(1)]{dfs}. In greater detail, let $\X$ be a formal
 scheme (not necessarily in $\bbF$, but as always, noetherian).
The homology localization functor $\BL_{\X}\colon\D(\X)\to\D(\X)$ is defined as
\[\BL_{\X}\set \R\sHomb(\R\iGp{\X}\co_\X,\,\boldsymbol{-})\]
where $\iGp{\X}$ is the torsion functor defined in \Ssref{ss:cat}. Note that the source and the target
of $\R\iGp{\X}$ is $\D(\X)$. For typographical convenience we write $\BL=\BL_{\X}$
and $\BG=\R\iGp{\X}$. By \cite[p.\,54,\,5.2.10.1]{dfs} we have $\BL$ is right adjoint to $\BG$. More precisely, 
we have a bifunctorial isomorphism 
\[\Hom_{\D(\X)}(\eE,\,\BL\eF) \iso \Hom_{\D(\X)}(\BG\eE,\,\eF)\qquad (\eE,\eF\in \D(\X)). \]
We are more interested in the sheafified version of the above, obtained by noting that $\BG$ and
$\BL$ are compatible with Zariski localization as is the above bifunctorial map (see also
the comment preceding \cite[p.\,24,\,(2.5.0.1)]{dfs}), namely:
\stepcounter{thm}
\begin{equation}\label{eq:lam-gam-1}\tag{\thethm}
\R\sHomb(\eE,\,\BL\eF) \iso \R\sHomb(\BG\eE,\,\eF) \qquad (\eE,\eF\in\D(\X)).
\end{equation}
The above also needs the fact that the adjointness between $\BG$ and $\BL$ is $\delta$-functorial,
i.e., it behaves well with translations.
There are three results (other than \eqref{eq:lam-gam-1})) that are important for us in this paper:
\begin{enumerate}
\item[(a)] According to \cite[p.\,68,\,Prop.\,6.2.1]{dfs}, a form of the Greenlees-May duality, we have
\stepcounter{thm}
\begin{equation}\label{eq:lam-F}\tag{\thethm}
\BL\eF \iso \eF \qquad (\eF\in \Dc(\X)).
\end{equation}
\item[(b)] For $\eE\in \D(\X)$ and $\eF\in \Dc(\X)$ we have
\stepcounter{thm}
\begin{equation}\label{eq:lam-gam-2}\tag{\thethm}
\R\sHomb(\eE,\,\eF) \iso \R\sHomb(\BG\eE,\,\BG\eF).
\end{equation}
\item[(c)] The functors $\BL$ and $\BG$ induce quasi-inverse equivalences between the categories
$\Dcs(\X)$ and $\Dc(\X)$. In other words we have (with $\D\,{\hat{}}\,(\X)$ the essential image of
$\D(\X)$ under $\BL$):
\begin{align*}
\eE \in \Dcs(\X) & \Longleftrightarrow \BL\eE\in\Dc(\X)\, {\text{and}}\, \eE\in\Dt(\X) ,\\
\eF\in \Dc(\X) & \Longleftrightarrow \BG\eF\in\Dcs(\X)\, {\text{and}}\, \eF\in\D\,{\hat{}}\,(\X)\\
\end{align*}
and 
\begin{align*}
\eE & \iso \BG\BL\eE \,\qquad (\eE\in \Dcs(\X))\\
\eF & \iso \BL\BG\eF  \qquad (\eF\in \Dc(\X))
\end{align*}
\end{enumerate}
Statement (b) is proven as follows. By \cite[pp.\.69--70,\,Rmks\,6.3.1\,(1)(c)]{dfs}, we have
\stepcounter{thm}
\begin{equation}\label{eq:lam-gam-3}\tag{\thethm}
\BL\BG \iso \BL.
\end{equation}
The relations \eqref{eq:lam-F}, \eqref{eq:lam-gam-3} and \eqref{eq:lam-gam-1} then give---for
$\eE\in\D(\X)$ and $\eF\in\Dc(\X)$---the sequence of isomorphisms
\begin{align*}
\R\sHomb(\eE,\,\eF) & \iso \R\sHomb(\eE,\,\BL\eF) \\
& \iso \R\sHomb(\eE,\,\BL\BG\eF)\\
& \iso \R\sHomb(\BG\eE,\,\BG\eF)
\end{align*}
thus establishing \eqref{eq:lam-gam-2}. We should point out that the map underlying the isomorphism
\eqref{eq:lam-gam-2} is the one induced
by the functor $\BG$. This can be seen by unravelling the definitions of the maps and isomorphisms
in \cite[pp.\,69--70,\,Rmks\,6.3.1\,(1)]{dfs}.

Statement (c) follows from the fact that $\BL$ and $\BG$ induce quasi-inverse equivalences between
$\Dt(\X)$ and $\D\,{\hat{}}\,(\X)$ (see the last line of \cite[Remarks\,6.3.1\,(1)]{dfs}).

\begin{rem}\label{rem:lam-gam}{\em{From \cite[Remarks\,6.3.1\,(1)]{dfs} it is apparent that
 $\BG$ and $\BL$ share
a symmetric relationship (e.g. $\BG^2 \cong \BG$, $\BL^2\cong \BL$, $\BG\BL\cong \BG$ and
$\BL \cong \BL\BG$). It is not hard to see that one can obtain an isomorphism ``dual" to
\eqref{eq:lam-gam-2} given by
\[\R\sHomb(\eE,\,\eF) \iso \R\sHomb(\BL\eE,\,\BL\eF) \qquad (\eE\in\Dcs(\X), \eF\in\D(\X)).\]
Again, the map underlying this isomorphism is the obvious one induced by the functor $\BL$.}}
\end{rem}

\subsection{c-Gorenstein and t-Gorenstein complexes}\label{ss:gorenstein} In this subsection we 
assume that our schemes and maps of schemes are in $\bbF$. 

\begin{defi}\label{def:gorenstein}Let $(\X,\,\De)\in\bbFc$.
A complex $\eF\in\D(\X)$ is said to be {\emph{t-Gorenstein}} with
respect to $\De$ if 
\begin{enumerate}
\item[(a)] $\eF\in\Dqct(\X)$;
\item[(b)] $\eF$ is Cohen-Macaulay with respect to $\De$; and 
\item[(c)] the Cousin complex of $\eF$ with respect to $\De$, $\Ed{\De}\eF$, 
consists of injective objects in $\Aqct(\X)$. 
\end{enumerate}
A complex $\eG\in \Dc(\X)$ is said to be {\emph{c-Gorenstein}} with respect to $\De$ if
$\BG\eG$ is t-Gorenstein.
\end{defi}

%\newpage
 \begin{rems}{\emph{
 \begin{enumerate}
 \item For ordinary schemes of finite Krull dimension, 
 the notions of c-Gorenstein and t-Gorenstein coincide, and in this
 situation we call such complexes simply Gorenstein.
 A t-dualizing complex is t-Gorenstein and a c-dualizing complex is c-Gorenstein (cf. \Dref{def:dualizing} and \Eref{ex:dualizing}).
 \item What is called simply {\emph{Gorenstein}} in \cite[p.\,179]{dcc} is the same as what we have called a t-Gorenstein complex in this paper.
 \item  If $\X$ has finite Krull dimension, a t-Gorenstein complex,
by this definition, is necessarily in $\Dqct^b(\X)$.
\item  We have the following two relations:\\
Let $\eF\in\Dcs(\X)$. Then $\eF$ is t-Gorenstein $\Longleftrightarrow  \BL\eF$ is c-Gorenstein.\\
Let $\eG\in \Dc(\X)$. Then $\eG$ is  c-Gorenstein  $\Longleftrightarrow \BG\eG$ is t-Gorenstein. 
\end{enumerate}
}}
\end{rems}

Our immediate aim is to prove that if $\eG_1$ and $\eG_2$ are c-Gorenstein and $\X$ is
of finite Krull dimension, then the complex
$\R\sHomb(\eG_1,\,\eG_2)$ is $\D(\X)$ isomorphic to a locally free $\co_\X$-module of
finite rank, i.e. \Pref{prop:gorenstein}. We need a preliminary discussion before we can do
this. To that end, let $\X$ be a (noetherian) formal scheme. The torsion functor $\iGp{\X}$ defined in
\Ssref{ss:cat} is a special case of the functor $\iG{\I}\colon \A(\X)\to \A(\X)$ for any 
$\co_\X$ ideal $\I$ (not necessarily an ideal of definition of the formal scheme $\X$).
This is defined, as in \cite[p.\,6,\,1.2.1]{dfs}, by the formula
\[
\iG{\I} \set \dirlm{n}\sHom_{\co_\X}(\co_\X/\I^n,\,{\boldsymbol{-}}).
\]
Note that if $\I$ is an ideal of definition for $\X$ then $\iG{\I}=\iGp{\X}$.

For $x\in\X$, let $k(x)$ denote
the residue field of the local ring $A\set\co_{\X,x}$, and let $E(x)$ denote the injective hull
(as an $A$-module) of $k(x)$. It is well known that $E(x)$ is also an $\widehat{A}$-module,
$\widehat{A}$ being the completion of $A$ with respect to the maximal ideal $\fm_A$ of $A$.
Moreover $E(x)$ is also the injective hull of $k(x)$ viewed as an ${\widehat{A}}$-module.
If $x$ is a {\emph{closed point}} of $\X$ we denote the ideal sheaf of $\{x\}$  by 
$\fm_x$. Let
\[\kappa\colon \widehat{\X}\to \X \]
be the completion of $\X$ along $\{x\}$, then, denoting the unique point of $\widehat{\X}$ by
${\widehat{x}}$, we have $\kappa^*i_xE(x) = i_{\hat{x}}E(x)$. We denote the common
$\co_{\widehat{\X}}$-module $\eE(x)$. For $\eG\in\D(\X)$ we have, by 
\cite[p.50,\,Proposition\,5.2.4]{dfs},  canonical functorial isomorphisms 
(each map deducible from the other via the adjoint pair $(\kappa^*,\kappa_*)$)
\stepcounter{thm}
\begin{align*}\label{iso:kappa-m}\tag{\thethm}
\kappa^*\R\iG{\fm_x}\eG & \iso \R\iGp{\widehat{\X}}\kappa^*\eG \\
\R\iG{\fm_x}\eG & \iso \kappa_*\R\iGp{\widehat{\X}}\kappa^*\eG.
\end{align*}
We point out that $\kappa$ being  a flat map of locally ringed spaces, $\kappa^*=\bL\kappa^*$.
%The adjoint relation between $\kappa^*$ and $\kappa_*$ (see \cite[p.\,147,\,Prop.\,6.7]{sp})
%gives us a bifunctorial map 
%\stepcounter{thm}
%\begin{equation*}\label{iso:kappa*}\tag{\thethm}
%\kappa^*\R\sHomb_{\co_\X}(\eG_1,\,\eG_2) \longrightarrow 
%\R\sHomb_{\co_{\widehat{\X}}}(\kappa^*\eG_1,\,\kappa^*\eG_2) \qquad (\eG_1,\eG_2\in\D(\X))
%\end{equation*}
%induced by the composite
%\[ \R\sHomb_{\co_\X}(\eG_1,\,\eG_2) \longrightarrow 
%\R\sHomb_{\co_\X}(\eG_1,\,\kappa_*\kappa^*\eG_2)
%\xleftarrow{\Iso} \kappa_*\R\sHomb_{\co_{\widehat{\X}}}(\kappa^*\eG_1,\,\kappa^*\eG_2).\]
Suppose now that $\De$ is a codimension function on $\X$ such that $(\X,\De)\in\bbFc$, and suppose
further that $\eG$ is c-Gorenstein with respect to $\De$. The complex $\eG$ being c-Gorenstein,
$\R\iG{\fm_x}\eG$ is isomorphic
to a direct sum of copies of the complex $i_xE(x)[-\De(x)]$, whence from \eqref{iso:kappa-m}.
%
%$\R\iGp{\widehat{\X}}\kappa^*\eG\in\Dcs({\widehat{\X}})$, and since it is isomorphic to 
%, it is isomorphic to a direct sum of copies of $\eE(x)[-\De(x)]$.
%
% This
%direct sum is necessarily finite since its ``Matlis dual" is coherent. More precisely, $\eE(x)$ is
%t-dualizing on ${\widehat{\X}}$, whence $\R\sHomb(\R\iGp{\widehat{\X}}\kappa^*\eG,\eE(x))$ is
%in $\Dc({\widehat{\X}})$, which
%From \eqref{iso:kappa-m} we therefore get
\stepcounter{thm}
\begin{equation*}\label{iso:kappa-E}\tag{\thethm}
\R\iGp{\widehat{\X}}\kappa^*\eG \cong \oplus_{i\in I} \eE(x)[-\De(x)].
\end{equation*}
where the index $i$ varies over a finite index set $I$.
To see the finiteness of $I$, note that $\eE(x)[-\De(x)]$ is t-dualizing on ${\widehat{\X}}$, 
whence from \cite[p.\,28,\,Prop.\,2.5.8\,(a)]{lns},
$\R\sHomb(\R\iGp{\widehat{\X}}\kappa^*\eG,\eE(x)[-\De(x)])\in \Dc({\widehat{\X}})$.  This amounts to
saying that the ${\widehat{A}}$-module $\prod_{i\in I}{\widehat{A}}$ is finitely generated, forcing $I$ to
be finite. We are now in a position to prove:

\begin{prop}\label{prop:gorenstein} Let $\X\in\bbF$ be of finite Krull dimension and $\De$ a 
codimension function on $\X$. 
\begin{enumerate}
\item[(a)] If $\eG_1$, $\eG_2$ are c-Gorenstein on $(\X,\De)$ then $\R\sHomb(\eG_1,\,\eG_2)$ is
$\D(\X)$-isomorphic to a locally free $\co_\X$-module of finite rank.
\item[(b)] If $\eF_1, \eF_2\in \Dcs(\X)$ are t-Gorenstein on $(\X,\,\De)$ then $\R\sHomb(\eF_1,\,\eF_2)$ is
$\D(\X)$-isomorphic to a locally free $\co_\X$-module of finite rank.
\item[(c)] Let $\eF\in\Dcs(\X)$ be t-Gorenstein, and for each $x\in\X$, let $r(x)$ be the number of
copies in $\Hr^{\De(x)}_x(\eF)$ of the injective hull of the residue field $k(x)$ (regarded 
as a $\co_{\X,x}$-module). Then
$r(x)$ is constant on connected components of $\X$.
\end{enumerate}
\end{prop} 

\proof Using the functors $\BL_\X$ and $\iGp{\X}$, especially  
\eqref{eq:lam-gam-2} 
and \Rref{rem:lam-gam}, we see that statements (a) and (b) are equivalent to each other. We will
therefore only prove (a).  Let $x\in\X$ be an arbitrary closed point, and as before, let
$\kappa\colon {\widehat{\X}}\to \X$ be the corresponding completion map. We will show below in 
\Lref{lem:kappa*} that since $\X$ is of finite Krull dimension (and since $\eG_1$ and
$\eG_2$ are c-Gorenstein), we 
have a natural isomorphism
\[\kappa^*\R\sHomb_{\co_\X}(\eG_1,\,\eG_2) \iso 
\R\sHomb_{\co_{\widehat{\X}}}(\kappa^*\eG_1,\,\kappa^*\eG_2).\]
From \eqref{iso:kappa-E} we  get isomorphisms $\R\iGp{\widehat{\X}}\kappa^*\eG_1\cong \oplus_{i\in I}\eE(x)[-\De(x)]$
and $\R\iGp{\widehat{\X}}\kappa^*\eG_2 \cong \oplus_{j\in J}\eE(x)[-\De(x)]$, where the index sets
$I$ and $J$ are finite. We thus have a series of isomorphisms:
\begin{align*}
\kappa^*\R\sHomb_{\co_\X}(\eG_1,\,\eG_2) & 
\xrightarrow[\phantom{\eqref{eq:lam-gam-2}}]{\Iso}
\R\sHomb_{\co_{\widehat{\X}}}(\kappa^*\eG_1,\,\kappa^*\eG_2)\\
& \xrightarrow[\eqref{eq:lam-gam-2}]{\Iso} 
\R\sHomb_{\co_{\widehat{\X}}}(\R\iGp{\widehat{\X}}\kappa^*\eG_1,\,\R\iGp{\widehat{\X}}\kappa^*\eG_2) \\
& \quad\cong\quad \R\sHomb_{\co_{\widehat{\X}}}(\oplus_{i\in I}\eE(x),\, \oplus_{j\in J}\eE(x))\\
& \quad\cong\quad \co_{\widehat{\X}}^{|I|\cdot |J|}
\end{align*}
Now, by \Lref{lem:kappa*}\,(a), $\eG_1$ and $\eG_2$ are in $\Dc^b(\X)$, whence 
$\R\sHomb(\eG_1,\,\eG_2)$ has coherent cohomology. Part (a) follows since 
$x\in\X$ was an arbitrary closed point.
Part (c) is an immediate consequence of the proof just given. Indeed $r(x)^2$ is the rank of
the locally free sheaf $H^0(\sHomb(\eF,\,\eF))$ on the connected component containing $x$.
\qed

\medskip

As a corollary we get the following theorem which
contains \cite[p.\,125,\,Thm.\,3.3]{dib}. (See \Rref{rem:Db-gor}.)

\begin{thm}\label{thm:gorenstein} Let $(\X,\,\De)\in\rbbFc$, i.e. $(\X,\De)\in\bbF$ and $\X$ possesses a 
c-dualizing
complex (equivalently $\X$ possesses a t-dualizing complex which is in $\Dcs(\X)$). Then
\begin{enumerate}
\item[(a)] $\eG$ is c-Gorenstein on $(\X,\,\De)$ if and only if there is a $\D(\X)$-isomorphism
\[\eG \cong \eD\otimes {\mathcal{V}} \]
where $\eD$ is a c-dualizing complex on $\X$ with associated codimension function $\De$ and
${\mathcal{V}}$ is a locally free $\co_\X$-module of finite rank.
\item[(b)] $\eF\in \Dcs$ is t-Gorenstein on $\X$ if and only if there is a $\D(\X)$-isomorphism
\[\eF \cong \eR\otimes {\mathcal{V}} \]
where $\eR$ is a t-dualizing complex on $\X$ with associated codimension function $\De$ and
${\mathcal{V}}$ is a locally free $\co_\X$-module of finite rank.
\end{enumerate}
\end{thm}

\begin{rem}\label{rem:Db-gor}{\em{This remark is intended towards commutative algebraists, 
and we adhere to
the terminology there. In particular, on a local ring, a complex of modules is called dualizing if it is residual. Let us assume, for this remark, that $A$ is a local ring possessing a dualizing complex,
$\De_{{}_F}$ is the associated fundamental codimension function and let $M$ be a finitely generated
$A$-module. 
In \cite[p.\,125,\,Thm.\,3.3]{dib}, Dibaei proves the following (with the notion
of $S_2$, unless otherwise stated, being the usual notion of $S_2$).:\\
{\emph{Suppose $A$ satisfies $S_2$ and suppose that it possesses a dualizing complex. 
Assume $M$ satisfies $S_2$ and 
$(0:_AM)=0$ and let $E(M)$ denote the Cousin complex $E(M)$ of 
$M$ with respect to the fundamental codimension function $\De_{{}_F}$ on $\Spec{(A)}$. 
Then $E(M)$ is an injective complex if and
only if $M$ is isomorphic to a direct sum of a finite number of copies of the canonical module
$K$ of the ring $A$.}}\\
Let us prove this using \Tref{thm:gorenstein}. In view of 6) of \Rsref{rem:Db-T-Kw}, 
evidently if $A$ is $S_2$ and $M$ is a direct sum of
a finite number of copies of $K$ (without the assumption that $M$ is $S_2$, as noticed by Dibaei
in his proof) we must have $M$ is also $S_2$ and $E(M)$ is an injective
complex. Conversely (without the assumption that $A$ is $S_2$, again noticed by Dibaei in {\em 
loc.cit.}), if $M$ is $S_2$ and $E(M)$
is an injective complex, then by \Tref{thm:gorenstein}, we have $E(M)$ is a direct sum of a finite
number of copies of the fundamental dualizing complex. Taking the zeroth cohomology we deduce
that $M$ is a direct sum of a finite number of copies of $K$. Note that as a consequence $A$ must
be $S_2$. Indeed $E(K)$, being a direct summand of $E(M)$, is such that $K=H^0(E(K))$. This
means $K$ is $S_2$ . It follows that $\Homb_A(K,E(K)) \neq 0$. Thus $K$ is $\De_{{}_F}$-$S_2$.
In view of \Tref{thm:main1} this forces $A$ to be $S_2$. As the reader must have noticed, $E(K)$ is
the fundamental dualizing complex in this case. 
}}
\end{rem}
\noindent{\em{Proof of \Tref{thm:gorenstein}.}} First note that since $\X$ carries a c-dualizing complex, it is necessarily of finite Krull dimension. 
Clearly (via a judicious use of $\BL_\X$ and $\R\iGp{\X}$) statement (b) is equivalent
to statement (a). We prove (a). To that end, suppose $\eG$ is
c-Gorenstein with respect to $\De$. Let $\eD$ be any c-dualizing complex whose associated
codimension function is $\De$.\footnote{By our hypothesis, such a $\eD$ exists, since any translation
of a c-dualizing complex is again c-dualizing, and any two codimension functions are translates of
each other.} Now $\eD$ is c-Gorenstein with respect to $\De$, whence by 
\Pref{prop:gorenstein}\,(a), $\R\sHomb(\eG,\eD)$
is $\D(\X)$ isomorphic to a locally free sheaf ${\mathcal{W}}$ of finite rank. Let ${\mathcal{V}}$
be the $\co_\X$-module dual (=vector bundle dual in this case) of ${\mathcal{W}}$. We then have
\[\eG\cong \R\sHomb({\mathcal{W}},\,\eD) \cong \eD\otimes{\mathcal{V}}.\]
The converse is obvious.  Indeed, if $\eG\cong \eD\otimes{\mathcal{V}}$
as in the statement, then $\eG$ is c-Gorenstein, since a c-dualizing complex is obviously
c-Gorenstein with respect to its associated codimension function.\qed

The proof of \Pref{prop:gorenstein} is incomplete since we need to prove \Lref{lem:kappa*}. So as before, assume $x\in\X$ is a closed point and $\kappa\colon {\widehat{\X}}\to \X$ the completion of
$\X$ along $\{x\}$.

The adjoint relation between $\kappa^*$ and $\kappa_*$ (see \cite[p.\,147,\,Prop.\,6.7]{sp})
gives us a bifunctorial map 
\stepcounter{thm}
\begin{equation*}\label{iso:kappa*}\tag{\thethm}
\kappa^*\R\sHomb_{\co_\X}(\eG_1,\,\eG_2) \longrightarrow 
\R\sHomb_{\co_{\widehat{\X}}}(\kappa^*\eG_1,\,\kappa^*\eG_2) \qquad (\eG_1,\eG_2\in\D(\X))
\end{equation*}
induced by the composite
\[ \R\sHomb_{\co_\X}(\eG_1,\,\eG_2) \longrightarrow 
\R\sHomb_{\co_\X}(\eG_1,\,\kappa_*\kappa^*\eG_2)
\xleftarrow{\Iso} \kappa_*\R\sHomb_{\co_{\widehat{\X}}}(\kappa^*\eG_1,\,\kappa^*\eG_2).\]

Suppose that $\eG_1$ and $\eG_2$ are both in $\Dc(\X)$ and one the following holds 
(cf.  \cite[p.\,8,\,(1),(2) and (3)]{ajl}):
\begin{enumerate}
\item $\eG_1\in\Dc^{-}(\X)$ and $\eG_2\in\Dc^{+}(\X)$; or
\item $\eG_2$ has finite injective dimension (i.e., $\eG_2$ is $\D(\X)$-isomorphic to a bounded
$\A(\X)$-injective complex; or
\item $\eG_1$ has local resolutions by finite locally free $\co_\X$-modules.
\end{enumerate}

Then using a ``way-out"  argument as in \cite[p.\,8]{ajl}---after replacing $\X$ by an affine 
neighborhood $\Spf{(A,I)}$ of our closed point if necessary---one sees  that \eqref{iso:kappa*}
is an isomorphism (one can reduce to the case where $\eG_1$ is bounded-above complex of
locally free $\co_\X$-modules and $\eG_2$ is a single coherent $\co_\X$-module).  In the event
$\eG_1$ and $\eG_2$ are c-Gorenstein \Lref{lem:kappa*} implies that condition (1) above is satisfied.  

\begin{lem}\label{lem:kappa*} Let $\X$ be of finite Krull dimension and $\De$ a codimension function
on it. 
\begin{enumerate}
\item[(a)] If $\eG\in\Dc(\X)$ is such that $\R\iGp{\X}\eG$ is Cohen-Macaulay on $(\X,\,\De)$ then
 $\eG\in\Dc^b(\X)$. 
 \item[(b)] If $\eG_1, \eG_2\in \Dc(\X)$ are such that $\R\iGp{\X}\eG_i$ is $\De$ Cohen-Macaulay for $i=1,2$,  
then \eqref{iso:kappa*} is an isomorphism for every closed point $x\in \X$.
\end{enumerate}
\end{lem}

\proof Suppose $\eF\set \R\iGp{\X}\eG$ is $\De$-Cohen-Macaulay. 
Then $\eF$ is bounded since
it is $\De$-Cohen-Macaulay and $\X$ has finite Krull dimension. Since $\X$ is noetherian, it is 
quasi-compact. By \cite[p.30,\,Lemma\,(4.3)]{ajl}  it follows that $\BL_{\X}\eF$ is bounded. But
$\eG\in\Dc(\X)$, whence $\BL_{\X}\eF \cong \eG$. This proves (a). Part (b) follows from the argument
given before the statement of the Lemma.
\qed

\subsection{Azumaya algebras}\label{ss:azumaya} We are interested in understanding the
algebra of endomorphisms of a t-Gorenstein complex (in $\Dcs$) on a formal schemes. We begin
a result that is well known over ordinary schemes (cf. \cite[p.\,57,\,Thm.\,5.1]{brauer}
 and \cite[p.\,163,\,Prop.\,6.11.1]{caen}) and lends itself to an easy generalization over formal
 schemes. Recall from \cite[p.124,\,10.3]{lns} that a map of
 formal schemes $f\colon\X\to \Y$ is {\em{\'etale}} if it is smooth of relative dimension zero. It is
 {\em{adic}} if $\I\co_\Y$ is an ideal of definition of $\X$ for some (and hence every)
 ideal of definition $\I\subset \co_\Y$ of $\Y$.
 
 \begin{thm}\label{thm:brauer5.1} Let $\X$ be a formal scheme,
 $\eA$ a coherent $\co_\X$-algebra. The following conditions are equivalent:
 \begin{enumerate}
 \item[(i)] $\eA$ is a locally free as a $\co_\X$-module, and, for every closed point $x\in\X$,
 $\eA(x)\set \eA\otimes k(x)$ is a central simple algebra, i.e., $\eA(x)\otimes_{k(x)}{\bar{k}(x)}$
 is isomorphic to a matrix algebra over the algebraic closure ${\bar{k}(x)}$ of $k(x)$.
 \item[(ii)] $\eA$ is a locally free $\co_\X$-module, and the natural homomorphism
 \[\eA\otimes_{\co_\X}\eA^\smcirc \to {\eE}nd_{\co_\X}(\eA)\]
 is an isomorphism, where
 $\eA^\smcirc$ is the opposite algebra of $\eA$.
 \item[(iii)] For every closed point $x\in \X$, there exists an integer $r\ge 1$, an open neighborhood
 $U$ of $x$, and a {\em{finite {\'e}tale surjective}} map $p\colon U'\to U$, such that $\eA_{U'}\set
 p^*\eA$ is isomorphic to the $r\times r$ matrix algebra $M_r(\co_{U'})$ over $\co_{U'}$.
 \item[(iv)] The same as (iii), with $p\colon U'\to U$ merely surjective and \'etale adic.
 \end{enumerate}
 \end{thm}
 
 \proof A couple of comments are in order. In (iii), the requirement that $U'\to U$ is finite
 forces the map to be adic. We have restricted ourselves to closed points in (i), (iii) and (iv)
 whereas \cite[thm.\,5.1]{brauer} has no such restriction, but very obviously for ordinary schemes,
 our formulation agrees with the classical formulation(s). The statements are local and therefore
 we may assume $\X=\Spf{(R,I)}$ where $(R,I)$ is an adic noetherian ring. Let $X=\Spec{(R)}$ and
 let $\kappa\colon \X\to X$ be the completion map. Suppose $A=\Gamma(\X,\eA)$. Then $A$
 is a finitely generated $R$-algebra and and defines an $\co_X$-algebra ${\widetilde{A}}$ (which
 is the same as the quasi-coherator  of $\kappa_*\eA$ 
 (cf. \cite[p.\,187,\,Lemme\,3.2]{illusie} and \cite[p.\,31, \S\,3.1]{dfs})). Since (i)--(iv) are equivalent
 conditions for ${\widetilde{A}}$ (i.e., with $(X,\,{\widetilde{A}})$ replacing $(\X,\,\eA)$), one can deduce
 the same for $\eA$. We leave the details to the reader. We point to \cite[\S\S\,2.1,\,pp.\,144-145]{dcc}
 for the relationship between the local rings of $X$ and $\X$, and we point out that since $(A,I)$
 is adic (i.e., $A$ is complete with respect $I$), $I$ is in the Jacobson radical of $A$.
 \qed
 
 \begin{rem}\label{rem:topos} {\em{See \cite[\S\,2,\,p. 51 and Remarque\,5.12, p. 60]{brauer}
 for generalizations to locally ringed toposes.}} 
 \end{rem}
 
 \begin{defi}\label{def:azumaya} Let $\X$ be a formal scheme and $\eA$ a sheaf of $\co_\X$-algebras on $\X$. We
 say $\eA$ is an {\em Azumaya algebra} on $\X$ if it satisfies any of the equivalent conditions
 of \Tref{thm:brauer5.1}.  An Azumaya algebra $\eA$ is said to be {\em split} if it is isomorphic 
 (as an $\co_\X$-algebra)
 to the algebra of endomorphisms ${\eE}nd_{\co_\X}({\mathcal{V}})$ of a locally free $\co_\X$-module
 of finite rank ${\mathcal{V}}$.
 \end{defi}

Let $(\X,\De)\in \bbFc$ with $\X$ of finite Krull 
dimension. For a t-Gorenstein complex $\eF\in \Dcs(\X)$ on $(\X,\De)$, we have a sheaf of
$\co_\X$-algebras $\eA=\eA(\eF)$ given by
\stepcounter{thm}
\begin{equation*}\label{def:A-F}\tag{\thethm}
\eA(\eF)\set \sHom_{\C(\X)}(\Ed{\De}(\eF),\Ed{\De}(\eF)).
 \end{equation*}

\begin{rems}\label{rem:azumaya} {\em{In what follows $\eF\in\Dcs(\X)$ is t-Gorenstein and $\eG=\BL_\X\eF$.

1) Very clearly, if $\eF\in\Dcs(\X)$ is t-Gorenstein, then
\[\eA(\eF)= \eA(\Ed{\De}(\eF)).\]
We could therefore have restricted ourselves to to complexes in $\Cozst_{\De}(\X)$ which are 
t-Gorenstein. However, allowing ourselves all t-Gorenstein complexes in $\Dcs(\X)$ gives us
greater flexibility. In practice (i.e., in our proofs) we will almost always identify t-Gorenstein
complexes $\eF$ with $\Ed{\De}(\eF)$ via the Suominen isomorphism $\Ssf\colon \eF \iso \Ed{\De}(\eF)$ \cite[p.42,\,Corollary 3.3.2]{lns}.

2) In $\D(\X)$, $\eA(\eF)$ can be identified with $\R\sHomb(\eF,\,\eF)$. In fact we have a sequence
of isomorphisms:
\begin{align*}
\eA(\eF) & \,\,\set  \,\,\sHom_{\C(\X)}(\Ed{\De}(\eF),\,\Ed{\De}(\eF)) \\
& \iso H^0(\R\sHomb(\eF,\,\eF))\\
& \iso \R\sHomb(\eF,\,\eF)
\end{align*}
The first isomorphism is deduced  by noting that  any $\C(\X)$-endomorphism of 
$\Ed{\De}(\eF)$ is equivalent to a
 $\D(\X)$-endomorphism of $\eF$ using Suominen's equivalence
of categories between Cohen-Macualay complexes and Cousin complexes (see
\cite[p.\,42,\, 3.3.1 and 3.3.2]{lns}), i.e., we have $\Hom_{\C(\X)}(\Ed{\De}(\eF),\Ed{\De}(\eF))$
is isomorphic to $\Hom_{\D(\X)}(\eF,\eF) = H^0(\R\Homb(\eF,\,\eF))$ which upon sheafifying gives us the
required isomorphism. The second isomorphism follows  from \Pref{prop:gorenstein}(b).

3) $\eA(\eF)$ is therefore a locally free $\co_\X$-module of finite rank by Prop.\,\ref{prop:gorenstein}(b).

4) By \Rref{rem:lam-gam} we have a $\D(\X)$ isomorphism
\[\eA(\eF) \cong \R\sHomb(\eG,\,\eG).\]

5) Let $x\in\X$ be a closed point, $\kappa=\kappa_x\colon {\widehat{\X}}\to \X$ the completion
of $\X$ along $\{x\}$, $\fm_x$ the ideal sheaf of $x$,
$A=\Gamma({\widehat{\X}},\,\co_{\widehat{\X}})={\widehat{\co}}_{\X,x}$,
$\fm=\fm_A$ the maximal ideal of $A$ (note ${\widehat{\X}}=\Spf{(A,\fm)}$), and $\De'$ the
codimension function on ${\widehat{\X}}$ which gives the unique point of ${\widehat{\X}}$
the value $\De(x)$. Set
\[{\widehat{\eF}} = \kappa^*\R\iG{\fm_x}\eF. \]
As in 1) above, we identify $\eF$ with $\Ed{\De}(\eF)$. Then ${\widehat{\eF}}$ is the sheafification
of the complex of $A$-modules 
$\Gamma_{\fm_x}\eF$ (= $\Hr_{\fm_x}^{\De(x)}(\eF)[-\De(x)]$).\footnote{Note 
that there is a difference between $\iG{\fm_x}$ and $\Gamma_{\fm_x}$. More precisely,
$ \Gamma_{\fm_x}\eF =\Gamma(\X,\iG{\fm_x}\eF)$} This shows that ${\widehat{\eF}}$ is
t-Gorenstein with respect to $\De'$ and it is clearly in $\Dcs({\widehat{\X}})$ since it can be
identified in $\D({\widehat{\X}})$ with $\R\iGp{{\widehat{\X}}}\kappa^*\eG$. We are now in a 
position to define a natural map of $\co_{\widehat{\X}}$-algebras
\stepcounter{thm}
\begin{equation*}\label{map:kappaA}\tag{\thethm}
\kappa^*\eA(\eF) \to \eA({\widehat{\eF}}).
\end{equation*}
Indeed, with $\eF=\Ed{\De}(\eF)$, any map $\eF\to \eF$ in $\C(\X)$ gives rise to a map of complexes
$\iGp{\fm_x}\eF\to \iGp{\fm_x}\eF$ in $\C({\widehat{\X}})$, and this correspondence 
is compatible
with Zariski localizations. In other words, we have a map of $\co_\X$-algebras
$\eA(\eF) \to \kappa_*\eA({\widehat{\eF}})$, which defines \eqref{map:kappaA} by the adjointness
of the pair $(\kappa^*,\,\kappa_*)$.
}}
\end{rems}

\begin{lem}\label{lem:azumaya} With the notations and hypotheses of \Rsref{rem:azumaya}, the map
$\kappa^*\eA(\eF) \to \eA({\widehat{\eF}})$ of \eqref{map:kappaA} is an isomorphism of 
$\co_{\widehat{\X}}$-algebras.
\end{lem}

\proof Let $\eG=\BL_{\X}\eF$. It is enough to show that \eqref{map:kappaA} is an isomorphism
in $\D({\widehat{\X}})$. To that end note that $\BL_{\widehat{\X}}{\widehat{\eF}}=\kappa^*\eG$.
By \Lref{lem:kappa*}, the natural map \eqref{iso:kappa*}
 is an isomorphism for $\eG_1=\eG_2=\eG$, i.e. we have a natural $\D(\X)$ isomorphism
 \[\kappa^*\R\sHomb(\eG,\,\eG) \iso \R\sHomb(\kappa^*\eG,\,\kappa^*\eG).\]
 The assertion follows from 2) and 4) of \Rsref{rem:azumaya} applied to $\eF$, $\eG$, ${\widehat{\eF}}$
 and $\kappa^*\eG$ (=$\BL_{\widehat{\X}}{\widehat{\eF}})$.
\qed

 \begin{prop}\label{prop:azumaya} Let $(\X,\De)\in \bbFc$ with $\X$ of finite Krull dimension,
 and suppose $\eF\in\Dcs(\X)$ is t-Gorenstein. Then the sheaf of
$\co_\X$-algebras $\eA(\eF)$ of \eqref{def:A-F} is an $\co_\X$-Azumaya algebra.
 \end{prop}
 
 \proof For a closed point $x\in \X$ we have, thanks to \Lref{lem:azumaya}, a $k(x)$ algebra
 isomorphism 
 \[\eA\otimes\otimes k(x) \iso \eA({\widehat{\eF}}) \otimes k(x),\]
 the tensor products being over appropriate structure sheaves. By \Tref{thm:brauer5.1}, we are
 therefore reduced to the case $\X=\Spf{(A,\fm)}$, $(A,\fm)$ being a complete local ring, and
 $\eF=\oplus_{i=1}^r\eE(x)[-\De(x)]$, a translate of a direct sum of a finite number of copies
 of the $\co_\X$-injective hull of $k(x)$. Since $\eE{nd}_{\co_\X}(\eE(x))$ is canonically
 isomorphic to $\co_\X$, it follows that $\eA(\eF)$ is isomorphic, as an $\co_\X$-algebra,
 to the matrix algebra $M_r(\co_\X)$. It follows that $\eA(\eF)\otimes k(x)$ is a matrix
 algebra over $k(x)$. The assertion follows from the definition of Azumaya algebras in
 \ref{def:azumaya} (see also \Tref{thm:brauer5.1}\,(i)).
 \qed
 
 \smallskip

 One might guess (based on results in \cite{ffgr}) that $\X$ has a
 c-dualizing complex if (and clearly only if) it has a t-Gorenstein complex $\eF\in\Dcs$ such that
 $\eA(\eF)$ is split. This may be overly optimistic. What can be shown is that in the event $\eA(\eF)$
 is split, $\X$ can be covered by Zariski open subschemes each of which admits a c-dualizing complex.
 However, if $\eA$ is split in a strong sense, i.e., when $\eA$ is isomorphic to a matrix algebra
 $M_r(\co_\X)$, then $\X$ does have a c-dualizing complex. 
%(see \Tref{thm:azumaya1}).  
  
 \begin{thm}\label{thm:azumaya1} Let $(\X,\,\De)\in \bbFc$ with $\X$ of finite Krull dimension.
 \begin{enumerate} 
 \item[(a)] Suppose $\X$ admits a c-Gorenstein complex $\eG$ with respect $\De$ such that 
 $\eA(\R\iGp{\X}\eG)$ is isomorphic as an $\co_\X$-algebra to $M_r(\co_\X)$ for some positive 
 integer $r$. Then $\X$ admits a c-dualizing complex.
 \item[(b)] Suppose $\X$ admits a t-Gorenstein complex $\eF\in\Dcs(\X)$  with respect to $\De$
 such that  $\eA(\eF)$ is isomorphic as an $\co_\X$-algebra to $M_r(\co_\X)$ for some positive
  integer $r$. Then $\X$ admits a c-dualizing complex.
 \end{enumerate}
 \end{thm}
 
 \proof Clearly (a) and (b) are equivalent. We will prove (b). To that end, we make, without loss
 of generality, the simplifying
 assumption that $\eF=\Ed{\De}(\eF)$. Now $M_r(\co_\X)$ has idempotent sections 
 $N_1,\dots, N_r$ such that $N_1+\dots+N_r=1$, $N_i^2=N_i$, $i=1,\dots,r$, and $N_iN_j=0$ for 
 $i\neq j$ and $1\le i,j \le r$. These give
 $r$ idempotent endomorphisms in $\C(\X)$ 
 \[\varphi_i \colon \eF \to \eF \qquad (i=1,\dots r) \]
 such that $\varphi_1+\dots+\varphi_r=1\in\eA(\eF)$,  $\varphi_i^2=\varphi_i$ for $i=1,\dots,r$, 
 and $\varphi_i\varphi_j=0$ for $i\neq j$ and
 $1\le i,j \le r$. These splitting idempotents break up $\eF$ into direct summands $\eF_i$,
 $i=1,\dots, r$. In fact $\eF_i$ is the image of $\eF$ under $\varphi_i$ as well as the
  kernel of the sum of the $\varphi_j$ distinct from $\varphi_i$. Direct summands of
  injectives being injective, and subcomplexes of Cousin complexes being Cousin,
  one sees that the $\eF_i$ are t-Gorenstein and for and $x\in\X$,
 $\eF(x)$ consists of only one copy of the injective hull of $k(x)$. We claim that $\eF_i\in\Dcs$
 for every $i$. But  $\BL_\X(\eF_i)$ can be identified with a direct summand of $\BL_\X(\eF)$,
 whence the former has coherent cohomology (since the latter has). This proves that
 $\eF_i\in \Dcs$ for $i=1,\dots, r$. 
 
 Using the matrices $M_{ij}$ which have $1$ in the $(i,j)$-th spot and zero elsewhere, we get maps
 $\varphi_{ij}\colon\eF\to \eF$.  Restricting $\varphi_{ij}$ to $\eF_j$ and projecting to $\eF_i$, we get
 maps $\psi_{ij}\colon \eF_j\to \eF_i$. It is easy to check that these are isomorphisms,
 and that the sheaf of endomorphisms, $\sHom_{\C(\X)}(\eF_i,\eF_j)$ is isomorphic
 to $\co_\X$ for every $1\le i,j \le r$. This proves that the
 $\eF_i$ are residual. Since they are in $\Dcs(\X)$, it follows that the $\BL_\X\eF_i$ are
 c-dualizing on $\X$.
 \qed
 
\smallskip

\Tref{thm:brauer5.1}(iii), \Tref{prop:azumaya} and \Tref{thm:azumaya1} gives immediately part (a) 
of the following theorem generalizing \cite[p.\,209,\,Cor.\,(4.8)]{ffgr}. Part (b) follows from \cite[p.\,109,\,Prop.\,9.3.5]{lns} and
part (a), using the fact that $f^*\Ed{\De}\cong \Ed{\De'}f^*$ for an \'etale adic map $f:(\X',\De')\to (\X,\De)$.

\begin{thm}\label{thm:azumaya2} Suppose $(\X,\,\De)$ has a c-Gorenstein complex which non-exact
on every connected component of $\X$
(or equivalently a t-Gorenstein complex in $\Dcs(\X)$, which is non-exact on every connected 
component of $\X$) and suppose $\X$ is of finite Krull dimension.
\begin{enumerate}
\item[(a)] For every closed point $x\in\X$
we can find a Zariski open neighborhood $U$ of $x$ and a finite \'etale map $U'\to U$ such
that $U'$ has a c-dualizing complex.
\item[(b)] We have  $\Ed{\De}\Dcs(\X)\subset \Dcs(\X)$. In particular, 
if $\X$ is an ordinary scheme, and $\eF$
has coherent cohomology, then $\Ed{\De}\eF$ has coherent cohomology.
\end{enumerate}
\end{thm}

\begin{rems} \label{rem:necessity}

{\em{1) Examining the proof of \Tref{thm:azumaya1} we notice that for a closed point
$x\in\X$, $\eA(\eF)\otimes k(x)$ is
already a split Azumaya algebra over $k(x)$ (i.e., a split central simple algebra over $k(x)$). This
means that the finite \'etale map $U'\to U$ ($U$ a Zariski neighborhood of $x$) can be chosen
so that $k(x')=k(x)$ for every point $x'\in U'$ lying over $x$. (See \cite[pp.\,58--59,\,5.4--5.8]{brauer}.)

2) Our main results concerning c-Gorenstein complexes (especially
\Tref{thm:azumaya1} and \Tref{thm:azumaya2}) rely crucially on the technical result that when $\X$ is
of finite Krull dimension \eqref{iso:kappa*} is an isomorphism for every closed point of $\X$
(see \Lref{lem:kappa*}). On the other hand, with the benefit of hindsight, if $\X$ is a scheme admitting
a c-Gorenstein complex, and is such that \eqref{iso:kappa*} is an isomorphism for every closed
point $x\in\X$, then the conclusion of \Tref{thm:azumaya2} holds. But any scheme containing
a c-dualizing complex is of finite Krull dimension, and a little thought shows that this forces the
quasi-compact scheme $\X$ to have finite Krull dimension.}}
\end{rems}

\Tref{thm:azumaya2} gives us the following (cf. again \cite[Cor.\,(4.8)]{ffgr}):

\begin{cor}\label{cor:hensel} Suppose $A$ is a Hensel local ring, complete with respect to an ideal
$I\subset A$, and let $\X=\Spf{(A,I)}$. The following are equivalent:
\begin{enumerate}
\item[(a)] $\X$ has a codimension function $\De$ and $(\X,\,\De)$ has a c-Gorenstein complex.
\item[(b)] $\X$ has a c-dualizing complex.
\item[(c)] $\X$ has a codimension function $\De$ and $(\X,\De)$ has a t-Gorenstein complex in 
 $\Dcs(\X)$.
 \item[(d)] $\X$ has a t-dualizing complex in $\Dcs(\X)$.
\end{enumerate}
\end{cor}

\proof Observe that if $\X$ has a c-dualizing complex then it necessarily
has a codimension function.  For the rest, first note that (a) and (c) are equivalent, as are (b) and (d).
Now let $k=A/\fm_{{}_A}$ be the residue field of $A$.  Now since $A$
is Hensel any \'etale neighborhood of $\Spec{(A)}$ having a $k$-valued point over the closed point
of $\Spec{(A)}$ admits a section through that point (whence so does any adic
\'etale  neighborhood of $\Spf{(A,I)}$ having a $k$-valued point). The equivalence of (a) and (b) is now
immediate, in view of the \Tref{thm:azumaya2} and 1) of \Rsref{rem:necessity}.
\qed

\begin{rems}\label{rem:comm-alg} {\em{We wish to re-interpret our results in commutative algebraic
terms. If $(A,I)$ is an adic  ring (always noetherian and of finite Krull dimension in these
remarks) then we will move from complexes
of $A$-modules and sheaves on $\Spf{(A,I)}$  in the usual fashion. The terminology used
in this remark is self-explanatory. For example, and if we say that a complex
$M^\bullet$ of $A$-modules lies in $\Dcs(A,I)$, then what is meant is that the corresponding
complex on $\Spf{(A,I)}$ is in $\Dcs(\X)$. In the event an ideal of definition $I$ in $A$ is not specified,
we take $I=0$, so that $\Spf{(A,I)}=\Spec{A}$. In these remarks, all rings and algebras occurring are 
commutative.

1) Suppose $(A,I)$ is an adic  ring, $\De$ a
codimension function on $(A,I)$, $M^\bullet$ a complex whose cohomologies are finite $A$-modules, and $M^\bullet \to J^\bullet$  an injective resolution of $M^\bullet$.
If the complex $\Gamma_IJ^\bullet$ is quasi-isomorphic to its Cousin complex  
$E^\bullet= \Ed{\De}(\R\Gamma_IM^\bullet)$ with respect to  $\De$ and $E^\bullet$ consists of
$A$-injectives, then  for every maximal ideal $\fm$ of $A$, there is an element $f\in A\setminus \fm$
and  a {\em finite \'etale} $A_f$-algebra $B$  such that $B$ possesses
a dualizing complex $D^\bullet_B$ (in the usual sense).  Moreover, in $\D(B)$ ,we have
\[ M^\bullet\otimes_AB \cong D^\bullet_B\otimes P\]
where $P$ is a finitely generated {\em projective} $B$-module.

2) With the above hypotheses, since $B$ possesses a dualizing complex, $\Spec{B}$ has a 
codimension function, and it is not hard to see that whence so has $\Spec{A}$. We add that 
$D^\bullet_B$ is
a c-dualizing complex on $(B,IB)$ and $\R\Gamma_ID_B^\bullet$ is t-dualizing on $(B,IB)$ and is
in $\Dcs(B,IB)$.

3) Suppose $A$ is a ring and $I\subset A$,  with $A$ not necessarily complete with respect to $I$. 
Suppose that $A/I$ admits a codimension function $\De$ and
that we have a complex $M^\bullet$ with finitely generated cohomology modules
satisfying the hypotheses in 1) with respect to $I$ and $\De$ (note that $\Gamma_IJ^\bullet$ being
supported on $V(A/I)$,  we can talk about its Cousin with respect to $\De$). 
Let $\widehat{A}$ be the completion of $A$ with respect to $I$. Then, for every maximal
ideal $\fm$ of ${\widehat{A}}$, there is an element $f\in A\setminus\fm$ and a  finite \'etale 
$A_f$-algebra $B$ such that  $B$ possesses a dualizing complex $D_B^\bullet$ and in $\D(B)$ 
we have
\[M^\bullet\otimes_AB\cong D^\bullet_B\otimes_B P\]
where $P$ is a finitely generated projective $B$-module. 

4) The case of $I=\fm_A(=\fm)$, a maximal ideal of $A$,
is then quite interesting. With $M^\bullet$, $J^\bullet$ as above, the previous statement reduces
to the following: if $\Hr^i(\Gamma_\fm J^\bullet)=0$
for $i\neq \De(\fm)$ and $\Hr^{\De(\fm)}(\Gamma_\fm J^\bullet)$ is an injective $A$-module, then
$M^\bullet\otimes_A{\widehat{A}}$ is isomorphic in $\D({\widehat{A}})$ to a direct sum of
a finite number of copies of a dualizing complex $D^\bullet$ on ${\widehat{A}}$.
}}
\end{rems}

\subsection{Gorenstein modules}\label{ss:sharp} Let $A$ be a noetherian ring of finite Krull dimension. 
Recall, from \cite[p.\,123]{sh1}, that a non-zero finite module $M$ is called Gorenstein if its Cousin 
complex $E(M)$, with respect to the height filtration of $M$, is an injective resolution (necessarily 
minimal) of $M$, i.e., $E(M)$ consists of injectives, and $H^0(E(M))=M$. In such a case, according to a 
result of Foxby \cite{foxby}, $\Hom_A(M,M)$ is a projective $A$-module. Moreover,
 $A$ is Cohen-Macaulay \cite[p.\,126,\,Cor.\,3.9]{sh1}. Further ${\text{Supp}}_AM=\Spec{(A)}$.
We end the paper with a discussion connecting
our results with these (well known) results on Gorenstein modules. The claim is not that our proofs 
are simpler, but we hope they are illuminating, as perhaps any set of proofs involving a different point
of view will tend to be. Since we work on the ordinary scheme $\Spec{(A)}$, the notions of
c-Gorenstein complexes and t-Gorenstein complexes in $\Dcs(A)$ coincide and we simply use
Gorenstein (without adornments) while referring to such complexes.

Fix $\fp\in\Spec{(A)}$. With $i={\text{ht}}_M(\fp)$, set $E(\fp)\set \Hr^i_\fp(M_\fp)$, i.e., $E(\fp)$ is
the summand of the total module of $E(M)_\fp$ which has $\fp$ as its associated prime.

One can show in the usual way, for example by using \cite[Lemma\,(3.1)]{bass}\footnote{The statement
is that for any finite $A$-module $N$, any non-negative integer $i$, and any immediate specialization 
$\fp\mapsto \fq$, the Bass number $\mu^i(\fp,N)\neq 0$ only if the Bass number 
$\mu^{i+1}(\fq,N)\neq 0$.}, that $\fp\mapsto{\text{ht}}_M(\fp)$ is a codimension function on 
${\text{Supp}}_A(M)$. In fact, as we shall shortly see, ${\text{Supp}}_A(M)=\Spec{(A)}$, whence
${\text{ht}}_M$ is a codimension function and $E(M)$ is a Gorenstein complex.

Arguing as we did in the proof of  \Pref{prop:gorenstein} (bearing in mind that $A$ has finite
Krull dimension) we see that for $\fm\in {\text{Supp}}_A(M)$ a maximal ideal, and ${\widehat{A}}$
the $\fm$-adic completion of $A$, we have
\[\Homb_A(E(M),\,E(M))\otimes_A{\widehat{A}} \iso \Homb_{\widehat{A}}(E(\fm),\,E(\fm)). \]
The right side has no cohomology in positive degrees therefore neither does the left side.
Taking the zeroth cohomology we get
\[\Hom_A(M,\,M)\otimes_A{\widehat{A}} \iso \Hom_{{\widehat{A}}}(E(\fm),\,E(\fm)). \]
But if $I=0:_AM$, then the left side is killed by  ${\widehat{I}}\set I{\widehat{A}}$. 
The right side is a free ${\widehat{A}}$-module. It follows that ${\widehat{I}}=0$, whence $I=0$. 
In other words,  ${\text{Supp}}_AM=\Spec{(A)}$. In particular ${\text{ht}}_M$ is a codimension function.

In order to show that $A$ is Cohen-Macaulay, it is enough to show this is so under the
assumption that $A$ has a dualizing complex. Indeed, $E(M)$ is a Gorenstein complex,
therefore, and therefore, according to  \Tref{thm:azumaya2}, we can find $f_1,\dots,f_n\in A$,
with $f_1+\dots+f_n=1$, and finite \'etale algebras $A_{f_i}\to B_i$, $i=1,\dots,n$, such that
$B_i$ possess dualizing complexes. Moreover, with $M_i=M\otimes_AB_i$, it is easy to see
that $E(M_i)=E(M)\otimes_AB_i$ is an injective resolution of $M_i$, whence $M_i$ is 
Gorenstein as a $B_i$-module. If $B_i$ is Cohen-Macaulay, so is $A_{f_i}$, and if this is
so for every $i=1,\dots,n$, clearly $A$ is Cohen-Macaulay (since $f_1+\dots+f_n=1$). 

We therefore assume, without loss of generality, that $A$ possesses a dualizing complex
$D^\bullet$ (as before, following commutative algebra conventions, $D^\bullet$ is residual). 
Further, by translating $D^\bullet$, we may assume that the associated codimension function
is ${\text{ht}}_M$. By \Tref{thm:gorenstein}, there is a projective $A$-module $P$ and a $\C(A)$-isomorphism:
\[E(M) \iso D^\bullet\otimes_AP.\]
Since $E(M)$ resolves $M$ and $P$ is projective, it follows that $H^i(D^\bullet)=0$ for $i\neq 0$
and the natural map $K=H^0(D^\bullet) \to D^\bullet$ is a resolution. It follows that $A$ is
Cohen-Macaulay.

If the reader wishes to test her/his hold on the techniques of \cite{dfs}, especially as used above, then
the following exercise is recommended:

{\emph{Exercise:}} Let $(A,\fm)$ be a complete local ring, and $M\neq 0$ a finite $A$-module such that 
$\Hr^j_\fm(M)$
is an injective $A$-module for some $j\ge 0$ and $\Hr^i_\fm(M)=0$ for $i\neq j$. Show that $M$ is
a Gorenstein $A$-module. [{\it{Hint:}} Show that
$\R\Gamma_\fm(M)$ is t-Gorenstein on the adic ring $(A,\fm)$ for a suitable codimension function
on $\Spf{(A,\fm)}$. Conclude, using $\BL_{{}_\X}\iGp{\X}(\eG)\cong \eG$ for $\eG\in\Dc(\X)$ ($\X$ any
formal scheme), that $M$ is c-Gorenstein on $(A,\fm)$. Use \Tref{thm:gorenstein} to conclude that
$M$ is a Gorenstein module.]

\begin{ack} We thank Joe Lipman and Amnon Yekutieli for stimulating discussions.
Yekutieli drew the second author's attention to the re-interpretation of
\cite[p.\,182,\,Theorem\,7.2.2]{dcc} using the correspondence between coherent
sheaves and Cohen-Macaulay complexes (with coherent cohomology). We would also like to
thank Rodney Sharp for providing us with appropriate references related to work done by him,
and by Dibaei and Tousi.

\end{ack}

%\newpage

\bibliographystyle{plain}

\begin{thebibliography}{ABC}

\bibitem[Ao]{ayma} Y.\:Aoyama, {\em Some basic results on canonical modules}, J. Math. Kyoto Univ. 
{\bf 23} (1983), 85--94.

\bibitem[AJL1]{ajl} L.\:Alonso Tarr{\'{\i}}o, A.\:Jerem{\'{\i}}as L{\'o}pez
and J.\:Lipman, Local homology and cohomology on schemes,
{\it Ann.\:Scient.\:{\'E}c.\:Norm.\:Sup.\;}{\bf 30}\,(1997), 1--39. See also
{\it  Correction,} on page 879 of vol.\;2 of the {\it Collected Papers of
Joseph Lipman,} Queen's Papers in Pure and Applied Math., Vol.\;{\bf 117},
Queen's University, Kingston, Ontario, Canada, 2000.

\bibitem[AJL2]{dfs}
\bysame, Duality 
and flat base change on formal schemes, 
{\emph{Contemporary Math.,}} Vol.\;{\bf 244}, Amer. Math. Soc., Providence, R.I.
(1999), 3--90

\bibitem[AJL3]{gm}
\bysame, Greenlees-May
Duality on formal schemes,
{\emph{Contemporary Math.,}} Vol.\;{\bf 244}, Amer. Math. Soc., Providence, R.I.
(1999), 93--112

\bibitem[C]{caen} S.\:Caenepeel, {\emph{Brauer Groups, Hopf Algebras and Galois Theory}},
Kluwer Academic Publishers, Dordrecht, 1998.


\bibitem[Db]{dib}
M.\,T.\,Dibaei, A study of Cousin complexes through the dualizing complexes,
{\emph{Comm. Algebra}} {\bf 33} (2005), no.\,1, 119--132.

\bibitem[DT]{dt}
\bysame, M.\,Tousi,
A generalization of the dualizing complex structure and its 
applications, \emph{J.\ Pure and Applied Algebra} {\bf 155} (2001), 17--28.

\bibitem[B]{bass} H.\:Bass, On the ubiquity of Gorenstein rings, {\emph{Math. Z.}} {\bf 82} (1963),
8--28.

%\bibitem[AJL3]{dfscorr}\bysame, Correction to the paper ``Duality
%and flat base change on formal schemes", {\it Proc. Amer. Math. Soc.,}
%Vol.\;{\bf 131}, 351--357. 
%
%\bibitem[AJL4]{dgm}\bysame, Greenlees-May duality on formal schemes,
%{\it Contemporary Math.,} Vol.\;{\bf 244}, Amer. Math. Soc., Providence,
%R.I. (1999), 93--112, 2002.
%
\bibitem[C]{conrad} B.\:Conrad, {\it Grothendieck Duality and Base Change},
Lecture Notes in Math., no.\,{\bf 1750}, Springer, New York, 2000.

%\bibitem[D1]{deligne} P.\:Deligne, Cohomology {\`a} support propre,
%et construction du foncteur $f^!$, appendix to R.\:Hartshorne's
%{\it Residues and Duality}, Lecture Notes in Math., no.\,{\bf 20},
%Springer--Verlag, Heidelberg, 1966.
%
%\bibitem[D2]{SGA4}\bysame, {\it Cohomology {\`a} supports propre}, SGA4 Tome\,3,
%Lecture Notes in Math., no.\,{\bf 305}, Springer--Verlag, New York, 1973.
%
%\bibitem[D3]{nagata3}
%\bysame, Deligne's notes on Nagata's Compactification. Notes by B.~Conrad.
%\newblock Unpublished.
%

\bibitem[F]{foxby} H.-B.\:Foxby, Gorenstein modules and related modules, {\emph{Math. Scand.}}, {\bf 31}
(1972), 367--384.

\bibitem[FFGR]{ffgr} R.\;Fossum, H-B.\;Foxby, P.\;Griffith and I.\;Reiten,
Minimal injective resolutions with applications to dualizing modules
and Gorenstein modules, {\emph{Publ.\,Math.\, IHES,}} {\bf 40} (1975),
193--215.

\bibitem[G]{brauer} A.\:Grothendieck, {\em {Groups de Brauer I, II, III}} in ``Dix Expos{\'e}s sur la
cohomologie des sch{\'e}mas", North Holland, Amsterdam (1968), 46--65, 66--87, 88--188.

%\bibitem[GD]{gd} A.\:Grothendieck and J.\:Dieudonn\'{e}, 
%{\it \'Elements de G\'{e}om\'{e}trie Alg\'{e}brique\/} {\bf I}, 
%Springer-Verlag, New York, 1971.
%
%\bibitem[EGA-III]{ega3} \bysame, {\it
%\'Elements de G\'{e}om\'{e}trie Alg\'{e}brique\/} {\bf III},
%Publications Math. IHES {\bf 11}, Paris, 1961.
%
%\bibitem[EGA-IV]{ega4} \bysame, {\it
%\'Elements de G\'{e}om\'{e}trie Alg\'{e}brique\/} {\bf IV},
%Publications Math. IHES {\bf 20}, Paris, 1964.
%

\bibitem[I]{illusie} L.\:Illusie, Existence de r\'esolutions globales, in {\em{Th\'eorie des Intersections
et Th\'eor\`eme de Riemann-Roch (SGA 6)}}, Lecture Notes in Math., no. {\bf 225}, Springer-Verlag,
New York, 1971, pp. 160--222.

\bibitem[Hrt]{RD} R.\:Hartshorne, {\it Residues and Duality,} Lecture
Notes in Math., no.\,{\bf 20}, Springer-Verlag, New York, 1966.

\bibitem[Kw]{kw}
T.\,Kawasaki, Finiteness of Cousin cohomologies, preprint,
\hfill\newline
{\tt <\ http://www.comp.metro-u.ac.jp/~kawasaki/articles.html\#preprint\ >}

\bibitem[LNS]{lns}
J.~Lipman, S.~Nayak, P.~Sastry,
\newblock Pseudofunctorial behavior of Cousin complexes on formal
schemes, {\emph{Contemporary Math.,}} Vol.\;{\bf 375}, Amer. Math. Soc.,
Providence, R.I. (2005), 3--133.

%\bibitem[Lu]{nagata2}
%W.~L{\"u}tkebohmert, On compactification of schemes,
%\newblock {\em Manuscripta Math.}, {\bf 80} (1993) 95---111. 
%
%\bibitem[M]{matsumura}
%H.~Matsumura,
%\newblock {\it Commutative ring theory,} Cambridge Univ.~press, Cambridge, 1986.
%
%\bibitem[N]{nagata}
%M.~Nagata,
%\newblock Imbedding an Abstract Variety in a Complete Variety,
%{\em J. Math. Kyoto Univ.}, {\bf 2} (1962) 1---10. 
%
\bibitem[Nay]{suresh}
S.~Nayak,
\newblock Pasting pseudofunctors, {\emph{Contemporary Math.,}} Vol.\;{\bf 375}, Amer. Math. Soc.,
Providence, R.I. (2005), 195--271.

%\bibitem[Hu1]{I-C} I-C.\:Huang, {\it Pseudofunctors on modules with zero
%dimensional support}, Memoirs, no.\,{\bf 548}, Amer. Math. Soc., 1995.
%
%\bibitem[Hu2]{IC2} \bysame, The Residue Theorem via an Explicit Construction
%of Traces, {\it Jour. Alg.}, {\bf 245} (2001), 310--354.
%
%\bibitem[S]{ast-208} P.\:Sastry, A pointwise criterion for dualizing pairs,
%Appendix to {\it An explicit construction of the Grothendieck residue
%complex} by A.\:Yekutieli, 117---126, Ast\'erisque {\bf 208} (1992).
%

\bibitem[S]{dcc}P.\,Sastry, Duality for Cousin Complexes, 
{\emph{Contemporary Math.,}} Vol,\;{\bf 375}, Amer. Math. Soc., Providence,
R.I. (2005), 137--192.

\bibitem[Sh1]{sh1}R.Y.\,Sharp, Gorenstein Modules, {\emph{Math.~Z.,}}
Vol,\;{\bf 115} (1970), 117--139.

\bibitem[Sh2]{sh2}\bysame, On Gorenstein modules over a complete
Cohen--Macaulay local ring, {\emph{Quart.~J.~Math.,}} (2), {\bf 22} (1971),
425--434.

\bibitem[Sh3]{sh3}\bysame, Cousin complex characterizations of two classes
of commutative noetherian rings, {\emph{J.~London Math.~Soc.,}} (2), {\bf{3}} 
(1971), 621--624.

\bibitem[Sh4]{sh4}\bysame, Finitely generated modules of finite injective
dimension over certain Cohen--Macaulay ring, {\emph{Proc.~London Math.~Soc.,}}
(3), {\bf 25} (1972), 303--328.

\bibitem[SS]{sh5}\bysame and P.\,Schenzel, Cousin complexes and Generalized 
Hughes complexes, {\emph{Proc.~London Math.~Soc.,}} (3), {\bf{68}} (1994),
499--517.

\bibitem[Sp]{sp}N.\,Spaltenstein, Resolutions of unbounded complexes, {\emph{Compositio 
Mathematica}} {\bf{65}} (1988), 121--124.

\bibitem[Su]{suominen}K.\:Suominen, Localization of sheaves and Cousin
complexes, {\it Acta mathematica}, {\bf 131} (1973), 1--10.

%\bibitem[Lp1]{ast-117} J.\:Lipman, {\it Dualizing sheaves, Differentials and Residues on
%Algebraic Varieties,}\hfill\\
%   Ast\'erisque {\bf 117} (1984).
%

%\bibitem[V]{f!} J. L. Verdier, Base change for twisted inverse image
%of coherent sheaves, in {\it Algebraic Geometry,} Oxford Univ.~press,
%1969, pp. 393--408.
%
%\bibitem[Y]{Ye} A. Yekutieli, Smooth formal embeddings and the residue complex,
%{\it Canadian J. Math.} {\bf 50}~(1998), \mbox{863--896.}

\bibitem[Y]{y}
A.\,Yekutieli,  Smooth formal embeddings and the residue complex, {\emph{Canadian J. Math}},
{\bf 50} (1998), no. 4, 863--896.

\bibitem[YZ]{yz}
A.\,Yekutieli,  J.\,J.\,Zhang,
Rigid dualizing complexes on schemes, preprint, math.AG/0405570.

\end{thebibliography}
%\printindex
\end{document}